\input amstex
\documentstyle{amsppt}
\input xypic
\magnification\magstep1
\vsize8.9 true in
\hsize6.5 true in

\document
\define\qq{{\Bbb Q}}
\define\zz{{\Bbb Z}}
\define\R{{\Bbb R}}
\define\F{{\Bbb F}}
\define\ff{\dot F/\dot F^2}
\define\<{\langle}
\define\ifff{{if and only if} }
\define\>{\rangle}
\define\lra{\longrightarrow}
\define\script{\Cal}
\define\wg{{\Cal G}_F}
\define\wga{{\Cal G}_{F_{1}}}
\define\wgb{{\Cal G}_{F_{2}}}
\define \Gal {\operatorname{Gal}}
\define \Br {\operatorname{Br}}
\define \ch {\operatorname{char}}
\define \vcd {\operatorname{vcd}}
\define \Span {\operatorname{Span}}
\define \Sper {\operatorname{Spec_{r}}}
\define \Ker {\operatorname{Ker}}

\def\today{\ifcase\month\or
      January\or February\or March\or April\or May\or June\or
      July\or August\or September\or October\or November\or December\fi
      \space\number\day, \number\year}

\NoBlackBoxes
\topmatter

\title\nofrills
Additive structure of  multiplicative subgroups of fields and Galois theory
\endtitle
\rightheadtext\nofrills{Additive structure of  multiplicative subgroups of
fields}
\leftheadtext\nofrills{L.Mah\'e, J. Min\'a\v c and T.L. Smith}

\author
Louis Mah\'e,
J\'an Min\'a\v c *$^\dagger$,  and Tara L. Smith $^\ddag$
\endauthor
\thanks{The two first authors were supported by the Volkswagen-Stiftung
(RIP-Program at
Oberwolfach)}
\endthanks
\thanks{* Research supported in part by the Natural
Sciences
and Engineering Research Council of Canada and by the special Dean of
Science Fund at the University of Western Ontario.}
\endthanks

\thanks{$^\dag$ Supported by the Mathematical Sciences Research
Institute, Berkeley}
\endthanks
\thanks{$^\ddag$ Research supported in part by the Taft Memorial Fund of the University of Cincinnati}
\endthanks
\address
IRMAR, Campus de Beaulieu,  F-35042 Rennes-Cedex France
\endaddress
\email
Louis.Mahe\@univ-rennes1.fr
\endemail

\address
Dept. of Mathematics, Univ. of Western Ontario, London,
Ontario N6A 5B7  Canada
\endaddress
\email
minac\@uwo.ca
\endemail

\address
Dept. of Mathematical Sciences, Univ. of
Cincinnati, Cincinnati, Ohio 45221-0025  U.S.A.
\endaddress
\email
tsmith\@math.uc.edu
\endemail

\subjclass Primary 11E81; Secondary 12D15 \endsubjclass

\abstract
One of the fundamental questions in current field theory, related
to Grothendieck's conjecture of birational anabelian geometry, is the
investigation
of the precise relationship between the Galois theory of fields and
the structure
of the fields themselves. In this paper we initiate the classification of
additive properties of multiplicative subgroups of fields containing all
squares, using pro-$2$-Galois groups of nilpotency class at most $2$, and
of exponent  at most $4$. This work extends some powerful methods and
techniques
from formally real fields to general fields of characteristic not $2$.
\endabstract

\toc
\head 1.  Introduction\endhead
\head 2.  Groups not appearing as  subgroups of $W$-groups\endhead
\head 3.  Galois groups and additive structures (1)\endhead
\head 4.  Maximal extensions, closures and examples \endhead
\head 5.  Cyclic subgroups of $W$-groups\endhead
\head 6.  Subgroups of $W$-groups generated by two elements\endhead
\head 7.  Classification of  rigid orderings\endhead
\head 8.  Construction of closures for  rigid orderings\endhead
\head 9.  Galois groups and additive structures (2)\endhead
\head 10.   Concluding remarks\endhead
\head {}  References\endhead
\endtoc
\endtopmatter

\head
\S 1.  Introduction
\endhead

Let $F$ be a field of characteristic not $2$ and $T$ be a
multiplicative subgroup of $\dot{F}=F\setminus\{0\}$ containing the squares.
By the additive structure of $T$, we mean a description of the
$T$-cosets forming $T+aT$.  The purpose of this article is to relate the
additive structure of such a group $T$, to some Galois pro-$2$-group
$H$ associated with $T$.  In the case when $T$ is a usual
ordering, the group $H$ is a group of order $2$.  In the general case,
$H$ is a pro-$2$-group of  nilpotency class  at most $2$, and of
exponent  at most $4$.  Therefore the structure of $H$ is relatively
simple, and this is one of the attractive features of this
investigation.

One of our main motivations is  to extend Artin-Schreier
theory to this general situation.  In  classical Artin-Schreier theory as modified by Becker, 
one studies euclidean closures and their relationship
with Galois theory~\cite{ArSch1, ArSch2, Be1}.  Recall that such a closure
is a maximal
$2$-extension of an ordered field to which the given ordering
extends. (See~\cite{Be1}.)

It came as a surprise to us that for a good number of isomorphism types
of groups $H$ as above, we
could provide a complete algebraic characterization of the multiplicative
subgroups of $\dot{F}/\dot{F}^2$ associated with $H$, entirely analogous to
the classical
algebraic description of orderings of fields. We thus obtain a
fascinating direct link between Galois theory and additive properties
of multiplicative subgroups of fields.

We obtain in particular a Galois-theoretic characterization of
rigidity conditions (Proposition~3.4 and
Proposition~3.5) and a full classification of rigid groups
$T$ (\S 7). We also know how to make closures (as defined below) with
respect to these rigid ``orderings'' (\S 8).

In \S 9 we refine the notion of $H$-orderings of fields.  We show that
under natural conditions, we can control the behaviour of the additive
structure of these orderings under quadratic extensions.  It is
worthwhile to point out that each finite Galois $2$-extension can be
obtained by successive quadratic extensions.  Therefore, it is
sufficient to investigate quadratic extensions.

We have in \S 2 a nice illustration of what a $W$-group can or cannot be.
Since the $W$-group of the field $F$, together with its level,
determines the Witt ring $W(F)$, it is clear that every result about
   the $W$-group of $F$ and its subgroups will provide information on
   $W(F)$.

   This fits together with one of the main ideas behind this work (see
   \S 10): obtaining new Local-Global Principles for quadratic forms, with
   respect to these new ``orderings''. This will be the subject of a
   subsequent article.

\bigskip
We now enter into more detail, fix some notation, and present a more
technical outline of
the structure of the paper.

\definition{Notation 1.1}
All fields in this paper are assumed to be of characteristic not $2$, with
any exceptions clearly pointed out.
Occasionally we denote a field extension $K/F$ as $F\lra K$. The
compositum of two fields $K$ and $L$ contained in a larger
field is denoted as $KL$.
Recall that the {\it level} of a field $F$ is
the smallest natural number $n>0$ such that $-1$ is a sum of $n$
squares in $F$ or $\infty$ if no such $n$ exists.
Given a field $F$, we denote by $F(\sqrt{\dot F})$ the compositum of
all quadratic extensions of $F$, and by $F^{(3)}$ the compositum of
all quadratic extensions of $F(\sqrt{\dot F})$ which are Galois over
$F$.  (The field $F(\sqrt{\dot F})$ was denoted by $F^{(2)}$ in
previous papers (e.g. \cite{MiSm2}), and this explains the notation
$F^{(3)}$.)  The W-group of the field $F$ is then defined as $\wg =
\Gal(F^{(3)}/F)$.
This W-group is the Galois-theoretic analogue of
the Witt ring, in that if two fields have isomorphic Witt rings, then
their W-groups are also isomorphic.  Conversely, if two fields have
isomorphic W-groups, then their Witt rings are also isomorphic
(provided that the fields have the same level when the quadratic form
$\< 1,1\>$ is universal over one of the fields (See \cite{MiSp2, Theorem 3.8})).

We denote by $\Phi(\wg)$ the Frattini subgroup of $\wg$.  The
Frattini subgroup is by definition the intersection of the maximal
proper subgroups
$H$ of $\wg$.
(This means that $H$ is a maximal subgroup of $\wg$ among the family
of all closed subgroups of $\wg$ not equal to
$\wg$. It is a basic fact in the theory of pro-$2$-groups that each
such subgroup of $\wg$ is a closed subgroup of $\wg$ of
index two.)
Notice that $\Gal(F^{(3)}/F(\sqrt{\dot F})) =
\Phi(\wg )$.  In the case of a pro-$2$-group $G$, the Frattini
subgroup is exactly the closure of the group generated by squares.
Observe that for each closed subgroup $H$ of $\wg$ we have 
$\Phi(H)\subseteq\Phi(\wg)\cap H$.  We say that a
closed subgroup $H \subseteq \wg$ satisfying $\Phi(H) = H \cap
\Phi(\wg)$ is an {\it essential} subgroup of $\wg$.  Two essential
subgroups $H_1, H_2$ are {\it equivalent} if $H_1\Phi(\wg) =
H_2\Phi(\wg)$.  In general, for a closed subgroup $H$ of $\wg$, we
have $H = \script E \times \prod_{i}(\zz/2\zz)_{i}$ where $\script E$ is
essential:
$\Phi (H) = \Phi (\script E)$ and $\Phi (\wg ) \cap H \cong \Phi
(\script E ) \times \prod_{i}(\zz/2\zz)_{i}$.  The equivalence class of
$\script E$
is that of $H$, and equivalent essential subgroups are always
isomorphic.  (See \cite{CrSm, Theorem 2.1}.)
\enddefinition

\remark{Remark}
For typographical reasons we are using two different notations for the
action of a Galois element $\sigma$ on a field element $x$: the
exponential notation $x^\sigma$ and the functional notation $\sigma(x)$.
This should not cause confusion, as in any given instance
the order in which the elements enter the products will be clear or is
irrelevant.
\endremark

We now give the field-theoretic interpretation of the notion of an
essential subgroup of $\wg$. Let $H$ be any closed subgroup of $\wg$
and let $L$ be the fixed field of $H$. Let  $N$ and $M$ be the
fixed fields of $\Phi(H)$ and $\Phi(\wg)\cap H$ respectively. Because
$\Phi(H)\subseteq\Phi(\wg)\cap H$, we see that $M\subseteq N$ and
equality holds for one of the inclusions if it holds for the other.
Finally observe that $M$ is the compositum of $F(\sqrt{\dot F})$ and
$L$, and that $N$ is the compositum of all quadratic extensions of $L$
contained in $F^{(3)}$. Summarizing the discussion above we obtain:

\proclaim{Proposition 1.2} Let $H$ be a closed subgroup of $\wg$ and
$L$ be the fixed field of $H$.  Then $H$ is an essential subgroup of
$\wg$ if and only if the maximal multiquadratic extension of $L$
contained in $F^{(3)}$ is equal to the compositum of $L$ and
$F(\sqrt{\dot F})$.
\endproclaim

   Kummer theory and Burnside's Basis Theorem allow us
to prove the following:

\proclaim{Proposition 1.3} For $H$ a closed subgroup of $\wg$, the
assignment $$H\mapsto u(H)=P_H := \{ a \in \dot F\,\mid\,
(\sqrt{a})^{\sigma} = \sqrt{a}, \quad\forall \sigma \in H\}$$ induces
a $1-1$ correspondence between equivalence classes of essential
subgroups of $\wg$ and multiplicative subgroups of $\ff$.
\endproclaim

\demo{Proof} Recall  from Kummer theory  that
$\Gal(F(\sqrt{\dot F})/F)$ is the Pontrjagin dual of the discrete
group $\ff$ under the pairing  $(g,[f])=g(\sqrt{f})/\sqrt{f}$ of
$\Gal(F(\sqrt{\dot F})/F)$ with $\ff$, with values in
$\zz/2\zz\cong\{\pm 1\}$. (See \cite{ArTa, Chapter 6}.)

Assume that $H_{1}$ and $H_{2}$ are two essential subgroups of $\wg$
such that $P_{H_{1}}=P_{H_{2}}=:P$.  This means
$\frac{H_{1}\Phi(\wg)}{\Phi(\wg)}= \frac{H_{2}\Phi(\wg)}{\Phi(\wg)}$
because they are both the annihilator of $P$ under the pairing above.
(See \cite{Mo, Chapter 5}.)

Therefore, by considering the natural map $\wg\lra\Gal(F(\sqrt{\dot
F})/F)$ and the inverse image of $\frac{H_{1}\Phi(\wg)}{\Phi(\wg)}
=\frac{H_{2}\Phi(\wg)}{\Phi(\wg)}$ in $\wg$, we can conclude that
$H_{1}\Phi(\wg)=H_{2}\Phi(\wg)$.
This proves that $u$ is injective on equivalent classes of
essential subgroups.

In order to show that it is surjective, consider any subgroup $P$ of $\dot
F$ containing ${\dot F}^2$.  In the $\F_{2}$-vector space $\dot F/\dot
F^2$, the subspace $P/{\dot F}^2$ may be written as $P=\bigcap_{i\in
I}P_{i}$ where each $P_{i}$ is a hyperplane and $I$ is minimal with
this property.  Then again by Kummer theory and Pontrjagin's duality,
there exist elements $\sigma_{i}\in\wg, i\in I$ such that
$P_{\langle\sigma_{i}\rangle}=P_{i}$ for each $i\in I$.  (Here
$\langle\sigma_{i}\rangle$ is the group generated by $\sigma_{i}$ in
$\wg$.)

Set $H:=\langle\sigma_{i\mid \, i\in I}\rangle$, the closure of
the subgroup of $\wg$ generated by $\sigma_{i}, i\in I$.  Then from
Burnside's Basis Theorem for pro-$2$-groups (see e.g. \cite{Koc, Chapter
6} or \cite{Hal, Chapter 12}), we see that $\{\sigma_{i}, i\in I\}$
form a minimal set of generators of $H$ and their images modulo
$\Phi(H)$ form a topological $\F_{2}$-basis of $H/\Phi(H)$.  From the
choice of the $\sigma_{i}$'s we have
$$\frac{H}{H\cap\Phi(\wg)}\cong\frac{H.\Phi(\wg)}{\Phi(\wg)}\cong
\prod_{I}\zz/2\zz=\prod_{i\in I}\<\bar{\sigma_{i}}\>$$ where
$\bar{\sigma_{i}}$ is the image of $\sigma_{i}$ in
$\frac{H.\Phi(\wg)}{\Phi(\wg)}$.  On the other hand
$\frac{H}{\Phi(H)}\cong\prod_{i\in I}\<\bar{\sigma_{i}}\>$.  Hence,
the natural homomorphism $\frac{H}{\Phi(H)}\lra\frac{H}{\Phi(\wg)\cap
H}$ is an isomorphism and $\Phi(H)=\Phi(\wg)\cap H$ as desired.  This
shows that $H$ is essential.  Since $P_{H}=\bigcap_{i\in I}P_{i}=P$,
$u$ is surjective and the proof is complete.\qed
\enddemo

The motivation for this study of essential subgroups grew out of the
observation in \cite{MiSp1} that for $H \cong \zz/2\zz$, if $P_H \neq
\ff$ (i.e. if $H \cap \Phi(\wg ) =
\{1\}$), then $P_H$ is in fact the positive cone of some ordering on $F$.
The reader is referred to  \cite{L2} for further details on orderings and
connections to quadratic forms. Some convenient references for basic facts
on quadratic forms are \cite{L1}
and \cite{Sc}.

Since the presence or absence of $\zz/2\zz$ as an essential subgroup of
$\wg$
   determines the orderings or lack thereof on $F$,
one wonders whether other subgroups of $\wg$ also yield interesting
information about $F$.  We make the following definition.

\definition{Definition 1.4}

(1) Let $\Cal C$ denote the category of pro-2-groups of exponent at
most $4$, for which squares and commutators are central.
(Observe that
since each commutator is a product of (three) squares,
it is sufficient to assume that all squares are central.)
All W-groups
are in category $\Cal C$.  See \cite{MiSm2} for further details. Note
that $\Cal C$ is a full subcategory of the category of
pro-$2$-groups. This allows us to freely use all of the properties of
pro-$2$-groups.

(2) Let $H$ be an isomorphism type of groups.  An embedding
$\varphi\colon H \lra \wg$ is an {\it essential embedding} if
$\varphi(H)$ is an essential subgroup of $\wg$.  Note that if $H$
embeds in $\wg$, then $H$ has to be in category $\Cal C$.

(3) An {\it $H$-ordering} on $F$ is a set $P_{\varphi(H)}$ where
$\varphi$ is an essential embedding of $H$ in $\wg$.

(4) Let $(F,T)$ be a field with an $H$-ordering $T$.  We say that
$(L,S)$ {\it extends} $(F,T)$ if $L$ is an extension field of $F$ in
the maximal Galois $2$-extension $F(2)$ of $F$, $S$ is a subgroup of
$\dot L$ containing $\dot L^2$, $T = S \cap \dot F$, and the induced
injection $L/S \lra \dot F/T$ is an isomorphism.  We also say $(L,S)$ is a
$T${\it -extension} of $F$.  (We will see in Propositions 4.1 and 4.2
that maximal $T$-extensions always exist, and that a maximal such
extension $(L,S)$ in $F(2)$ has $S = \dot L^2$.)  An extension $(L,S)$
of $(F,T)$ is said to be an $H${\it -extension} if $S$ is an
$H$-ordering of $L$.

(5) An extension $(L,S)$ of $(F,T)$ is called an $H${\it -closure} if
it is a maximal $T$-extension which is also an $H$-extension.  Note
this implies $T = \dot L^2$ and $\Cal G_L \cong H$. Note also that we
will not consider maximal $H$-extensions $(K,S)$, because in general they
need not satisfy $S=\dot K^{2}$.
\enddefinition

We set the following notation: $C_n$ denotes the cyclic group of order $n$,
$D$ denotes the dihedral group of order $8$, $Q$ denotes the quaternion group
of order $8$.

   If $G_1$ and $G_2$ are in $\Cal C$, we denote by $G_1*G_2$ the free
   product (i.e. the coproduct) of the two groups in category $\Cal C$.
   Then $G_{1}$ and $G_{2}$ are canonically embedded in $G_1 * G_2$ and
   the latter can be thought of as $(G_1 \times [G_1,G_2]) \rtimes G_2$
   with the obvious action of $G_{2}$ on the inner factor.  (See
   \cite{MiSm2}.)  For example, $D \cong C_2 * C_2$.

Let $a \in \dot F \backslash \dot F^2$.  By a $C_4^a$-extension of a
field $F$, we mean a cyclic Galois extension $K$ of $F$ of degree 4,
with $F(\sqrt{a})$ as its unique quadratic intermediate extension.
Let $a,b \in \dot F$, independent modulo $\dot F^2$.  By a
$D^{a,b}$-extension of $F$ we mean a dihedral Galois extension $L$ of
$F$ of degree 8, containing $F(\sqrt{a},\sqrt{b})$, for which
$\Gal(L/F(\sqrt{ab})) \cong C_4$.  Observe that any $C_{4}$-extension
is a $C_{4}^{a}$-extension for an $a\in F$, and that any $D$-extension
is a $D^{a,b}$-extension for suitable $a,b\in \dot F$.

The following result is not hard to prove, and is a special case of
more general results in~\cite{Fr}. (See also~\cite{L1, Exercise VII.8}.)

\proclaim{Proposition 1.5} There exists a $C_4^a$-extension of $F$ if
and only if $a\in\dot F\setminus F^{2}$ and the quaternion algebra
$\left(\frac{a,a}{F}\right)$ is split.  There exists a $D^{a,b}$-extension
of $F$
if and only if $a,b\in\dot F$ are independent modulo squares and the
quaternion algebra $\left(\frac{a,b}{F}\right)$ is split.
\endproclaim

   This proposition is actually one of the main tools we use to link the
   Galois-theoretic properties of an essential subgroup $H$ of $\wg$ to
   the algebraic properties of an $H$-ordering.  Since we will need to
   refer to such extensions often in the sequel, we sketch the subfield
   lattice of a $D^{a,b}$-extension $L/F$.

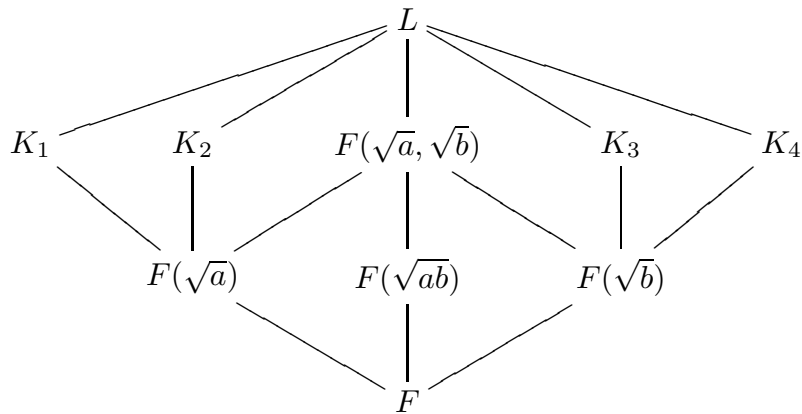
\bigskip

$$\diagram
& & L\dllline\dlline\dline\drline\drrline & & \\
K_{1}\drline & K_{2}\dline & F(\sqrt{a},\sqrt{b})\dlline\dline\drline &
K_{3}\dline & K_{4}\dlline\\
& F(\sqrt{a})\drline & F(\sqrt{ab})\dline & F(\sqrt{b})\dlline & \\
& & F & &
\enddiagram$$

\bigskip

The paper is organized as follows.

In \S 2, we show that the only abelian groups which can appear as
essential subgroups of a W-group are $C_2$ and $(C_4)^I$ where $I$ is
some nonempty set.  We also determine the possible nonabelian subgroups
   generated by two elements. In Theorem~2.7 we provide a
strong restriction on possible finite subgroups of a $W$-group. Some
of these results are important in determining the cohomology rings of
$W$-groups.

In \S 3 we show how properties of an $H$-ordering $T$, such as stability
under addition or rigidity, may be described in a Galois-theoretic way.
The definition and first properties of extensions and closures are given
in \S 4. We illustrate with Proposition~4.9 that even in a very
geometric situation, we cannot expect that every $H$-ordering $T$ admits
a closure. We also point out (Proposition~4.10) that this leads to a negative
answer to a strong version of the question asked in~\cite{Ma}: we
produce an example of a field $F$ having no field extension $F \lra
K$ with $W_{red}(K)\cong W(K)$, such that the induced map $W_{red}(F)
\lra W_{red}(K)$ is an isomorphism. (See also \cite{Cr2, Theorem 5.5},
from which one can also extract such examples. We address Craven's result
in
\S 4.) Later in
\S 8 we are able to provide a similar example
of a field $F$ with a subgroup $T$ of $\dot F$ such that the associated
Witt ring $W_T(F)$ is isomorphic to
$W(\qq_p), p\equiv 1(4)$ but again there is no field extension $F\lra K$
inducing the isomorphism
$W_T(K)\cong W(K)$. This example is interesting because $\mid\dot F/T\mid$
is finite. (For details see Example~8.14 and
the remark following this example.)

In \S 5 and \S 6 we study the case of essential subgroups $H$ generated
by $1$ or $2$ elements, and show that they admit closures.

In \S 7 we give a complete Galois-theoretic, as well as an algebraic
classification of rigid orderings, and in \S 8 we show that they admit
closures, provided that in the case of $C(I)$, the associated
valuation is not dyadic.  (See Theorem~8.15 and Example 8.14.)  In
Example~6.4 we see that the link between the additive structure of an
$H$-ordering and the Galois-theoretic properties of $H$ is not as
tight as we might have expected.  This leads us to investigate this
question more thoroughly in \S 9.  Actually, with a few natural extra
requirements on the Galois groups we are considering, this can be fixed.
We are then able to obtain a perfect identification between the two
aspects.

As we have already said, application of this theory to local-global
principles for quadratic forms will constitute the core of a
subsequent paper.  In the conclusion we illustrate by an easy example,
what we intend to do in this direction.

The authors would like to acknowledge Professors A. Adem,
J.-L. Colliot-Th\'{e}l\`{e}ne, T.~Craven, B. Jacob, D. Karagueuzian, J.
Koenigsmann,
T.-Y. Lam, D. Leep and H. W. Lenstra, Jr.  for valuable discussions
concerning the results in this
paper; and also the hospitality of the Mathematical Sciences Research
Institute at Berkeley, the Department of Mathematics at the University of
California at Berkeley, and the Mathematisches Forschungsinstitut at
Oberwolfach, which the authors were privileged to visit during the
preparation of this paper.

\head
\S 2.  Groups not appearing as  subgroups of $W$-groups
\endhead

In this section we show that no essential subgroup of $\wg$ can have
$C_2$ as a direct factor (except in the trivial case where the
subgroup is $C_2$), nor can $Q$ appear as a subgroup of $\wg$.  These
two facts will then be used to show that the four nonabelian groups
$C_{2}*C_{2}=D, C_{2}*C_{4},C_{4}\rtimes C_{4}$ and $C_4 * C_4$, together
with the abelian group $C_4 \times C_4$, comprise all of the possible
two-generator essential subgroups of W-groups.  Thus we have a good
picture of the minimal realizable and unrealizable subgroups.  We
further show that every finite subgroup of a W-group is in fact an
\lq\lq S-group\rq\rq\ as defined in \cite{Jo}.  (We shall call such
groups \lq\lq split groups\rq\rq\ here.)  The fact that $Q$ is not a
subgroup of $\wg$ is actually a consequence of this last result.

We often use the fact that there is a perfect $\F_2$-vector space
duality between the relations among the generators of $\wg$ and
the $F$-quaternion algebras in $\Br(F)$, the Brauer
group of $F$. (For a detailed exposition of this duality, see \cite{MiSp2,
Theorem~2.20}.)

Briefly, this duality occurs as follows: let $a_i, i
\in I$ form a basis for the $\F_2$-vector space $\ff$, where by abuse
of notation we identify an element $a \in \dot F$ with its image in
$\ff$.  We can choose a minimal set of (topological) generators
$\sigma_i, i \in I$ for $\wg$ having the property that
$\sigma_i(\root\of{a_i}) = - \; \root\of{a_i}$, and $\sigma_i$ fixes
$\root\of{a_j}$ for $i \neq j$.  Then the quaternion algebra
$(\frac{a_i,a_i}{F})$ is viewed as \lq\lq corresponding to\rq\rq\ the square
$\sigma_i^2$, and the quaternion algebra $(\frac{a_i,a_j}{F}), i \neq j$ is
viewed as \lq\lq corresponding to\rq\rq\ the commutator
$[\sigma_i,\sigma_j]$.
(Because the quaternion algebras $(\frac{a_i,a_j}{F})$ and
$(\frac{a_j,a_i}{F})$ are isomorphic and $[\sigma_i,\sigma_j] =
[\sigma_j,\sigma_i]$ we see that the order of $i$ and
$j$ is irrelevant, and we consider each quaternion algebra and each
commutator only one time.)

In order to explain this pairing in a more detailed way, we set 
$U\colon=\prod_{i}(C_{2})_{i}\times\prod_{(i,j)}(C_{2})_{(i,j)}$ to be
a topological group with the product topology, where each
$(C_{2})_{i}, i \in I$ is a discrete group of order $2$ with a formal
generator $\sigma_{i}^{2}$ and each $(C_{2})_{(i,j)}$ is a discrete
group of order $2$ with a formal generator
$[\sigma_{i},\sigma_{j}],i,j\in I,i\neq j$, and with the understanding
that $[\sigma_i,\sigma_j]=[\sigma_j,\sigma_i]$ so that
$(C_2)_{(i,j)}=(C_2)_{(j,i)}$.

We also set $P$ to be a set of degree $2$ homogeneous polynomials in
variables $Z_i, i \in I$ over a field $\Bbb{F}_2$.

Then by Pontrjagin's duality we have a perfect pairing $U\times P \lra
\{\pm 1\}$ such that the topological basis
$\{\sigma_i^2,[\sigma_i,\sigma_j],i\in I, (i,j)$ is an unordered pair
of distinct elements in $I$\} of $U$ is orthogonal to the vector basis
$\{Z_i^2, Z_i Z_j, i \in I, (i,j)\in I\times I, i\neq j\}$.  Then we
have a homomorphism $\psi\colon P\lra \Br (F)$ such that
$\psi(Z_i^2)$ is the class of $(\frac{a_i,a_i}{F})$ in $\Br (F)$ and
$\psi(Z_i,Z_j)$ is the class of $(\frac{a_i,a_j}{F})$ in $\Br (F)$.  The
kernel $Q$ of $\psi$ may be thought of as the group of relations between
the products of quaternion algebras $(\frac{a_i,a_i}{F})$ and
$(\frac{a_i,a_j}{F}), i, j \in I, i\notin j$ in $\Br (F)$.

The key fact proved in \cite{MiSp2, Theorem~2.20} tells us that the
group of relations $\nu$ between the products of $\sigma_i^2$ and
$[\sigma_i,\sigma_j],i\in I, i, j\in I, i \neq j$ is the annihilator
of $Q$ under the pairing above.  This allows us to conclude that the
pairing between $U$ and $P$ induces Pontrjagin's duality between $\nu$
and the group of quaternion algebras in $\Br (F)$.  (See \cite{MiSp2,
Corollary~2.21}.)

Using this Pontrjagin duality, we can say informally that the relations
among the generators of $\wg$ which
may all be expressed as products of squares and commutators of elements
$\sigma_i, i \in I$, are
\lq\lq dual\rq\rq\ to the corresponding product of quaternion algebras
in $\Br (F)$.  In particular, this means that if $\sigma_i,\sigma_j$
correspond to linearly independent elements $a_i,a_j$
in $\ff$ over
$\F_2$, respectively,
under this dual relation, then if $(\frac{a_i,a_i}{F}) = 1 \in \Br (F)$,
$\sigma_i^2$ will
not appear in any of the relations for $\wg$, and if $(\frac{a_i,a_j}{F}) =
1 \in
Br(F), i\neq j$, then
$[\sigma_i,\sigma_j]$ will not appear in any of the relations
for $\wg$.

\proclaim{Lemma 2.1}\cite{Mi1}, \cite{CrSm}  The groups $C_2\times C_2$
and
$C_4\times C_2$ cannot be realized as essential subgroups of $\wg$ for any
field $F$.
\endproclaim
\demo{Proof}  Assume $H = \< \sigma , \tau \mid \sigma^2 = \tau^2 =
[\sigma ,\tau
] = 1\> \subseteq \wg$ or $H = \< \sigma, \tau \mid \sigma^2 = [
\sigma, \tau ] = 1,
\tau^4 = 1 \>$, and assume $\sigma ,\tau ,\sigma\tau \notin \Phi
(\wg)$.   Then $-1
\notin \dot F^2$, for if $-1 \in \dot F^2$, we would have $(\frac{a,a}{F}) = 1
\in Br(F)$ for all $a \in \dot F$.  This means, in the relations for
$\wg$, no \lq\lq squared terms\rq\rq\ appear.  But if $H$ is a subgroup of
$\wg$, then $\sigma ^2$ appears as a relation in $\wg$.

Now consider a $D^{a,-a}$-extension $L/F$, where $\sqrt{a}$ is not fixed by
$\sigma$.  Such an extension exists since
$\sigma \notin \Phi (\wg)$, $|\ff| \ge 4$, and $-1 \notin \dot F^2$.
Consider $\<\bar\sigma ,\bar\tau\>$, the image of $H$ in $\Gal(L/F)$. We have
$\bar\sigma ^2$ = 1, so the fixed field of $\bar\sigma$ is of index 2 in
$L$ and
does not contain $\sqrt{a}$.  This means it cannot contain $\sqrt{-1}$ either,
but must be one of the two extensions of $F$ of degree 4 sitting over
$F(\sqrt{-a})$, so $(\sqrt{-1})^\sigma = -\sqrt{-1}$.  Now choose an element $b
\in \dot F \backslash\dot F^2$ for which $\sqrt{b}^\sigma = \sqrt{b}$ and
$\sqrt{b}^\tau = -\sqrt{b}$. Such an element $b$ exists since  $\sigma
,\tau ,\sigma\tau
\notin \Phi (\wg)$.  Consider the image $\<\bar\sigma ,\bar\tau\>$ of $H$
inside
the Galois group $G$ of a $D^{b,-b}$-extension $K$ of $F$. (Because
$(\sqrt{-1})^{\sigma} = -\sqrt{-1}$ we see that $-b$
is not a square in $F$, and we can conclude that the elements $b$ and
$-b$ are linearly independent when they are considered
as elements in $\ff$.)
The fixed field
$K_\sigma$ of $\bar\sigma$ cannot contain $\sqrt{-b}$, so it must be one of the
two subfields of index 2 in $K$ not containing $\sqrt{-b}$.  On the other hand,
the fixed field $K_\tau$ of $\tau$ cannot contain $\sqrt{b}$, so considering
the subfield lattice, we see that $K_\sigma \cap K_\tau = F$.
Then the
image of $H$ in $G$ generates $G$, which means $\sigma$ and $\tau$ cannot
commute.  This is a contradiction, so $H$ cannot exist as an essential subgroup
of $\wg$.
\qed
\enddemo

  From the lemma above we immediately obtain the following result, which
is used in \cite{AKMi} to investigate those fields $F$ for which the
cohomology ring
$H^*(\wg)$ is Cohen-Macaulay.

\proclaim{Corollary 2.2} Let $\sigma$ be any involution in $\wg \backslash
\Phi(\wg)$ and set $E_\sigma = \Phi(\wg) \times \langle \sigma\rangle$.  Then
the centralizer $Z(E_\sigma)$ of $E_\sigma$ in $\wg$ is $E_\sigma$ itself.
\endproclaim

\demo{Proof} If $\tau \in Z(E_\sigma)\setminus E_\sigma$ then $[\tau
,\sigma ] = 1$ and
$\langle\tau,\sigma\rangle = C_2 \times C_2$ or $C_4 \times C_2$,
where $\<\tau, \sigma\>$ is an essential subgroup of
$\wg$. From Lemma
2.1, this is a contradiction, and we see $\tau \in E_\sigma$ as desired.
\qed
\enddemo

\proclaim{Corollary 2.3} No essential subgroup of $\wg$ can have $C_2$ as
a direct factor (except in the trivial case where the subgroup is $C_2$).
\endproclaim
\demo{Proof}  Since $\Phi (H \times C_2) = \Phi (H)$, if $H\times C_2$ is a
subgroup with $\Phi (H\times C_2) = (H \times C_2) \cap \Phi (\wg )$, then the
$C_2$-factor is not in $\Phi (\wg )$.  Take any single element $\sigma \in
H\backslash \Phi (H)$. Then $\<\sigma\>\times C_2 \cong C_2\times C_2$ or
$C_4\times C_2$, which cannot be an essential subgroup. Therefore neither can
$H\times C_2$.
\qed
\enddemo

\proclaim{Proposition 2.4} The quaternion group $Q$ cannot appear as a
subgroup of $\wg$.
\endproclaim
\demo{Proof}  Suppose $Q = \<\sigma ,\tau | \sigma ^2 = \tau^2 = [\sigma ,\tau
]\> \subseteq \wg$.  Then as in the lemma above, $-1\notin\dot F^2$, since we
have $\sigma^2\tau^2 , \sigma^2[\sigma ,\tau ], \tau^2[\sigma ,\tau ]$ in the
relations for $\wg$.  Consider the image of $Q$ in any dihedral extension
of $F$ of order 8.  Since $Q$ is not isomorphic to $D$, this image must be a
proper quotient of $Q$, and therefore is elementary abelian.  Then the same
argument as in the lemma above shows that  $(\sqrt{-1})^\sigma = -\sqrt{-1}$.
Indeed we see again that there exists a
$D^{a,-a}$-extension $L/F$, where $\sqrt{a}$ is not fixed by
$\sigma$. As we have just observed, the image
$\bar{\sigma}$ of $\sigma$ in $\Gal(L/F)$ has order at most $2$,
and the order must be $2$ because
$\sigma(\sqrt{a}) = -\sqrt{a}$. Therefore we can again conclude that
$\sigma(\sqrt{-a}) = \sqrt{-a}$ and consequently
$\sigma(\sqrt{-1}) = -\sqrt{-1}$.
Choosing an element $b\in \dot F \backslash\dot F^2$ for which
$\sqrt{b}^\sigma = \sqrt{b}$ and $\sqrt{b}^\tau = -\sqrt{b}$ as before, we can
again take a $D^{b,-b}$-extension of $F$, and observe that the image of $Q$ in
the Galois group of this extension must generate the entire group, which is a
contradiction.
\qed
\enddemo

\proclaim{Theorem 2.5} The only groups generated by two elements which can
arise as essential subgroups of $\wg$ are the five groups $C_2*C_2$, $C_2*C_4$,
$C_4*C_4$, $C_4\times C_4$, and $C_4\rtimes C_4$.
\endproclaim
\demo{Proof} Let $H$ be generated by $x,y$.  We have an exact sequence
$$1\to\Phi (H) \to H \to C_2\times C_2 \to 1,$$
where $\Phi (H) \cong (C_2)^k$ is generated by $x^2, y^2, [x,y]$, so $k \leq
3$.  Then $|H| = 2^{k+2}$, so $|H| \leq 32$, and $|H| = 32$ if and only if
$|\Phi (H)| = 8$, if and only if $H \cong C_4*C_4$.  Otherwise $|H| = 8$ or
$16$, and there are only a few groups to consider.  If $|H| = 8$, necessarily
$H \cong C_2*C_2$, as all other groups of order 8 and exponent at most 4 either
have $C_2$ as a direct factor or are isomorphic to $Q$.

There are fourteen groups of order 16; among these, five are abelian, and by
Lemma 2.1 only $C_4 \times C_4$ among these can be an essential subgroup of
$\wg$.  Among the nine nonabelian groups, two have $C_2$ as a direct factor,
and four more have exponent 8. The remaining three are the groups $C_2*C_4$,
$C_4\rtimes C_4$, and $DC$, the central product of $D$ and $C_4$ amalgamating
the unique central subgroup of order 2 in each group.  This group, however, has
$Q$ as a subgroup (see \cite {LaSm}), so cannot be an essential subgroup
of
$\wg$.
\qed
\enddemo

That the group $Q$ cannot appear as a subgroup of any W-group is a special case
of a more general description of the kinds of groups which can appear as
essential subgroups of W-groups.  All finite subgroups must in fact be \lq\lq
split groups\rq\rq , which we define next.  These are the same as \lq\lq
S-groups\rq\rq\ as defined in \cite {Jo}.  The quaternion group $Q$ is not such
a group.

\definition {Definition 2.6} Let $G$ be a nontrivial finite group and $X = \{
x_1, x_2,
\dots , x_n\}$ be an ordered minimal set of generators for $G$.  We say
that $G$
satisfies the {\it split condition with respect to} $X$ if $\langle x_1
\rangle \cap [G,G]\langle x_2, \dots , x_n\rangle = \{ 1 \}$.  The group
$G$ is called a {\it split group} if it has a minimal generating set with
respect to which it satisifies the split condition.  We also take the trivial
group to be a split group.
\enddefinition

\proclaim {Theorem 2.7} Let $\Cal G_F$ be a W-group, and let $G$ be any finite
subgroup of $\Cal G_F$.  Then $G$ is a split group.
\endproclaim
\demo {Proof}  Each finite subgroup $H$ of $\Cal G_F$ can be written as $H = G
\times \prod _1^m C_2$ for some $m \in \Bbb N \cup \{ 0 \}$, where $G$ is
an essential subgroup of $\Cal G_F$ \cite {CrSm}.  Thus it is enough to prove
the theorem for $G$ a finite essential subgroup of $\Cal G_F$.

Then let $G$ be such a group and let $P_G$ be the associated $G$-ordering.
Let $\dot F / P_G = \langle a_1P_G, \dots , a_nP_G\rangle$ so that the cosets
$a_iP_G$ give a minimal generating set for $\dot F/P_G$.   Further set $\{
\sigma _1, \dots , \sigma _n\}$ to be a minimal generating set for $G$
such that $\sigma _i(\sqrt{a_j}) = (-1)^{\delta _{ij}}\sqrt {a_j}$ where
$\delta _{ij}$ is the Kronecker delta. (This is possible because $G$ is an
essential subgroup of $\Cal G_F$, so that a minimal set of generators for $G$
can be extended to a minimal (topological) generating set of $\Cal G_F$.)

Assume first that we can choose the representatives
$a_i$ in such a way that $a_1t_1 + a_1t_2 = f^2
\in \dot F^2$ for some $t_1, t_2 \in P_G$. (Note that this is equivalent to
saying
that $a_1 \in P_G + P_G$.)  In this instance, there are two cases to consider.

First, suppose that $t_1, t_2$ are congruent mod $\dot F^2$.  Then there
exists $g \in \dot F$ such that $a_1t_1 + a_1t_1g^2 = f^2$, and so $a_1t_1f^2
= (a_1t_1)^2 + (a_1t_1g)^2$, and $a_1t_1$ is a sum of two squares in $F$ which
is not itself a square.  Thus we have a $C_4^{a_1t_1}$-extension $L$ of $F$.
We claim that $G$ satisfies the split condition with respect to $\{ \sigma _1,
\dots , \sigma _n\}$.  Checking this condition is equivalent to
showing $\sigma _1^2 \notin [G,G]\langle \sigma _2, \dots , \sigma
_n\rangle$.
Suppose it is not true.  Then we have an identity $\sigma _1^2 \prod _{1
\leq i < j
\leq n}[\sigma _i,\sigma _j]^{\epsilon _{ij}}\prod _{k = 2}^n\sigma
_k^{2\epsilon _k} = 1$ in $G$, where $\epsilon _{ij}, \epsilon _k \in \{
0,1\}$.  Restricting to $L$ we see that $\sigma _1^2 |_L = 1$.  This cannot be
the case as $\sigma _1$ does not fix $\sqrt {a_1t_1}$.  Thus in this case $G$
is a split group.

Next suppose that $t_1\dot F^2 \neq t_2\dot F^2$.  In this case we can find a
$D^{a_1t_1, a_1t_2}$-extension $L/F$.  Assuming again that $G$ does not
satisfy
the split condition with respect to $\{ \sigma _1, \dots , \sigma _n\}$, we
again have an identity $\sigma _1^2 \prod _{1 \leq i < j
\leq n}[\sigma _i,\sigma _j]^{\epsilon _{ij}}\prod _{k = 2}^n\sigma
_k^{2\epsilon _k} = 1$ in $G$, where $\epsilon _{ij}, \epsilon _k \in \{
0,1\}$.  Since each of the $\sigma _i, i = 2, \dots , n$ acts
trivially on $F(\sqrt{a_1t_1}, \sqrt{a_1t_2})$, we see that each $\sigma _i, i
>   1$ is central when restricted to $L$.  Thus again $\sigma _1^2 |_L
= 1$.  But
$\sigma _1 |_L$ generates $\Gal (L/F(\sqrt{a_1t_1}\cdot \sqrt{a_1t_2}))
\cong C_4$. Hence $G$ is a split group.

Finally, assume that we cannot choose $a_1 \in P_G + P_G$.  Then necessarily
$P_G + P_G \subseteq P_G \cup \{0\}$.  If $-1 \in P_G$, then $P_G = \dot F$ and
$G = \{ 1\}$ which is a split group.  Otherwise $P_G$ is a preordering in $F$,
and we may write $P_G = \cap _{i = 1}^n P_i$ where each $P_i$ is an ordering,
and each $P_i = \{\,f \in \dot F \mid \sqrt{f}^{\sigma _i} = \sqrt{f}\,\}$.
Then $\{\sigma _1 , \dots , \sigma _n\}$ is a minimal generating set for
$G$.  Furthermore, each $\sigma _i^2 = 1$.  (See \cite {MiSp1} for
details. The definition of a preordering in a field
$F$ can be found in \cite {L2, Chapter 1}, together with the basic
properties of preordered rings.)
Thus again we see that $G$ is a split group.
\qed
\enddemo

\proclaim {Corollary 2.8} Each nontrivial finite subgroup $G$ of a
W-group $\Cal G_F$ can
be obtained inductively from copies of $C_2$ and $C_4$ by taking
semidirect products at each step.  Thus we have $G = G_n \supseteq G_{n-1}
\supseteq \dots \supseteq G_1 \supseteq G_0$ where $G_0 \in \{ C_2, C_4\}$,
and $G_i = G_{i-1}\rtimes C_2$ or $G_i = G_{i-1}
\rtimes C_4$ for each $i = 1, \dots , n$.
\endproclaim
\demo {Proof}  We proceed by induction on the number of generators of $G$.
The statement clearly holds for any group $G$ generated by a single element.
Let $G$ be any (nontrivial) finite subgroup of
the W-group $\Cal G_F$.  Then we can write $G = H \times \prod _1^m C_2$
where $H$ is essential, and $G$, if not equal to $H$, is clearly built up as
described from $H$, where the action in the semidirect product is trivial. We
can choose a minimal set of generators $\{ \sigma _1, \dots , \sigma _n\}$ for
$H$ such that $H$ satisfies the split condition with respect to these
generators.  Clearly $N := [H,H]\langle \sigma _2, \dots , \sigma _n\rangle$
is a
normal subgroup of $H$, and $H \cong N \rtimes \langle \sigma _1\rangle$, where
$\langle \sigma _1\rangle \cong C_2$ or $C_4$.  Since $N \cong
\langle \sigma_2,\dots \sigma_n\rangle \times \prod_1^kC_2$ (for some
positive integer $k$),  we  finish by induction.
\qed
\enddemo

\example{Example 2.9} Consider the W-group $\Cal G_2$ of the 2-adic
numbers
$\qq _2$.  It has the presentation $\langle \sigma, \tau, \rho \mid
\sigma^2[\tau,\rho ]\rangle$ in the category $\Cal C$ of groups of exponent at
most four with squares and commutators central.
(See \cite{MiSp2, Example 4.4}.)
A basis for $\ff$ is given by
$\{ -1, 2, 5\}$, and $\sigma$ may be chosen to fix $\sqrt{2}$ and
$\sqrt{5}$ but
not $\sqrt{-1}$, $\tau$ to fix $\sqrt{-1}$ and $\sqrt{5}$ but not $\sqrt{2}$,
and $\rho$ to fix $\sqrt{-1}$ and $\sqrt{2}$ but not $\sqrt{5}$.  Then $\Cal
G_2$ can be constructed inductively from copies of $C_4$ and $C_2$
using semidirect products as follows:
$$
\align G_0 &= \langle \rho\rangle \cong C_4\\
G_1 &= G_0 \times \langle [\sigma ,\rho ]\rangle \cong G_0 \times C_2\\
	G_2 &= G_1 \rtimes \langle \sigma \rangle \cong G_1 \rtimes C_4\\
G_3 &= G_2 \times \langle [\sigma ,\tau ]\rangle \cong G_0 \times C_2\\
\Cal G_2 &= G_3 \rtimes \langle \tau \rangle \cong G_3 \rtimes C_4
\endalign
$$
Thus $\Cal G_2 \cong \{ [ ( C_4 \times C_2 ) \rtimes C_4 ]
\times C_2 \} \rtimes C_4$.
\endexample

Corollary 2.8 is an interesting generalization of the known structure
of W-groups associated with Witt rings of finite
elementary type. In fact, all W-groups associated with Witt rings of
finite elementary type can
easily
be seen to be built up from cyclic groups of order 2 or 4, using only
semidirect products.  First one checks that the groups associated with basic
indecomposable groups are such groups.  Then the group ring construction for
Witt rings corresponds directly to taking a semidirect product with a cyclic
group of order 4, while the direct product construction for Witt rings
corresponds to taking a free product of W-groups in the appropriate category.
But this in turn just involves taking a direct product with an appropriate
number of copies of $C_2$ (representing the necessary commutators) and
then taking a semidirect product with the generators of one of the initial
W-groups.  See \cite{MiSm2} for details.

Corollary 2.8 is quite useful for the investigation of cohomology rings of
W-groups.  This is important in light of the recent proof of the Milnor
Conjecture by Voevodsky \cite{Vo}.  In particular, Voevodsky's result shows
that the cohomology rings of absolute Galois groups with $\Bbb
F_2$-coefficients carry no more information about the base field than
Milnor's K-theory $\mod 2$.  On the other hand, the cohomology rings
of W-groups
carry substantial additional information.  (See \cite {AKMi}.)

Using \cite {Jo: Cor, p. 370} and Theorem 2.7 above, we immediately obtain the
following.

\proclaim {Corollary 2.10} Let $G$ be any nontrivial finite subgroup of a
W-group $\Cal G_F$.  Then the cohomology ring $H^*(G,\Bbb F_2)$ contains
nonnilpotent elements of degree 2, and hence of every even degree.
\endproclaim

\head
\S 3.   Galois groups and additive structures (1)
\endhead

In this section we give a simple Galois-theoretic characterization
of two important additive properties of $H$-orderings: stability
under addition and rigidity.  This generalizes the results on rigidity
and on the realizability of certain Galois groups obtained in
\cite{MiSm1}.

 For the rest of this paper, unless otherwise mentioned,
subgroups of $\wg$ will always be {\it essential}.
   Throughout this paper we write $T + aT = \{ t_1 + at_2 \mid t_1, t_2 \in
   T \cup \{ 0\} , t_1 + at_2 \neq 0 \}$, so $T$ and $aT$ are always
   subsets of $T + aT$, and $T + aT \supseteq \dot F^2$. (Here $T$ is
any subgroup of $\dot F$ containing all squares in
$\dot F$.)

\proclaim{Proposition 3.1}  Let $H$ be a subgroup of $\wg$, and $T$ its
associated $H$-ordering.  Then $H$ has $C_4$ as a quotient if and only if $T +
T \neq T$.
\endproclaim
\demo{Proof}  First assume there exists $a \in T + T$ which is not in $T$. Let
$K$ be the fixed field of $H$ in $F^{(3)}$. We construct a $C_4^a$-extension
$F_1$ of $F_0 = F(\sqrt{T}) = K \cap F^{(2)}$ inside $F^{(3)}$.  Then $L =
KF_1$ is a $C_4^a$-extension of $K$ in $F^{(3)}$, showing $H$ has $C_4$ as a
quotient.  We may write $a = t_1 + t_2$, so $a^2 - at_1 = at_2$.  Let $y = a -
\sqrt{a}\sqrt{t_1} \in F_0(\sqrt{a})$, so $N_{F_0(\sqrt{a})/F_0}(y) = [a] \in
\dot F_0/\dot F_0^2$.  Then $F_1 = F_0(\sqrt{a},\sqrt{y})$ is a
$C_4^a$-extension of $F_0$.  Since $yy^\sigma = y^2$ or $at_2 \in (\dot
F_0(\sqrt{a}))^2$   for all $\sigma \in \Gal (F_0(\sqrt{a})/F)$, we
see $F_1$ is
Galois over $F$, and hence is contained in $F^{(3)}$.

Conversely, assume $T + T = T$. If $-1\in T$, then $T=\dot F$ and
$H=\{1\}$. If $-1\notin T$, then $T$ is a preordering, so $T $ is
intersection of orderings, and $H$ is generated by involutions. (See
\cite{CrSm, Proposition~3.1}.)
Thus $H$ cannot have $C_4$ as a quotient.
\qed
\enddemo

\remark {Remark}  If $H$ has a $C_4$-quotient, then there exists a
$C_4^a$-extension of $F_0$ where we may take $a$ to be in $F$.  However, it is
not necessarily the case that $a \in T + T$.  That is, the quaternion algebra
$\left(\frac {a,a}{F(\sqrt{T})}\right)$ is split, so $a$ can be represented
as the sum of
two squares in $F(\sqrt{T})$, but not necessarily as the sum of two elements in
$T$. This can be seen in Example~6.4.
\endremark

The following definition generalizes the notion of the rigidity of a
field, and introduces the notion of the level of $T$.
(See \cite{Wa, page 1349}.)

\definition {Definition 3.2} Let $T$ be a subgroup of $\ff$.
We say that $T$ has {\it level} $s$ if $-1$ is a sum of $s$ elements
of $T$, and not a sum of $s-1$ elements of $T$. We say that this
level is infinite if $-1$ is not such a sum for any natural number $s$.
We say that the field $F$ is {\it T-rigid}, or equivalently that $T$ is
rigid, if for every $a \notin  T\cup -T$, we have $T + aT \subseteq T \cup aT$.
\enddefinition

We have the following easy-to-prove but important property of rigid
$H$-orderings:

\proclaim{Proposition 3.3} Let $T$ be  a rigid $H$-ordering on $F$. Then

(1) The level of $T$ is $1,2$ or infinite.

   (2) If the level of $T$ is $2$, then $T+T=T\cup-T$.
\endproclaim

\demo{Proof}  Let $T$ be an $H$-ordering of finite level $s>1$ and let us
write
$-1=a+a_{s}$ with $a=a_{1}+\ldots+a_{s-1}$  and $a_{i}\in T$ for $i=1,\ldots,s$.
If $a\in T\cup -T$ then since $a \notin -T$ we see $a \in T$ and
$s$ must be $2$. Thus we may assume  $a\notin
T\cup -T$. If $T$ is rigid, then $-1=a+a_{s}\in T+aT=T\cup aT$. This is a
contradiction, proving (1).

Assume the level of $T$ is $2$.  Then $-1\in T+T$ and $T\cup
-T\subseteq T+T$. Suppose there is $a\in (T+T)\setminus (T\cup -T)$ and
let us write $a=s+t, s,t\in T$. Then of course $-a\notin T\cup -T$
and we have $-t=s-a\in T+(-a)T=T\cup -aT$ by rigidity. But $-t\notin
T$ because the level is $2$, and $-t\notin -aT$ because $a\notin T$.
This is again a contradiction, proving (2).\qed
\enddemo

\proclaim{Proposition 3.4}Let $H$ be a subgroup of $\wg$, and let $T$
be an $H$-ordering.  Assume $-1 \in T$.  The following are equivalent.
\roster
\item $F$ is $T$-rigid.
\item $D$ is not a quotient of $H$.
\item $H$ is abelian.
\endroster
\endproclaim
\demo{Proof} We will show $(2) \implies (1) \implies (3) \implies (2)$.
For the first implication, we show the contrapositive.  Thus assume that $F$ is
not $T$-rigid.  Let $K$ be the fixed field of $H$, and let $F_0 = K \cap
F(\sqrt{\dot{F}})
= F(\{ \sqrt{t} : t \in T\})$.  We will construct a $D$-extension
$F_1$ of $F_0$ inside $F^{(3)}$.  Then $L = KF_1$ will be a $D$-extension
of $K$ in $F^{(3)}$, showing that $H$ has $D$ as a quotient.  Since $F$ is
not $T$-rigid and $-1 \in T$, there exist $a,b \in \dot F \backslash T$ such
that
$b = t_1 - at_2$, where $t_1, t_2 \in T$ but $b \notin T \cup aT$. Let $y =
\sqrt{t_1} + \sqrt{a}\sqrt{t_2} \in F_0(\sqrt{a})$, and let $F_1 =
F_0(\sqrt{a},\sqrt{b},\sqrt{y})$.  Notice that $yy^\sigma \in \{\pm y^2,\pm
b\} \subseteq F_0(\sqrt{a},\sqrt{b})^2 \,$ for all $\sigma \in
\Gal(F_0(\sqrt{a},\sqrt{b})/F)$, so $F_1/F$ is Galois, and $F_1 \subseteq
F^{(3)}$.  Then the usual argument (see \cite{Sp} or \cite{Ki,
Theorem 5}) shows $\Gal(F_1/F_0)
\cong D$.

Now assume $F$ is $T$-rigid. To see that $H$ is abelian, it is
sufficient to show that
for all $\sigma , \tau \in H$, the restrictions of $\sigma , \tau$ to any
$D$-extension $L$ of $F$ commute. (This is because $F^{(3)}$ is the compositum
of all quadratic, $C_4$- and $D$-extensions of $F$. (See \cite{MiSp2,
Corollary 2.18}.)
Thus if $\sigma , \tau$
commute on all $D$-extensions, they commute in $\wg$.)  Let $D^{a,b}$ be some
dihedral quotient of $\wg$, and let $L$ be the corresponding extension of $F$.
Denote as $\bar{\sigma} , \bar{\tau}$ the images of $\sigma$ and $\tau$ in
in $D^{a,b}$ and suppose $[\bar{\sigma}, \bar{\tau}]
\neq 1$.  Then $\sigma , \tau$ must each move at least one of
$\sqrt{a},\sqrt{b}$, and they cannot both act in the same way on these square
roots.  That implies $a, b, ab \notin T$. But $(\frac{a,b}{F})$ splits,
so $b \in
F^2 - a F^2 \subseteq T - aT = T + aT = T \cup aT$ by (1). Since $b
\notin T$, we have $b \in aT$, which contradicts the fact that $ab \notin T$.
Thus $[\sigma , \tau] = 1$.

The final implication is trivial.
\qed
\enddemo

\proclaim{Proposition 3.5}
Let $H$ be a subgroup of $\wg$, and let $T$ be an $H$-ordering.  Assume $-1
\notin T$. Let $K$ be the fixed field of $H$, and let $H_0$ be the subgroup of
$H$ which is the Galois group of $F^{(3)}/K(\sqrt{-1})$. The following are
equivalent.
\roster
\item $F$ is $(T \cup -T)$-rigid.
\item $D$ is not a quotient of $H_0$.
\item $H_0$ is abelian.
\item Every $D$-extension of $K$ in $F^{(3)}$ contains $K(\sqrt{-1})$.
\endroster
\endproclaim
\demo{Proof} Let $S = T \cup -T$.  Then $S$ is clearly an $H_0$-ordering, and
the equivalence of the first three statements follows from the preceding
proposition.  If there exists a $D$-extension $L$ of $K$ not containing
$K(\sqrt{-1})$, then $L(\sqrt{-1})$ will be a $D$-extension of $K(\sqrt{-1})$,
and $H_0$ will have $D$ as a quotient.  This shows $(2) \implies (4)$.
Finally, assume there exist $\sigma, \tau \in H_0$ which do not commute. Then
there exists some $D^{a,b}$-extension $M$ of $F$ such that $\Gal(M/F) = \langle
\bar{\sigma},\bar{\tau}\rangle$, where we denote $\bar{\sigma}$ and
$\bar{\tau}$ the images of $\sigma$ and $\tau$ in
$\Gal(M/F)$.
Then $\sigma$ and $\tau$ each move one of $a,b$
and cannot act in the same way on each.  Thus $a,b, ab \notin S$, but
$(\frac{a,b}{F})$ splits, so $b \in F^2 - a F^2 \subseteq T - aT$.  This
gives a $D$-extension $MK$ of $K$, which, since $a, b, ab \notin S$, does not
contain $\sqrt{-1}$. This shows $(4) \implies (3)$.
\qed
\enddemo

\head
\S 4. Maximal extensions, closures and examples
\endhead

Given any $C_2$-ordering $P$ on a field $F$, one can find a real closure of
$F$ with respect to that ordering, i.e. a real closed field $L$, algebraic over
$F$, with $P = \dot L^2 \cap F$.  Specifically, set $E = F(\sqrt{P})$.
(This means that $E$ is the compositum of all field
extensions $F(\sqrt{p}), p \in P$.)
Then $E$
is formally real, and a real closure $L$ of $F$ in $\bar F$ contains $E$,
$\dot{L}=\dot{L}^2\cup
-\dot{L}^2$, and $\dot{L}^2$ is an ordering of $L$.
($\bar{F}$ here means an algebraic closure of $F$.)
Then $\Gal(\bar F/L)$ is $\<\tau\> \cong C_2$, and we have $P =
\{ a
\in \dot F | \sqrt{a}^\tau = \sqrt{a}\}$.  Notice that for our purposes nothing
is lost by considering a real closure of $E$ inside $F(2)$, i.e. the
euclidean closure, rather than a real closure within the algebraic closure
$\bar
F$ of $F$. (See \cite{Be1}.)
This observation motivated the definition of $H$-closure given in
Definition 1.4.  The following two propositions show that maximal
$T$-extensions
always exist, i.e. that given any subgroup $T$ of $\dot F$, containing $\dot
F^2$, we can find a $T$-extension $(L,\dot L^2)$ of $(F,T)$ in $F(2)$. Thus
the real problem is in showing that $H$-closures exist, i.e. in showing that
$\Cal G_L \cong H$, or that we can find an $H$-extension which is a maximal
$T$-extension.

\proclaim {Proposition 4.1} Let $T$ be a subgroup of
$\dot{F}/\dot{F}^2$.  Then $(F,T)$ possesses a maximal $T$-extension.
\endproclaim

\demo{Proof} Let $\Cal S$ be the set of $T$-extensions
$(L,S)$ of $(F,T)$ inside $F(2)$, where $\Cal S$ is ordered under
inclusion. (See Definition~1.4 (4).)
More precisely we may say that $(L_1, S_1) \leq (L_2, S_2)$ if $L_1
\subset L_2$ and $S_2 \cap L_1 = S_1$. Observe that since both $(L_1,S_1)$
and $(L_2,S_2)\in \Cal S$, we automatically have a natural isomorphism
$L_1/S_1\lra L_2/S_2$. Then $\Cal S$ is nonempty, since $(F,T) \in \Cal
S$.  Now consider a totally ordered family $(F_j,T_j)$ in $\Cal S$.  Let
$K =
\cup F_j$, $S = \cup T_j$.  We will show $(K,S)$ is an upper bound for
the family $(F_j,T_j)$ in $\Cal S$.  First observe $T = S \cap F$ by
definition.  Thus $\dot
F/T \cong \dot F_j/T_j \to \dot K/S$ is one-to-one.  This map is also
onto, since if $b \in \dot K$, then $b \in F_j$ for some $j$, and $[b]_S$
is the image of $[b]_{T_j}$.  Then by Zorn's Lemma $\Cal S$ contains a
maximal element, which is a maximal $T$-extension of $(F,T)$.  \qed
\enddemo

\proclaim {Proposition 4.2} Let $(K,S)$ be a maximal $T$-extension of
$(F,T)$.  Then $S = \dot K^2$.
\endproclaim
\demo{Proof}  Let $\{ a_i : i \in I\}$ be a basis for $\dot F/T$ which
lifts to a basis for $\dot K/S$, which we can do because $S\cap F = T$
and $\dot K/S \cong \dot F/T$.  Assume $S \neq \dot K^2$, and choose $c \in S -
\dot K^2$.  Let $L = K(\sqrt{c})$, so $\dot L^2 \cap K = \dot K^2 \cup
c\dot K^2
\subseteq S$, and $\{ a_i : i \in I\}$ remain independent in $\dot L/\dot
L^2$. Let $\{ a_i : i \in I\} \cup \{ b_j : j \in J\} \cup \{ c\}$ be a basis
for $\dot K/\dot K^2$ such that $\{ b_j : j \in J\} \cup \{ c\}$ forms a basis
for $S$.  Then $\{ a_i : i \in I\} \cup \{ b_j : j \in J\}$ can be extended to
a basis $\{ a_i : i \in I\} \cup \{ b_j : j \in J\} \cup \{ c_{j'} : j' \in
J'\}$ for $\dot L / \dot L^2$.  Let $S'$ be the subgroup of $\dot L / \dot L^2$
generated by $\{ b_j : j \in J\} \cup \{ c_{j'} : j' \in
J'\}$.  Then $S' \cap K = S$, so $S' \cap F = T$, and $\dot L/S' \cong \dot K/S
\cong \dot F/T$, contradicting the maximality of $(K,S)$.  Thus we conclude $S
= \dot K^2$.
\qed
\enddemo

\proclaim{Corollary 4.3} An $H$-ordered field $(F,T)$ is an $H$-closure if
and only if
$T = \dot F^2$.
\endproclaim
\demo{Proof}  If $(F,T)$ is an $H$-closure, then it is also a maximal
$T$-extension, and $T = \dot F^2$ by the preceding proposition.
Conversely, suppose $T = \dot F^2$.  Let $L \supset F$ be any proper
extension of $F$ in $F(2)$.  Then $L$ contains a quadratic extension
of $F$, so $\dot L^2 \cap F \supsetneq \dot F^2$ and $L$ cannot extend
$(F,T)$.  This shows that $(F,T)$ is its own maximal $T$-extension,
and as it is an $H$-ordering, it is an $H$-closure.
\qed
\enddemo

Thus we see that our main task will be to show that there exists a
maximal $T$-extension $(K,\dot K^2)$ for an $H$-ordered field, which is
itself
$H$-ordered, i.e.  for which $\Cal G_K \cong H$.  The rest of the
section is devoted to the study of usual preorderings.
We will see in particular that in some important cases, usual
preorderings do not admit
closures.
Although this is in some sense a negative result, we shall see that
these examples are very interesting, and
that they deserve careful analysis. (See Proposition~4.10 below.)

\bigskip

Suppose $F$ is a formally real field equipped with a preordering $T$.
By \cite{CrSm, Proposition~3.1},  $T$ is an $H$-ordering for an $H$
generated by involutions, and conversely, any $H$-ordering
with $H$ generated by involutions has to be a preordering.
Thus, if $(L,P)$ is an $H$-extension of $(F,T)$, $P$ is a preordering
in $L$. In the same direction we have the following.

\proclaim{Lemma 4.4} Let $(F,T)$ be an $H$-ordered field with $T$ a
preordering, and assume $(K,S)$ is an $H$-closure of $(F,T)$. Then for
any intermediate extension $L/F$ of $K/F$, the pair
$(L,L\cap S)$ is a
$T$-extension of $(F,T)$ and $L\cap S$ is a preordering of $L$.
\endproclaim
\demo{Proof} Because the composite of the injective maps
$\dot{F}/T\lra\dot{L}/(L\cap S)\lra\dot{K}/S$ is bijective, each
injection is bijective and the intermediate extension is a
$T$-extension. Because $T$ is a preordering, $H$ is generated
by involutions and $S$ is also a preordering, forcing $L\cap S$ to be
a preordering as well.
\qed
\enddemo

We fix some notation.  For any field $k$ let $X(k)$ denote the
space of usual orderings of $k$. Within this section we will sometimes
assume that usual orderings contain $0$. This
will be clear from the context and should not cause confusion.
For $U\subseteq k$, let $\hat{U}$ be
the set of orderings of $k$ containing $U$.  When
$U=\{g_{1},\ldots,g_{n}\}\subset k$, we will denote the Harrison open
set $\hat{U}=\{\beta\in X(k)\,\mid\, g_{1},\ldots,g_{n}\in\beta\}$ by
$D_k(g_{1},\ldots,g_{n})$.  A set of the form $D_{k}(g)$ will be
called {\it principal}.  If $F\lra L$ is an extension, we denote by
$\pi\colon X(L)\lra X(F)$ the restriction map $\beta\mapsto \beta\cap
F$.

Let us  recall a few basic properties of the real spectrum which will be
useful in the sequel. They can be found in \cite{BCR, Chapter 7}.
Two other very nice introductions to the real spectrum are given in
\cite{Be2} and \cite{L3}.  The {\it real spectrum} $\Sper A$ of a ring $A$
is the set of pairs $({p},\alpha)$ with ${p}$ a prime
ideal of $A$ and $\alpha$ an ordering on the residue field $k({p})$
(the quotient field
of
$A/p$).
It is equipped with a topology generalizing the Harrison topology,
given by the subbasis of sets $D(f):=\{(p,\alpha)\in\Sper A\,\mid\,
f(\alpha)>0\in k({p})\}$, $f$ being any element of $A$. When $A$ is a
field $k$, $\Sper A$ is just the space of orderings
$X(k)$. Because $p=\alpha\cap-\alpha$ we see that
$\alpha$ already determines the prime $p$, and we will use
$\alpha$ instead of $(p,\alpha)$. As in the field case,
$\alpha$ may also be thought of as a subset of $A$, called the
``positive cone'' of the elements $f\in\alpha$ such that
$f(\alpha)>0$ in $k(p)$. Then we may write either $f(\alpha)>0$ or
$f\in\alpha$, according to our needs.

When $V$ is an affine algebraic variety over a real closed field $R$,
then the set of $R$-points $V(R)$ embeds continuously (with respect to
the euclidean topology of $V(R)$) as a dense subset in $\Sper R[V]$.
This embedding induces a $1-1$ correspondence $C\mapsto \tilde{C}$
between the semi-algebraic sets of $V(R)$ and the constructible sets
of $\Sper R[V]$.  ( See \cite{BCR, Theorem 7.2.3}.)

On the other hand, the map $R[V]\lra R(V)$ induces an embedding
$X(R(V))=\Sper R(V)\lra\Sper R[V]$ wich has the following properties
(See \cite{BCR, \S 7.6}):

(1) If $V$ is smooth, the embedding is dense.

(2) For every constructible open set $C\subseteq X(R(V))$ there exists a
constructible open set $D\subseteq\Sper R[V]$ such that $D\cap X(R(V))=C$.

(3) If $D_{1},D_{2}$ are constructible sets in $\Sper R[V]$ coinciding on
$X(R(V))$, then they coincide on $V(R)$ up to a positive codimensional
set.

\proclaim{Lemma 4.5}Let $V$ be an algebraic variety over a real
closed field $R$ and let $F$ denote the function field $R(V)$ of $V$.
Any nonempty open set $U$ of $X(F)$ contains a nonempty principal open
set $D_F(u)$.
\endproclaim

\demo{Proof} Since we are dealing with the function field $F$ of $V$,
we may always assume $V$ is affine and smooth.  As any nonempty open set
contains a constructible nonempty open set, we may assume $U\subseteq
X(F)$ is constructible.  Then we know that there exists a
nonempty constructible open set $U_{1}$ in $\Sper R[V]$ such that
$U=U_{1}\cap X(F)$, and
$U_{2}:=U_{1}\cap V(R)$ is not empty.
Denote by $(x_{1},\ldots,x_{n})$ the
coordinates of the ambient space $\R^{n}$ of $V$.
Let $a=(a_{1},\ldots,a_{n})\in
U_{2}$.  There is a ball of some radius $\epsilon$, centered on $a$
and contained in $U_{2}$.  If
$u=\sum_{i=1}^{n}(x_{i}-a_{i})^{2}-\epsilon$, then the ball $D(-u)$
satisfies $\emptyset\neq D(-u)\subseteq
U_{1}$ in $\Sper R[V]$ and we also have $D_{F}(-u)\subseteq U$ in
$X(R(V))$.  Since the embedding $X(F)\lra\Sper R[V]$ is dense,
$D_{F}(-u)\neq\emptyset$.
\qed
\enddemo

\definition{Remark 4.6} Note that Lemma~4.5 is not true for a general
formally real field. Consider for example the field $F=\R((X))((Y))$ of
iterated
   power series. We have $|X(F)|=4$ and each singleton is open
   and does not contain any principal set because it is not principal.
   \enddefinition

\proclaim{Lemma 4.7} Let $F$ be the function field of an algebraic
$R$-variety over a real
closed field $R$.
Let $T$ be a preordering of $F$ such that $\hat{T}$ is open in $X(F)$,
let $s\in T\setminus F^{2}$ and
$L=F(\sqrt{s})=F[Z]/(Z^{2}-s)$.

(1) For any preordering $P$ on $L$ such
that $T=P\cap F$,  the restriction  map $\pi\colon X(L)\lra X(F)$
induces a surjection from $\hat{P}$ to $\hat{T}$.

(2) If $(L,P)$ is a $T$-extension of $(F,T)$, then $\pi$ induces a
bijection between $\hat{P}$ and $\hat{T}$.
\endproclaim

\demo{Proof}
(1) Suppose there exists $\alpha\in \hat{T}$ such that
$\pi^{-1}(\alpha)=\{\beta_{1},\beta_{2}\}$ does not intersect $\hat{P}$.
Then neither $\beta_{1}$ nor $\beta_{2}$ contains $P$, and there exist
$f_{1}\in P\setminus\beta_{1},f_{2}\in P\setminus\beta_{2}$. One of
the three elements
$f\in\{f_{1},f_{2},f_{1}f_{2}\}$ must satisfy $f\in
P\setminus\{\beta_{1},\beta_{2}\}$, so
$f(\beta_{1})<0,f(\beta_{2})<0$. This shows that $D_{L}(-f)$ contains
$\{\beta_{1},\beta_{2}\}$ and $D_{L}(-f)\cap\hat{P}=\emptyset$. Write
$f=a+bz,a,b\in F$, and  denote by $N(f)=a^{2}-b^{2}s$  the norm of
$f$ from $L$ down to $F$. Then $D_{L}(-f)=D_{L}(N(f),-a)\cup
D_{L}(-N(f),-bz)$.
Then $N(f)(\alpha)>0$ because $-bz$ is positive in exactly one
ordering in $\{\beta_1, \beta_2\}$ but
$-f \in \beta_1 \cap \beta_2$. Hence we also have
$N(f)(\beta_{1})>0,N(f)(\beta_{2})>0$. Thus we have
$\{\beta_{1},\beta_{2}\}\subseteq
D_{L}(N(f),-a)$ and also $\alpha\in D_{F}(N(f),-a)$. By the
preceding lemma and because $\hat{T}$ is open, there is a $u\in F$ such
that $\emptyset\neq
D_{F}(-u)\subseteq D_{F}(N(f),-a)\cap\hat{T}$.
Hence $D_{L}(-u) \subset D_L(-f) \subset
X(L)\setminus\hat{P}$.
Thus $u>0$ on $\hat{P}$ and $u\in P\cap F=T$. But since
$D(-u)\cap\hat{T}\neq\emptyset$, this is a contradiction, which proves
(1).

(2) We show that $\pi$ is also injective. Let $\alpha\in X(F)$.  If
$\beta_{1},\beta_{2}$ are
the two points in $\pi^{-1}\{\alpha\}\subset X(L)$,
then
$z$ has opposite signs at $\beta_{1}$ and $\beta_{2}$.
But we also have $zf\in P$ for some $f\in F$, and if
$\beta\in\hat{P}$, then $zf\in\beta$. Since  $f$ has the same sign at
$\beta_{1}$ and  $\beta_{2}$ (given by the sign at $\alpha$),
   $z$ would have the
same sign at $\beta_{1}$ and $\beta_{2}$ if both were in $\hat{P}$: a
contradiction. This shows that only
one of the $\beta_{i}$'s can be in $\hat{P}$ and that $\pi$ is
injective on $\hat{P}$.
\qed
\enddemo

\definition{Remark 4.8}
Concerning Lemma 4.7 (2), we should point out that
one can construct a field $F$ with a preordering $T$ and
a quadratic extension $L$ with a preordering $P$ such that:
\roster
\item $P \cap F = T$,
\item There is a $1-1$ correspondence between $C_2$-orderings of $L$
containing $P$ and $C_2$-orderings in $F$ containing
$T$ induced by the restriction map,
\item The natural map $\dot{F}/T \lra \dot{L}/P$ is not surjective.
\endroster

One possible example is $F=\Bbb{R}(X,Y)$, $T$ the set of nonzero sums
of squares in $F$ and $L=F(\sqrt{s})$, where $s=1+X^2$.  Then set
$\hat{P}=\{\alpha\in X(F)\mid \sqrt{s}\in\alpha\}$, and
$P=\bigcap_{\alpha\in\hat{P}}\alpha$.  We claim that $P$
satisfies the conditions (1), (2) and (3) above.  (1) is valid because for
each ordering $\gamma\in X_F$, there is an ordering $\alpha\in X_L$
such that $\alpha\cap F=\gamma$.  (2) is also valid because for each
ordering $\gamma\in X_F$ there is exactly one ordering $\alpha\in X_L$
such that $\alpha\cap F=\gamma$.  Finally using the proof of
Proposition~4.10 below (which does not utilize Remark~4.8), we
show that (3) is valid as well.  Suppose to the contrary that the natural
map $\dot{F}/T \lra \dot{L}/P$ is surjective. Then by (1) it is
an isomorphism.  In the proof of Proposition~4.10 we show that
there is no such $P$ in $L$.  \qed
\enddefinition

\proclaim{Proposition 4.9} Let $F$ be the function field of an
algebraic variety over a real closed field $R$. Let $T$ be a
preordering such that $\hat{T}$ is open. Let $L=F[Z]/(Z^{2}-s)$ for $s\in
T\setminus F^{2}$.
Then the following conditions are equivalent:

(1) $(F,T)$ admits a $T$-extension $(L,P)$ with $P$ a preordering,

(2) There exists an element $f\in \dot{F}$ such that for all $a,b\in
F$ there exists $g\in \dot{F}$ for which
$[D_{F}(a^{2}-b^{2}s,a)\cup
D_{F}(-a^{2}+b^{2}s,bf)]\cap\hat{T}=D_{F}(g)\cap\hat{T}$.
\endproclaim

\demo{Proof} Denote by $N$ the norm from $L$ down to $F$, by $z$ the
class of $Z$ in $L$ and by $\pi\colon X(L)\lra X(F)$ the restriction
map.

Let us prove that condition (1) implies condition (2). Assume $(L,
P)$ is a $T$-extension of $(F,T)$ with $P$ a preordering of $L$.  Then
$\dot F/T\cong \dot L/P$ and for any $h\in \dot L$ there is a $g\in
\dot F$ such that $gh\in P$.  Let $f$ be an element of $\dot F$
such that $fz\in P$.

We show first that $\pi(D_{L}(h)\cap\hat{P})=D_{F}(g)\cap\hat{T}$.
By the preceding lemma we know that $\pi$ induces a bijection from
$\hat{P}$ onto $\hat{T}$.  For $\alpha\in \hat{T}$, let $\beta$ be the
unique element of $\hat{P}$ such that $\pi(\beta)=\alpha$.  Then we
have $\alpha\in D_{F}(g)\cap\hat{T}$ if and only if $g(\alpha)>0$ and $
\alpha\in\hat{T}$, if and only if $g(\beta)>0$ and $ \beta\in\hat{P}$,
because
$g\in F\subset L$.  Since $gh\in P\subset\beta$, this is also
equivalent to $h(\beta)>0$ and $\beta\in\hat{P}$, i.e.  $\beta\in
D_{L}(h)\cap\hat{P}$.  Thus we have proved that
$\pi(D_{L}(h)\cap\hat{P})=D_{F}(g)\cap\hat{T}$. In particular, if
$b\in F$, then
$\pi(D_{L}(bz)\cap\hat{P})=D_{F}(bf)\cap\hat{T}$.

On the other hand, for $h=a+bz\in L, a,b\in F$, one has
$D_{L}(h)=D_{L}(N(h),a)\cup D_{L}(-N(h),bz)$.  Since $\pi$ is a
bijection between $\hat{P}$ and $\hat{T}$, it preserves subset
intersections (and of course unions). Then
$\pi(D_{L}(h)\cap\hat{P})=(\pi(D_{L}(N(h),a))\cap\hat{T})\cup
\pi(D_{L}(-N(h),bz)\cap\hat{P})$. The first set of this
union is $D_{F}(N(h),a)\cap\hat{T}$, and the second is
$D_{F}(-N(h))\cap\pi(D_{L}(bz)\cap\hat{P})=D_{F}(-N(h),bf)\cap\hat{T}$.
This proves $(1)$ implies $(2)$.

Conversely, we show that condition (2) implies condition (1).
Suppose there is an $f$ satisfying  condition (2).
Define $S:=D_{L}(fz)\cap\pi^{-1}(\hat{T})$ and set
$P:=\bigcap_{\beta\in S}\beta$. We want to show that $(L,P)$ is a
$T$-extension of $(F,T)$. An element $a\in F$ is in $P$ if and only if
$a(\beta)>0$ for $\beta\in S$, if and only if $a(\alpha)>0$ for
$\alpha\in\pi(S)$. Since $\hat{T}\subseteq\pi(D_{L}(fz))$, we have $P\cap
F=T$.

Let $h=a+bz\in\dot{L}$ with $a, b\in F$. By our assumption there exists
$g\in\dot{F}$ such that
$[D_F(N(h),a)\cup D_F(-N(h),bf)]\cap\hat{T}=D_F(g)\cap\hat{T}$. We claim
that $gh\in P$.

For the sake of simplicity write $V=D_F(N(h),a)\cup D_F(-N(h),bf)$. We have
$D_L(h)=D_L(N(h),a)\cup D_L(-N(h),bz)$. Hence $D_L(h)\cap S=\pi^{-1}(V)\cap
S=(\pi^{-1}(D_F(g)\cap\hat{T}))\cap S=
D_L(g)\cap S$. Therefore $gh\in P$ as required. This shows that $(L,P)$ is
a $T$-extension of $(F,T)$, and the proof of the proposition is complete.
\qed
\enddemo

\remark{Remark} Observe that in the proof that (2) implies (1) in
Proposition~4.9, we did not use all hypotheses
in this proposition. In fact we proved:
\endremark

{\it{Let $F$ be a formally real field and $T$ be a preordering of $F$. Let
$L=F(\sqrt{s})$ for $s\in T \setminus F^2$.
Suppose that there exists an element $f \in \dot{F}$ such that for all
$a,b\in F$ there exists $g\in \dot{F}$ for
which $D_F(a^2-b^2s,a) \cup D_F(-a^2+b^2s,bf) \cap \hat{T}=D_F(g)\cap
\hat{T}$. Then $(F,T)$ admits a $T$-extension
$(L,P)$ with $P$ a preordering.}}

\bigskip

As an application of the material discussed above we have the
following illustration.

\proclaim{Proposition 4.10} Let $F=\R(X,Y)$ and let $T$ be the set of
nonzero sums of squares in $F$. If $H$ is a subgroup of $\wg$ such that
$T=P_{H}$, then the
$H$-ordered field $(F,T)$ does not admit
an  $H$-closure.
\endproclaim

\demo{Proof} The hypotheses of Proposition~4.9 are obviously
satisfied, because $\hat{T}$ is the whole space.  Assume $(K,S)$ is an
$H$-closure of $(F,T)$.  Let $s\in T\setminus F^2$,
$L=F(\sqrt{s})=F[Z]/(Z^2-s)$,
and let $P=L\cap\dot K^2$.
Then by Lemma~4.4, $(L, P)$ is a $T$-extension of $(F,T)$ with $P$ a
preordering of $L$.  By Proposition 4.9, there exists an $f\in F$ such
that for every $u,v\in F$ the open sets $D_{F}(u^{2}-v^{2}s,u)\cup
D_{F}(-u^{2}+v^{2}s,vf)$ are principal.  We show that this is not
true for $s=1+X^{2}$.

Take $h=Y+c+bz\in L$ with $c,b\in\R, b>0$.
Assume that the corresponding set $D_{F}(N(h),Y+c)\cup D_{F}(-N(h),f)$
is the principal set $D_{F}(g)$ for a given square-free polynomial
$g\in F$.  Note that the equation $N(h)=0$ in $\R^{2}$ defines the
hyperbola ${\Cal H}$ of equation $(Y+c)^2=b^2(1+X^2)$.  Set
$A :=\{(X,Y)\in\R^2\,\mid\, N(h)>0,\quad Y+c>0\}$ (respectively
$B :=\{(X,Y)\in\R^2\,\mid\, N(h)>0, \quad Y+c<0\}$) the open region of
the plane above
(respectively  below) the upper (respectively  lower) branch of ${\Cal H}$.  By
assumption, we know that $g> 0$ on $\tilde{A}\cap X(F)=D_{F}(N(h),Y+c)$ and
$g< 0$ on
$\tilde{B}\cap X(F)=D_{F}(N(h),-(Y+c))$.  This implies that $g\ge 0$ on $A$
and $g\le
0$ on $B$ (see \cite{BCR}, \S 7.6) and that
$A$ and $B$ are separated by a branch (i.e. a $1$-dimensional
irreducible connected component) of $g=0$.
Moreover, no branch of $g=0$ can go inside $A\cup B$,
or else $g$ would change sign on $A$ or $B$. (This is due to
the fact that $g$ is square free, and thus every branch is a
sign-changing branch).  Set $C:=\R^2\setminus A\cup B$.
Then $\tilde{C}\cap X(F)=D_{F}(-N(h))$.  Since
$D_{F}(g,-N(h))=D_{F}(bf,-N(h))=D_{F}(f,-N(h))$, we know that $f$ and $g$
have the
same sign on $C$, up to a  $0$-dimensional set.  Thus $f=0$ must also have
a sign-changing branch
contained in $C$, and since $f$ may be chosen square free, any branch
of $f=0$ having a nonempty intersection with the interior of $C$ must
be contained in $C$.

Suppose this is true at the same time for $h=h_{1}=Y+z$ and
$h=h_{2}=Y+4+2z$.  Then
   \roster
\item no branch of $f=0$ is allowed to cross a branch of the
hyperbolas ${\Cal H}_{i}, \, i=1,2$, and
\item there is a branch of $f=0$ splitting the plane into two
connected components, each of them containing one branch of ${\Cal H}_{i}$.
   \endroster
As the upper branch of ${\Cal H}_{2}$ crosses the two branches of
${\Cal H}_{1}$, this is impossible.  This provides a contradiction to the
existence of an $H$-closure for $T$, finishing the proof of 
Proposition~4.10.
\qed
\enddemo

\definition{Remark 4.11} Associated to the group $\dot F/T$ of the
preceding proposition is the ``abstract Witt ring'' of
$T$-forms (see \cite{Ma}), which is actually the reduced Witt ring
$W_{red}(F)$. (See also \cite{L2, Chapter 1} for the definition of
$W_{red}(F)$.)
Proposition~4.10 shows there is no extension
$F\lra K$ such that $W_{red}(F)$ becomes isomorphic to $W(K)$.

This can be viewed as a weak version of the ``unrealizability'' of
$W_{red}(F)$ as a ``true'' Witt ring (See \cite{Ma}, as well as
\cite{Cr2}, and the
remarks on Craven's results below). Note that $W_{red}(F)$ might
actually be isomorphic to $W(K)$ for some field $K$ not related to
$F$, as shown in Example~8.14. We shall now make these remarks more precise.
\enddefinition

\proclaim{Proposition 4.12} Let $F = \Bbb{R}(X,Y)$. Then there is
no field extension $F\lra K$ with $W(K)\cong W_{red}(K)$ such that the
induced map $W_{red}(F)
\lra W_{red}(K) \cong W(K)$ is an isomorphism.
\endproclaim

\demo{Proof} Suppose on the contrary that there exists a field
extension $K/F$ such that 
the inclusion $F\lra K$
induces an isomorphism $W_{red}(F)\cong W_{red}(K) \cong W(K)$.

Because $W_{red}(F)$ is a
torsion-free ring and $W_{red}(F) \cong W(K)$ we see that $W(K)$ is
torsion-free as well.
Thus $K$ is a  pythagorean field.  (See \cite{L1, Chapter 8}.)
\comment
(Recall that a field $K$ is called pythagorean if each
sum of two squares in $K$ is again a square in
$K$.)

Indeed if $k = k_1^2 + k_2^2 \in \dot{K}$ we see that $2 \<1,-k\>
\cong \<1,1,-k,-k\>$
is a four-dimensional
isotropic Pfister form over $K$. Hence this form is hyperbolic, and
therefore it represents the zero element in
$W(K)$. (See \cite{L1, Chapters 2 and 10}.)
Because $W(K)$ is torsion-free, it follows that $\<1,-k\> = 0$ in $W(K)$,
and hence $k$ is a
square in $\dot{K}$. This means  that $K$ is indeed a pythagorean
field. 
\endcomment
Observe also that $-1 \notin \dot{K}^2$ because
otherwise $K$ would be a
quadratically closed field and $W_{red}(F)$ would not be
  isomorphic to $W(K)$. Hence $K^2$ is
a preordering in $K$. Set $T$ to be the set of nonzero sums of 
squares in $F$. It is well known that the group of units in
$W_{red}(F)$ is $\{fT\colon f \in \dot{F}\} =
\dot{F}/T$ and the group of units of $W(K)$ is
$\dot{K}/\dot{K}^2$ (because $W(K)$ is reduced). (See
\cite{L2, Proposition 1.24}.) Therefore the isomorphism $W_{red}(F)
\cong W(K)$ induces
an isomorphism $\dot{F}/T \cong
\dot{K}/\dot{K}^2$.

Now let $F(2)$ be a quadratic closure of $F$ and set $L = K \cap
F(2)$. (We assume that both fields $K$ and
$F(2)$ lie in some fixed field extension of $F$.) Let $l_1,l_2 \in L$
and $l_1^2 + l_2^2 \in \dot{L}$. Then because $K$ is
a pythagorean field we see that there exists an element $k \in
\dot{K}$ such that $k^2 = l_1^2 + l_2^2$. Since $k$ also
belongs to $F(2)$ we see that $k \in L$ and $L$ is a pythagorean
field. (Observe that in general any intersection of
pythagorean fields is a pythagorean field.)
We also see that $\dot{L}^2 \cap F = \dot{K}^2 \cap F = T$,
because for each $t \in T, \sqrt{t} \in F(2)$.

 Finally we claim that
the natural homomorphism
$\varphi\colon\dot{F}/T \lra \dot{L}/\dot{L}^2$ is in fact
an isomorphism. Because $\dot{L}^2 \cap F=T$,
we see that
$\varphi$ is injective.
Consider now an element $l \in \dot{L}$. Because the natural map
$\dot{F}/T \lra \dot{K}/\dot{K}^2$ is surjective, we see
that there exist elements $f \in \dot{F}$ and $k \in \dot K$ such that
$lf^{-1} = k^2
\in
\dot{K}^2$. Because $lf^{-1} \in \dot{L} \subset
F(2)$ we see that $k \in F(2) \cap K=L$. Therefore  the map $\dot{F}/T
\lra \dot{L}/\dot{L}^2$ is surjective.

 From the proof of Proposition 4.10, we see that there is no field
extension $L/F, L
\subset F(2)$, such that $\dot{L}^2$ is additively closed, $-1\notin
\dot{L}^2$,
and the
natural homomorphism
$\dot{F}/T \lra \dot{L}/\dot{L}^2$ is an isomorphism. Thus we have
arrived at a
contradiction.
\qed
\enddemo

T. Craven kindly called our attention to \cite{Cr2, Theorem 5.5},
which can be applied in the construction of formally real  fields $F$ such that
$W_{red}(F)$ is not isomorphic to $W(K)$ for any field extension
under the natural map
induced by the inclusion $F\lra K$. (As observed in the proof of
Proposition 4.12 above,
if we want
$W_{red}(F)\cong W(K)$ then $K$ must  be a formally real pythagorean
field, and the inclusion $F \lra K$ induces a  natural
homomorphism $W_{red}(F)\lra W_{red}(K) = W(K)$.) The following
proposition we attribute to
T.~Craven, as it is an  immediate corollary of
\cite{Cr2, Theorem 5.5}.

\proclaim{Proposition 4.13 (Craven)} Let $F=L(X)$ where $L$ is a
formally real field, which is not a pythagorean field.
Then for each pythagorean field extension $K/F$, the natural
homomorphism $W_{red}(F)\lra W_{red}(K) = W(K)$ induced by the inclusion
map
$F\lra K$ is not an isomorphism.
\endproclaim

\demo{Proof} Assume that $K$ is a pythagorean field extension of
$F=L(X)$, where $L$ is a formally real field which is not
pythagorean, and suppose that the field extension
$F\lra K$ induces an isomorphism $W_{red}(F)\lra
W_{red}(K)$.

Because $L$ is not a pythagorean field, there exists an element
$l=l_1^2+l_2^2, l_1, l_2 \in L$ such that $l \notin
\dot{L}^2$. Because $K$ is a pythagorean field, there exists an
element $k \in \dot{K}$ such that $k^2=l$. Hence the
polynomial $f(X)=X^2-l$ has a root in $K$. Then from
\cite{Cr2, Theorem 5.5(b)}, we see that
$f(X)$ has exactly one root in every real closure of $L$. Of course
this is not true, as each real closure of $L$ must
contain both roots of $f(X)$. Hence we have arrived at a
contradiction, completing the proof. \qed
\enddemo

\remark{Remark} We can say that $W_{red}(F)$
is not realizable as the Witt ring of an extension $K/F$.
\endremark
\smallskip
In the other direction we present a case below, where $(F,T)$ admits a
maximal preordered $T$-extension
$(\dot{K},\dot{K}^2)$. We recall that a preordering $T$ in $F$ is SAP
(Strong Approximation Property)
if and only if for each set of
elements $a_1,\ldots,a_n \in \dot{F}$ there exists an element $a\in
\dot{F}$ such that
$D_F(a_1,\ldots,a_n)\cap\hat{T}=D_F(a)\cap\hat{T}$. (Here as above,
$\hat{T}$ is the set of all orderings $\alpha\in F$
such that $T \subset \alpha$.) If $T$ is SAP and $R$ is a
preordering of $F$ containing $T$, then
$R$ is SAP as well. (See \cite{L2, Theorem~17.12 and Corollary~16.8}.) 
The definition of SAP implies that $(D_F(a) \cup D_F(b)) \cap \hat T =
D_F(c) \cap \hat T$ for some $c\in \dot F$. Thus condition (2)  of
Proposition 4.9 holds (and hence also condition (1), by the remark
following Proposition 4.9).

\proclaim{Proposition 4.14} Let $F$ be a formally real field, and
let $T$ be a SAP preordering in $F$. Then $(F,T)$ admits a maximal
preordered $T$-extension
$(K,\dot K^{2})$.
\endproclaim

\demo{Proof} Let $F$ be a formally real field and let $T$ be a SAP
preordering in $F$. Using Zorn's lemma we see that
there exists a $T$-extension $(L,S)$ of $(F,T)$ which is maximal among
the preordered $T$-extensions. We claim that
$S$ is a SAP preordering in $L$. In order to show this, pick any elements
$a_1,\ldots,a_n\in\dot{L}$. Because $(L,S)$ is a
$T$-extension of $(F,T)$ we see that there exist elements
$b_1,b_2,\ldots,b_n\in \dot{F}$ such that $b_i a_i \in S$ for
each $i=1,2,\ldots,n$. Because $T$ is  SAP there exists an element
$b\in\dot{F}$ such that
$D_F(b_1,\ldots,b_n)\cap\hat{T}=D_F(b)\cap\hat{T}$.

We have $D_L(b)\cap\hat{S}=D_L(b_1,\ldots,b_n)\cap\hat{S}$. Indeed let
$\alpha\in D_L(b)\cap\hat{S}$. Then
$b\in\alpha$ and $\alpha\cap F\in\hat{T}$. Therefore
$b_1,\ldots,b_n\in\alpha$ and $\alpha\in D_L(b_1,\ldots,b_n)\cap\hat{S}$.
Assume now that $\alpha\in D_L(b_1,\ldots,b_n)\cap\hat{S}$. Then $b_1,\ldots,b_n\in\alpha$ and $\alpha\cap F\in\hat{T}$.
Hence $b\in\alpha$ and $\alpha\in D_L(b)\cap\hat{S}$.

Finally observe that since $b_i a_i\in S$ for all $i=1,\ldots,n$ we have
$D_L(b_1,\ldots,b_n)=D_L(a_1,\ldots,a_n)$.
Therefore $D_L(b)\cap\hat{S}=D_L(a_1,\ldots,a_n)\cap\hat{S}$ as required.

Now we claim that $S=\dot{L}^2$. Suppose that this is not true. Then there
exists an element $s\in S \setminus \dot{L}^2$
and we can set $E=L(\sqrt{s})$. By the remark following Proposition~4.9 we
see that one can find a preordering $R$ in
$E$ such that $(E,R)$ is a $T$-extension of $(F,T)$. This is a
contradiction with the fact that $(L,S)$ is a maximal
$T$-extension of $(F,T)$ such that $S$ is a preordering in $L$. Therefore
we can set $L=\dot{K}$ and $S=\dot{K}^2$
to complete the proof.
\qed
\enddemo

The preceding proposition will apply in particular when $F$ is a formally real
field of transcendence degree $1$ over a real closed field, because
those fields are known to have stability index $1$, which implies
Strong Approximation Property (\cite{L2, Corollary~17.11}).

\definition{Remark 4.15} Note that SAP is not a necessary condition
for
the existence of closures for preorderings. If $F$ is pythagorean,
then it is its own closure with respect to its minimal preordering. But
there are pythagorean fields which are not SAP: for example the
field of iterated power series $\R((X))((Y))$. ( See also~\cite{Cr1}
   for  more examples.)
\enddefinition

\head
\S 5.  Cyclic subgroups of $W$-groups
\endhead

In this section we consider the subgroups $H$ of $\wg$ which are the
easiest to understand in terms of their associated $H$-orderings,
namely the two cyclic groups $C_2$ and $C_4$.  As mentioned earlier,
$C_2$ in many ways is the motivating example for this entire theory,
and we cite here the results previously given in \cite{MiSp1} for this
group, as a means of illustrating the results we are attempting to
generalize in this paper.  As any single element of $\wg$ necessarily
generates a cyclic subgroup of order 2 or 4, those which generate
subgroups of order 4 are precisely those not associated with usual
orderings on the field $F$.  These are the so-called half-orders of
$F$, as investigated in \cite{K1}; this concept was first introduced
by Sperner \cite{S} in 1949, in a geometrical context.

\definition{Definition 5.1} A {\it nonsimple involution} of $\wg$ is an element
$\sigma \in \wg$ such that $\sigma ^2 = 1$ and $\sigma \notin \Phi (\wg )$.  In
other words, a nonsimple involution is an element of $\wg$ which generates an
essential subgroup of order 2.
\enddefinition

\proclaim{Theorem 5.2} \cite{MiSp1} The field $F$ is formally real if and
only if $\wg$ contains a nonsimple involution.  There is a one-one
correspondence between orderings on $F$ and nontrivial cosets of $\Phi (\wg )$
which have an involution as a coset representative.
\endproclaim

We have the well-known characterization of those subgroups of $\dot F$ that are
orderings, which we include here for the sake of completeness.

\proclaim{Proposition 5.3}  A subgroup $S$ of $\dot F$ containing $\dot F^2$ is
a
$C_2$-ordering of
$F$ if and only if the following conditions hold.
\roster
\item $|\dot F/S| = 2$ and
\item $1 + s\in S\ \forall s \in S$.
\endroster
\endproclaim

We can now characterize those subgroups $S$ of $\dot F$ which
are $C_4$-orderings. They are precisely those subgroups of index $2$ which fail
to be orderings.  We also see that $C_4$-ordered fields always admit a
closure.

\proclaim{Proposition 5.4} A subgroup $S$ of $\dot F$ containing $\dot F^2$
is a
$C_4$-ordering of
$F$ if and only if the following conditions hold.
\roster
\item $|\dot F/S| = 2$ and
\item $\exists s \in S$ such that $1 + s \notin S$.
\endroster
\endproclaim
\demo{Proof} We know $S$ is a $C_4$-ordering of $F$ if and only if there
exists $
\sigma \in \wg$ such that $S = \{ a \in \dot F | \sqrt{a}^\sigma = \sqrt{a}\}$
where $\sigma ^2 \neq 1$. Now any subgroup of index $2$ in $\dot F$ is of the
form $\{ a \in \dot F | \sqrt{a}^\sigma = \sqrt{a}\}$ for some $\sigma \in
\wg$, so we need only guarantee that $S$ is not an ordering, which condition
(2) does.
\qed
\enddemo

\definition{Remark 5.5}\roster
\item Note that it is easy to see that condition (2) above can be replaced
by (2') $S+S=\dot{F}$.

\item There are actually two kinds of $C_{4}$-orderings, distinguished by
whether or not they contain $-1$. If $S$ is a $C_{4}$-ordering
such that $-1\in S$, we  say that $S$ has level $1$. The
prototype is given by $\F_{p}^2$ when $p \equiv 1\bmod 4$. If $-1\not\in S$,
then necessarily $-1\in S+S$, and we  say that $S$ has level $2$.
The model is $\F_{p}^2$ when $p \equiv -1\bmod 4$. It is clear that every
$C_{4}$-extension preserves the level.
\endroster
\enddefinition

\proclaim{Proposition 5.6}
   Let $(K,\dot K^2)$ be a maximal $T$-extension of a $C_4$-ordered
field $(F,T)$. Then \roster
\item $K$ is characterized by the condition of being maximal in $F(2)$
among fields $L
\supseteq F$ such that $\sqrt{a} \notin  L\ \forall a \in \dot F \backslash T$.
\item $\Cal G_K \cong C_4$.
\item $\Gal(K(2)/K) \cong \zz _2$, the group of $2$-adic integers.
\endroster
In particular, every maximal $T$-extension of a $C_{4}$-ordered
field $(F,T)$ is a $C_{4}$-closure, and thus $C_{4}$-closures always
exist.
\endproclaim
\demo{Proof} Let $(K, \dot K^2)$ be a maximal $T$-extension of the
$C_4$-ordered
field
$(F,T)$.  Since $\dot K^2 \cap F = T$,  we see that for any $a \in \dot F
\backslash T$, we have $\sqrt{a} \notin K$, while for any $a \in T$, we have
$\sqrt{a}
\in K$.  Now if $L \supsetneq K$ in $F(2)$, then $L \supseteq K(\sqrt{a})$ for
some $a \in \dot K\backslash\dot K^2$.  Since the cosets of $\dot K^2$ in $\dot
K$ correspond naturally to the cosets of $T$ in $\dot F$, we see that $L$
contains
$\sqrt{a'}$ for some $a' \in \dot{F}\setminus T$, and thus $K$ is maximal
among such extensions
of $F$ in $F(2)$. Conversely, suppose $K$ is maximal in $F(2)$ among fields $L
\supseteq F$ such that $\sqrt{a} \notin  L \ \forall a \in F \backslash T$.
Then we see that $\dot K^2
\cap F = T$. We need to see that $|\dot K /\dot K^2| = 2$.  Suppose
it is not true.
    Fix $a \in \dot F \backslash T$, so that $a \notin \dot K^2$, and suppose
there exists some $b \in \dot K$ such that $a, b$ are linearly independent in
$\dot K/\dot K^2$.  Then certainly $b \notin aT$, and setting $L = K(\sqrt{b})$
contradicts the maximality of
$K$.  Thus we have that $(K, \dot K^2)$ is a maximal $T$-extension for
$(F,T)$, and this proves (1).

Now observe that $\Cal G_K$ is generated by one generator, since $|\dot
K/\dot K^2| = 2$, so $\Cal G_K \cong C_2$ or $C_4$.  It cannot be $C_2$, or
else $T$ would be an ordering on $F$.  Thus $\Cal G_K \cong C_4$.  Finally,
$\Gal(K(2)/K)$  is cyclic and cannot be finite, since it is not $C_2$
(see \cite{Be1}). Thus $\Gal(K(2)/K) \cong \zz _2$.
\qed
\enddemo

\head \S 6.  Subgroups of $W$-groups generated by two elements
\endhead

As we saw in Theorem 2.5, a  group generated by two elements appearing as
a subgroup of $\wg$ may only be one in the list $C_2*C_4,
C_4*C_4,C_{2}*C_{2}, C_4\times C_4, C_4\rtimes C_4$. The last two
   are particular cases of the groups studied in \S~7 and
\S~8, and we will focus in this section on the first three. The third one is
better known as the dihedral group $D$.

We will give an algebraic characterization for the orderings
associated with these groups and show
that it is always possible to make closures.  Portions of the proofs
rely on the characterizations of $C_4 \times C_4$- and $C_4 \rtimes
C_4$-orderings
obtained in \S~7; but since the results in \S~7 do not rely on those in
\S~6, we freely use these results where needed.

\proclaim{Lemma 6.1} Let $T$ be a subgroup of $\dot F$ such that $\dot F^2
\subseteq T$ and $|\dot F/T| = 4$. If $-1 \notin T$, then $F$ is $(T \cup
-T)$-rigid.
\endproclaim
\demo{Proof}  Let $\dot F/T = \{ 1, -1, a, -a\}$.  Then $(T \cup -T) + a(T \cup
-T) \subseteq (T \cup -T) \cup a(T\cup -T) = \dot F$.
\qed \enddemo

\proclaim{Proposition 6.2}  A subgroup $T$ of $\dot F$ is a $C_2*C_4$-ordering
if and only if $\dot F^2 \subseteq T$, $|\dot F/T| = 4$, and the following two
conditions hold.
\roster
\item $T + T \neq T$, and
\item $-1\not\in\sum T$, where $\sum T$ denotes the set of all finite
sums of elements of $T$.
\endroster
\endproclaim

\demo{Proof} The conditions $\dot F^2 \subseteq T$ and $|\dot F/T| = 4$ are
necessary and sufficient for $T$ to be a $G$-ordering for some
essential subgroup $G \subseteq \wg$ generated by two elements 
$\sigma, \tau$,
independent mod $\Phi(\wg )$. We next show the necessity of conditions (1) and
(2).  Let $G \cong C_2*C_4$ be a subgroup of $\wg$, where $T = P_G$. We assume
$G$ is generated by two noncommuting (hence independent mod $\Phi(\wg)$)
elements $\sigma, \tau$ such that $\sigma^2 = 1, \tau^4 = 1$.  If $T +
T = T$, then by Proposition~6.14, $T$ would be a
$D$-ordering (this  is independent of previous results). Since it is not,
we see that (1) holds.  Also $-T \nsubseteq \sum
T$, since $\sum T \subseteq P_\sigma$, which is an ordering because $\sigma$ is
an involution. Thus $P_{\sigma}$ cannot contain $-T$ and condition (2) holds.

We now show the sufficiency of the conditions.  Since $T$ is a $G$-ordering for
some essential subgroup with two generators, it must be isomorphic to one of
the five groups listed in Theorem 2.5. Since $-1 \notin T$ by (2), it cannot be
$C_4\times C_4$ by Proposition 7.2 in the next section.  Also (1) shows that
$G$ cannot be isomorphic to $D \cong C_2*C_2$ by Proposition~6.14, and (2)
shows that $G$ cannot be isomorphic to $C_4 \rtimes C_4$ by Proposition 7.6.
Finally, from (1) and (2) we can see that $\sum T$ is an ordering on $F$, since
it is clearly a proper subgroup of $\dot F$, which properly contains $T$, so
must be of index $2$ in $\dot F$; it does not contain $-1$, and it is closed
under addition.  Then $\sum T = T \cup aT$ for some $a \notin T$, and $G$ is
generated by elements $\sigma, \tau$ where the intersection of the fixed field
of $\sigma$ with $F^{(2)}$ is $K(\sqrt{a})$, and the intersection of the fixed
field of $\tau$ with  $F^{(2)}$ is $K(\sqrt{-1})$.  Then $P_\sigma = \sum T$ is
an ordering, so $\sigma$ is an involution.  This shows $G$ cannot be isomorphic
to $C_4*C_4$.  Thus the only remaining possibility is $G \cong C_2*C_4$.
\qed
\enddemo

\proclaim{Proposition 6.3}  A subgroup $T$ of $\dot F$ is a
$C_4*C_4$-ordering if and only if $\dot F^2 \subseteq T$, $|\dot F/T| = 4$,
and one of the following two conditions hold.
\roster
\item $-1 \in T$ and $F$ is not $T$-rigid, or
   \item $-1 \notin T$, $-1\in\sum T$, but $T + T \neq	 T \cup -T$.
\endroster
\endproclaim
\demo{Proof}If $-1 \in T$, the only possible subgroups $H$ of $\wg$ with two
generators for which $T$ can be an $H$-ordering are $C_4 \times C_4$ and
$C_4*C_4$.  The other three are eliminated by Propositions 6.14, 6.2, and
7.6.
Also, if $-1
\in T$, then
$F$ is $T$-rigid if and only if $T$ is a $C_4 \times
C_4$-ordering by Proposition 7.2.  This leaves $C_4 * C_4$ as the only
possibility.

If $-1 \not\in T$, there are three possibilities to consider:
$-1\not\in\sum T$, $T + T = T \cup - T$, or $-1\in\sum T$ but $T
+ T \neq T \cup -T$.  The first case occurs if and only if $T$ is either a
$D$-ordering (by Proposition 6.14) or a $C_2*C_4$-ordering (by Proposition
6.2).  The second case occurs if and only if $T$ is a $C_4 \rtimes
C_4$-ordering by Proposition 7.6 and Lemma 6.1.  Thus, the third case must
occur if and only
if $T$ is a $C_4*C_4$-ordering as claimed.
\qed
\enddemo

The following example constructs a $C_4 * C_4$-ordering of $\qq _2$. It is
illustrative, in that it shows how even in a relatively \lq\lq small\rq\rq\
setting, the additive structure of $T$ can behave quite differently from the
additive structure of $\dot F(\sqrt{T})^2$.  In particular, it shows that
$\langle 1, 1\rangle$ may represent elements in $F(\sqrt{T})$ which are not in
$T + T$.  In this example, $T + T$ is not multiplicatively closed, but of
course the form $\langle 1,1\rangle$, being a Pfister form, is multiplicative
in $F(\sqrt{T})$.

\example{Example 6.4} In $F = \qq _2$ consider the subgroup $T = \dot F^2 \cup
5\dot F^2$ of the square class group.  Using the notation for $\Cal G_2$ as in
Example 2.9, we see that the corresponding subgroup of $\Cal
G_2$ is $H =
\langle \sigma, \tau\rangle \cong C_4 * C_4$.  This is a W-group associated
with the Witt ring $\zz/4\zz \times_{M} \zz/4\zz$, where the product
``$\times_{M}$" is taken in
the category of Witt rings (see \cite {Ma} and \cite{MiSm2}).
The fixed field of $H$ is $K=\qq_2(\sqrt{5})$.
The form
$\langle 1,1\rangle$
   represents $-1$ over $K$, and this can be shown as follows.
   It is well known and easy to show that for any quadratic field extension
   $F\lra K=F(\sqrt{a})$, one has $( K^{2}+ K^{2})\cap \dot
   F=(F^{2}+F^{2})(F^{2}+aF^{2})$. If $F=\qq_{2}$ and $a=5$, we have
   $30=5\times 6\in (K^{2}+aK^{2})$ and $2\in (K^{2}+K^{2})$. Then
   $15\in K^{2}+ K^{2}$, and since $15$ is congruent to $-1\bmod 16$: it is a
   negative square in $\qq_{2}$. This shows that $-1\in K^{2}+K^{2}$.

   However, when one considers which elements of
$\ff$ are in $T + T$, one finds only the six classes represented by $1, 2, 5,
10, -2, -10$.  In particular, $-1 \notin T + T$, and $T + T$ is not
multiplicatively closed (so forms mod $T$-equivalence do not behave as
quadratic forms over a field behave).  Nonetheless, it is easy
to see
that $-1 \in T + T + T$, so that $T + T \neq T \cup -T$, but $-1
\in \sum T$, consistent with the proposition above.
\endexample

In \S 9 we introduce natural conditions for a subgroup $H$ of
$\wg$ in order to keep track of the additive properties of $\dot F/T$
under $2$-extensions.  We shall see in \S 9 that the group $H \subset
\wg$ above does not possess one of the key properties we require.

\proclaim{Theorem 6.5}  A $(C_2*C_4)$-ordered field $(F,T)$ admits a closure.
\endproclaim
\demo{Proof}  Let $\Cal S$ be the set of extensions $(L,S)$ of $(F,T)$ inside
$F(2)$ satisfying the additional condition that $-1 \notin \sum S$.  As in
the proof of Proposition 4.1, we see that $\Cal S$ has a maximal element
$(K,T_0)$ with $\dot K/T_0 \cong \dot F/T, T = T_0 \cap F$, and $-1 \notin \sum
T_0$.  Then $(K,T_0)$ is a $(C_2*C_4)$-ordered field. To see this we need only
show that conditions (1) and (2) of Proposition 6.2 hold, and condition (2) is
given by construction of $(K,T_0)$.  Condition (1) holds since if $T_0 + T_0 =
T_0$, then $T + T \subseteq (T_0 + T_0) \cap F = T_0 \cap F = T$, contradicting
the fact that $T$ is a $C_2*C_4$-ordering on $F$.

To conclude, we must show $T_0 = \dot K^2$.  Notice $\sum T_0$ is an ordering
on $K$, so $K$ is formally real.  We may write $\dot K/T_0 = \{ \pm T_0, \pm
aT_0\}$, where $a \in T + T$.  If $T_0 \neq \dot K^2$, we can choose $c \in T -
\dot K^2$, and consider $L = K(\sqrt{c})$. Since $-c \notin \sum T_0$, $\sum
T_0$ extends to an ordering $S_0$ on $L$. Then $S_0 \cup -S_0 = \dot L$ and $a
\in S_0$.  Let $S$ be a subgroup of $S_0$ containing $T_0$ and maximal with
respect to excluding $a$.  Then $\dot L/S = \{ \pm S, \pm aS\} \cong \dot K/T_0
\cong \dot F/T$.  Also $S\cap K \supseteq T_0$ by construction, and if there
exists $b \in S \cap K, b \notin T_0$, then $b \in aT_0 \cup -T_0 \cup -aT_0$,
which implies either $a \in S$ or $-1 \in S$, which leads to a contradiction in
either case.  Thus $S \cap K = T_0$, and $(L,S)$ is an extension
contradicting the maximality of $(K,T_0)$.  We conclude $T_0 = \dot K^2$.
\qed
\enddemo

\proclaim{Theorem 6.6} A $(C_4*C_4)$-ordered field $(F,T)$ admits a
$(C_4*C_4)$-closure $(K,\dot K^2)$.
\endproclaim
\demo{Proof} Let $(K,\dot K^2)$ be a maximal $T$-extension for
$(F,T)$.  First assume $-1 \in T$.  We must show $K$ is not a rigid
field.  Let $\{ 1, a, b, ab\}$ be a set of representatives for $\dot
F/T$ which lifts to a set of representatives for $\dot K / \dot K^2$.
Since $F$ is not $T$-rigid, we may, without loss of generality,
assume $b \in T + aT$.  Then $T + aT \subseteq \dot K^2 + a\dot K^2$,
but $b \notin \dot K^2 \cup a\dot K^2$, so $K$ is not rigid, and $\dot
K^2$ is a $(C_4*C_4)$-ordering on $K$.

Now assume $-1 \not\in T=F\cap\dot{K}^2$.  Then 
$-1\not\in\dot{K}^2$, and $-1\in\sum T\subseteq\sum\dot{K}^2$.
Letting $\{ 1, -1, a, -a\}$ be a set of representatives for $\dot
F/T$, this again lifts to a set of representatives for $\dot K / \dot
K^2$.  Since $T + T \neq T \cup -T$, but clearly also $T + T \neq T$,
we may assume $a \in T + T$, so $a \in \dot K^2 + \dot K^2$ as well.
This shows $\dot K^2$ is a $(C_4*C_4)$-ordering on $K$.
\qed
\enddemo
\definition{Remark 6.7} We have defined in Definition~3.2  the level of
an $H$-ordering.  It is then easy to see that the level of a
$(C_{4}*C_{4})$-ordering $T$ is at most $4$.  The level of the closure $K$
(which is the \lq\lq usual" level) is less than or equal to the level
of $T$.  The level of a $(C_4*C_4)$-closure is either $1$ or $2$, as
any field of finite level with at most four square classes has level
at most $2$.  The level of $T$ is $1$ if and only if the level of $K$
is $1$, but in the other cases the level may actually decrease:
Example~6.4 shows that $T$ has level $3$ and that its closure has
level $2$.
\enddefinition

Now we turn our attention to $D$-orderings.  We showed in \S~2 that $C_2 \times
C_2$ cannot be an essential subgroup of $\wg$, so if $H$ is an essential
subgroup of $\wg$ generated by two elements of order 2, necessarily $H \cong
D$. Recall that according to~\cite{Br}, a {\it $2$-element fan}  in $F$ is a
set of two distinct orderings $P_1$, $P_2$ on $F$, and it can be
identified with the preordering $T=P_{1}\cap P_{2}$.

\proclaim{Lemma 6.8} The dihedral group $D$ is a subgroup of $\wg$ if
and only if there is a $2$-element fan in $F$.  In this case,
$T\subseteq\dot{F}$ is a $D$-ordering if and only if $T$ is a $2$-element
fan in $F$.
\endproclaim
\demo{Proof} Let $H = \<\sigma , \tau | \sigma^2 = \tau^2 = [\sigma ,\tau ]^2
= 1\> \cong D$ be a subgroup of $\wg$.  Then $P_\sigma$ and $P_\tau$ are
positive cones of two distinct orderings on $F$, and $P_H = P_\sigma \cap
P_\tau$.  Conversely, if $P_1$, $P_2$ are positive cones corresponding to
distinct orderings on $F$, then there exist nontrivial involutions $\sigma
,\tau \in \wg$, in distinct cosets of $\Phi (\wg )$, such that $P_1 =
P_\sigma$ and $P_2 = P_\tau$.  Then $H = \<\sigma ,\tau\>$ is an essential
subgroup of $\wg$, and $H \cong D$.
\qed
\enddemo

In \cite{BEK}, a field $F$ with two orderings $P_1, P_2$ is defined to be
maximal with respect to $P_1, P_2$ if for any algebraic extension $K$ of
$F$, at
least one of the two orderings cannot be extended to $K$.  Since we prefer to
work inside $F(2)$, we modify this as follows.

\definition{Definition 6.9} A field $F$ with two orderings $P_1, P_2$ is
{\it maximal with respect to} $P_1, P_2$ if for any $2$-extension $K$ of
$F$, at
least one of the orderings does not extend to $K$.
\enddefinition

\proclaim{Proposition 6.10} $(F, P_1, P_2)$ is maximal if and only if
$(F, T_F)$
is a $D$-ordered field, where $T_F = P_1 \cap P_2$, and there exists no proper
$D$-ordered extension field $(L,T_L) \subseteq F(2)$ with $T_L \cap F = T_F$.
\endproclaim
\demo{Proof}
Suppose that the field $(F, P_1, P_2)$ is maximal.  Let $\sigma _1, \sigma _2$
be involutions in $\wg$ such that $P_i = \{ a \in \dot F | \sqrt{a}^{\sigma _i}
= \sqrt{a}\}, i = 1,2$.  Then the subgroup $\< \sigma _1, \sigma _2\> \subseteq
\wg$ is isomorphic to $D$, as we have seen, and $(F, T_F)$ is a $D$-ordered
field as claimed.

Now suppose that $L$ is a $D$-ordered field containing $F$
inside $F(2)$, such that $T_L \cap F = T_F$.  Then $\Cal G_L$ contains a
subgroup isomorphic to $D$, which we can take to be generated by two
involutions $\tau _1, \tau _2$ such that $T_L = Q_1 \cap Q_2$, where $Q_i = \{
a \in \dot L | \sqrt{a}^{\tau _i} = \sqrt{a}\}, i = 1,2$ are distinct orderings
of $L$.  Now $Q_i \cap F \supseteq T_L \cap F = T_F$, so $Q_i \cap F$ is an
ordering of $F$ which contains $T_F, i = 1,2$.  Thus $\{ Q_1 \cap F,
Q_2\cap F\}
= \{ P_1, P_2\}$. Then by maximality of $(F,P_1,P_2)$, we see $L = F$.

Conversely, suppose that $F$ is a $D$-ordered field contained in no proper
$D$-ordered extension field as described.  Then $F$ has at least two distinct
orderings $P_1$ and $P_2$ corresponding to the two involutions generating the
subgroup $D$ of $\wg$, and since there is no proper $D$-ordered extension
field, we see that it is not possible for both orderings to extend to any
extension of $F$.  Thus $(F, P_1, P_2)$ is maximal, as claimed.
\qed
\enddemo

By  Zorn's Lemma  we immediately see the following.

\proclaim{Proposition 6.11} \cite {BEK, Prop.3} Given a field $F$ with two
orderings $P_1, P_2$, there always exists an algebraic extension $\tilde F$ of
$F$ which is maximal with respect to $\tilde P_1, \tilde P_2$, where
$\tilde P_1,
\tilde P_2$ are extensions of $P_1, P_2$ to $\tilde F$.
\endproclaim

\proclaim{Theorem 6.12} A field $(F, P_1, P_2)$ is maximal if and only if
\roster
\item there exist exactly two orderings on $F$ and
\item $F$ is pythagorean, i.e. any sum of squares is a square in $F$.
\endroster
\endproclaim
\demo{Proof} \cite{BEK} Suppose three different orderings $P_1, P_2, P_3$
are possible in $F$. Let $x \in \dot F$ be such that $x$ is positive with
respect to the first two orderings, and negative with respect to $P_3$.  Then
$\sqrt{x} \notin F$, so $F(\sqrt{x})$ is a proper algebraic extension of $F$,
and since $x$ is positive with respect to $P_1$ and $P_2$, they extend to
$F(\sqrt{x})$, and $(F, P_1, P_2)$ cannot be maximal.  Similarly, if $\alpha,
\beta$ are elements of $F$ such that $\sqrt{\alpha^2 + \beta^2} \notin F$, then
$P_1$, $P_2$ can be extended to the proper extension $F(\sqrt{\alpha^2 +
\beta^2})$ of $F$, again contradicting maximality.  Thus conditions (1) and (2)
are necessary.

Conversely, one can show that any field $F$ satisfying conditions (1) and (2)
has $\ff = \{ 1, -1, a, -a\}$ for some $a\in\dot F$.  Now let $F$ be such a
field and let $P_1, P_2$ be the two unique orderings in $F$, so that $a$ is
positive with respect to $P_1$ and negative with respect to $P_2$.  Suppose
$(F,
P_1, P_2)$ were not maximal, and let $K = F(\sqrt{b})$ be a proper quadratic
extension of $F$ such that both $P_1$ and $P_2$ extend to $K$.  Since $K$ is an
ordered proper extension of $F$, $b \neq 1,-1 \in \ff$, so $b = a$ or $-a$.
Then
either $\sqrt{a} \in K$ or $\sqrt{-a} \in K$, so that not both $P_1$ and
$P_2$ extend to $K$. This is a contradiction.
\qed
\enddemo

\proclaim{Corollary 6.13} The $D$-ordered field $(F,T)$ is a maximal
$D$-ordered field if and only if
$\wg \cong D$.  Thus any $D$-ordered field admits a $D$-closure.
\endproclaim
\demo{Proof}  By the preceding theorem, if $F$ is maximal, it has exactly two
orderings, so $\wg$ has exactly two involutions which are independent mod
$\Phi (\wg)$.  Also $F$ is pythagorean, so by \cite{MiSp1}
$\wg$ is generated by involutions.  Thus $\wg$ is generated by two elements
of order 2, and since $\wg$ is necessarily an essential subgroup of itself,
we see that $\wg \cong D$.

Conversely, if $\wg \cong D$, then $F$ is a $D$-ordered field, and since
orderings on $F$ correspond to independent involutions in $\wg$, we see that
$F$ has precisely two distinct orderings.  Also, since $\wg$ is generated by
these involutions, we see that $F$ is pythagorean.  Thus, by the preceding
theorem, $F$ is a maximal $D$-ordered field. Then we see that for any
$D$-ordered field $(L, P_H)$, a maximal $D$-ordered extension $(F, \dot F^2)$
containing $(L, P_H)$ will be a closure for $(L, P_H)$.
\qed
\enddemo

\proclaim{Proposition 6.14} A subgroup $S$ of $\dot F$ containing $\dot F^2$ is
a $D$-ordering of $F$ if and only if  $|\dot F/S| = 4$ and
$1 + s \in S$ whenever $s \in S$.
\endproclaim
\demo{Proof} All that is necessary for $S$ to be a $D$-ordering of $F$ is
that it be a $2$-element fan in $F$. In other words,
$S$ must be a preordering of index 4 in $F$.  A subgroup $S$ of $\dot F$ is
such a preordering if and only if the conditions in the statement of
the proposition are met.
\qed
\enddemo

\head \S 7. Classification of  rigid orderings
\endhead

This section will provide a full Galois-theoretic and algebraic
characterization of all
possible rigid orderings.
We start with the following definition.

\definition{Definition 7.1} Let $I$ be a possibly empty index set.
We call $G$ a
$C(I)$-group if $G$ is isomorphic to $(C_4)^I\times C_4$, an
$S(I)$-group if $G$ is isomorphic to $(C_4)^I\rtimes C_4$, and a
$D(I)$-group if $G$ is isomorphic to $(C_4)^I\rtimes C_2$, the
semidirect product being defined with the nontrivial action of $C_4$
or $C_{2}$ on each inner factor in the last two cases, when $I$ is nonempty.
A $G$-ordering on $F$ is called a $C(I)$- (respectively
$S(I)$-, $D(I)$-) ordering if $G$ is a $C(I)$- (respectively $S(I)$-,
$D(I)$-) group.  When $I=\emptyset$ the $C(I)$- and $S(I)$-orderings
are the $C_{4}$-orderings, and the $D(I)$-orderings are the
$C_{2}$-orderings, that is the usual orderings. Observe that
$C(\emptyset)$- and $S(\emptyset)$-orderings both correspond 
to the same group $C_4$.  The difference between them is that a
$C(\emptyset)$-ordering has level $1$, while an $S(\emptyset)$-ordering has
level $2$. (See Remark 5.5 for comparison.)
When $|I|=1$, we
obtain the groups generated by two elements which are respectively
$(C_4)\times C_4$, $(C_4)\rtimes C_4$ and $D$.
\enddefinition

In this section we will characterize $C(I)$-orderings,
$S(I)$-orderings and $D(I)$-orderings in terms of their algebraic
properties as subgroups of $\dot F$.  We will see in particular that
they are all rigid, and that they constitute the whole class of rigid
orderings.  The group $\coprod_{i \in I}G_i$ will denote the direct sum of
the groups $G_i, i \in I$.

\proclaim{Proposition 7.2} A subgroup $T$ of $\dot F$ containing $\dot F^2$
is a
$C(I)$-ordering if and only if the following three conditions hold.
\roster
\item $-1 \in T$,
\item $F$ is $T$-rigid, and
\item $\dot F /T \cong \coprod_{i\in I\cup\{x\}}(C_2)_i$.
\endroster

In other words, the $C(I)$-orderings are exactly the rigid orderings
of level $1$.
\endproclaim
\demo{Proof} If $I=\emptyset$, the result follows from Proposition~5.4 and
Remark~5.5,
so we shall assume $I\neq\emptyset$.  We begin by showing
that
the three conditions above are necessary.  Let $G \cong C(I)$ and let $T$
be a $G$-ordering.  Suppose $-1 \notin T$.  Let $\{\sigma_i, i \in I;
\sigma_x\}$ generate $G$.  Then $T = \cap_ {i\in
I\cup\{x\}}P_{\sigma_i}$ and $|\dot F/T| \geq 4$.  Thus there are at
least four classes mod $T$, which we can represent as $1, -1, a, -a$ for some
$a \in \dot F$, and there exists a $D^{a,-a}$-extension $L$ of $F$.  Hence
there exist elements $\sigma , \tau \in G$ such that $a \in P_\sigma \backslash
P_\tau$ and $-a \in P_\tau\backslash P_\sigma$. It then follows that
the restriction of $\sigma\tau$ to $L$ has order $4$, so that $\sigma|_L,
\tau|_L$ generate $\Gal(L/F)\cong D$, and hence cannot commute.  Yet $\sigma ,
\tau \in G$, which is an abelian group. This is a contradiction, so $-1 \in T$,
and (1) holds.

Since $-1 \in T$, we have $T \cup -T = T$.  Suppose we have a nonrigid element
$c \in \dot F \backslash T$, so that we have $t_1, t_2 \in T$ with $t_1 + ct_2
\notin T \cup cT$. Then $b = 1 + ct_2/t_1 \notin T \cup cT$. Let $a = -ct_2/t_1
\notin T$. Then $a + b = 1$, so $(\frac{a,b}{F})$ splits. Since $b \notin T
\cup cT = T\cup aT$, $a$ and $b$ are independent mod $T$ and thus mod $\dot
F^2$.  Hence we have a $D^{a,b}$-extension $L$ of $F$, and by the same argument
as above, we find $\sigma , \tau \in G$ which do not commute, leading
to a
contradiction.  Thus $F$ is $T$-rigid and (2) holds.  Finally, by Kummer theory
we know that $\dot F /T $ is isomorphic to the dual $(G/\Phi (G))^* \cong
\coprod_{i\in I\cup\{x\}}(C_2)_i$, giving (3).

We now show that the three conditions are sufficient for $T$ to be a
$C(I)$-ordering. By (3) we see that $T = \cap_{i\in I\cup\{x\}}P_i$ where
$P_i$ is the kernel of the projection $\dot F \to \dot F/T\cong
\coprod_{i\in I\cup\{x\}}(C_2)_i \to (C_2 )_i$.  Further, for each
$P_i$ we have a $\sigma_i \in \wg$ such that $P_i = P_{\sigma_i}$.  Let $G$ be
the closed subgroup of $\wg$ generated by $\{\sigma_i | i\in I\cup\{x\}\}$.
Then $G \subseteq \{\sigma | P_\sigma \supseteq T\}$ because every element of
$G$ must fix every $\sqrt{a}$ left fixed by the $\sigma _i$. So we also have $T
= \cap_{\sigma \in G}P_\sigma$, and $T$ is a $G$-ordering. It remains to show
that $G$ is a $C(I)$-group.

Since $-1 \in T\subseteq P_{\sigma_i}$, none of the $P_{\sigma_i}$ can be
usual
orderings on $F$, so each $\sigma_i$ must have exponent $4$ in $G$.  Since $-1
\in T$ and $F$ is $T$-rigid, we see by Proposition 3.4 that $G$ is abelian.
Then $G$ is a compact abelian group of exponent 4, and $(G/\Phi
(G))^* \cong \coprod_{i\in I\cup\{x\}}(C_2)_i$ is a discrete group of exponent
2.   Then  $((G/\Phi (G))^*)^* \cong G/\Phi (G) \cong \prod_{i\in I\cup\{x\}}
(C_2)_i$, and $G \cong \prod_{i\in I\cup\{x\}}(C_4)_i$, so $G$ is
a $C(I)$-group as claimed.
\qed
\enddemo

In order to characterize the subgroups of $\dot F$ which are $S(I)$-orderings,
we will first prove three lemmas.  Let $G$ be an $S(I)$-group.  It will be
helpful to fix the following notation: write $G = G_1 \rtimes G_2$ where $G_1
\cong \prod_{i \in I} (C_4 )_i$ and $G_2 \cong C_4$.  Let
$\tau$ be a
generator of $G_2$ and $P_\tau = \{ a \in \dot F | \sqrt{a}^\tau = \sqrt{a}\}$.

\proclaim{Lemma 7.3} Let $T$ be a $G$-ordering.  Then $T$ has index $2$ in
$P_{G_1}$.
\endproclaim
\demo{Proof} If $P_{G_1} \subseteq P_{\tau}$, we would have $T = P_{G_1} \cap
P_{\tau} = P_{G_1} = P_G$.  But by Kummer theory and the Burnside Basis
Theorem,
that would imply $G = G_1$. Thus $P_{G_1} \nsubseteq P_\tau$, $T
\subsetneq P_{G_1}$, and $|P_{G_1}/T| \geq 2$.  On the other hand, since $T =
P_{G_1} \cap P_{\tau}$, we have $|P_{G_1}/T| \leq 2$, and so $|P_{G_1}/T| = 2$.
\qed
\enddemo

\proclaim{Lemma 7.4} For any group homomorphism $\theta : G \to C_4 = \langle
\sigma \rangle$, we have $\theta (G_1) \subseteq \< \sigma ^2\>$.
\endproclaim
\demo{Proof} If $a \in G_1$, writing multiplicatively, we have $$\theta
(a^{-1})
= \theta (\tau a\tau^{-1}) = \theta (\tau )\theta (a) \theta (\tau )^{-1} =
\theta (a),$$ so $\theta (a)^2 = 1$.
\qed
\enddemo

\proclaim{Lemma 7.5} We have $T + T \subseteq P_{G_1}$.
\endproclaim
\demo{Proof} Let $a \in T + T, a \notin T$, and consider the
following three cases.

Case 1: $a = x^2 + y^2$.  Then there exists a $C_4^a$-extension $L$ of $F$, and
we have a map $\theta : G \to \Gal(L/F) \cong C_4$, and by Lemma 7.4
$\theta (G_1)$ has order at most $2$.  Thus $\theta (G_1)$ fixes $\sqrt{a}$ and
$a \in P_{G_1}$.

Case 2:  $a = x^2 + t, t \in T\backslash \dot F^2$.  We have $a^2 = ax^2 + at$,
and $a, at$ are independent modulo $\dot F^2$.  Thus there exists a
$D^{a,at}$-extension $L$ of $F$, and $\Gal(L/F(\sqrt{t})) \cong C_4$.  Since $t
\in T$, we have $\sqrt{t}^\sigma = \sqrt{t}$ for $\sigma \in G$, which means we
have a homomorphism $\theta : G \to \Gal(L/F(\sqrt{t})) \cong C_4$. Again
applying Lemma 7.4, $\theta (G_1)$ has order at most $2$, so $G_1$ must fix
$\sqrt{a}$ and $a \in P_{G_1}$.

Case 3:  $a = s + t$, $s, t \in T\backslash \dot F^2$.  We can write $as^{-1} =
1 + ts^{-1}$, and then we are in one of the previous two cases. Hence $as^{-1}
\in P_{G_1}$, and it follows that $a \in P_{G_1}$.
\qed
\enddemo

\proclaim{Proposition 7.6} A subgroup $T$ of $\dot F$ containing $\dot F^2$ is
an $S(I)$-ordering if and only if the following four conditions hold.
\roster
\item $-1 \notin T$,
\item $F$ is $(T\cup -T)$-rigid,
\item $T + T = T \cup -T$, and
\item $\dot F /T \cong \coprod_{i\in I\cup\{x\}}(C_2)_i$.
\endroster
\endproclaim
\demo{Proof} When $I=\emptyset$ the result follows from
Proposition~5.4 and Remark~5.5. Thus we may assume that $I\neq\emptyset$.
We begin by showing the conditions above are necessary. Condition (4)
follows from Kummer theory. Condition (1) follows from Lemma 7.5 above, for if
$-1 \in T$, we would have $\dot F \subseteq \dot F^2 - \dot F^2 \subseteq T -T
= T+T \subseteq P_{G_1}$, but as $|I| \geq 1$, we cannot have $P_{G_1}$ being
all of $\dot F$.

To show the necessity of condition (3), first observe that $-1 \in
P_{G_1}$, $-1 \notin
T$, and $|P_{G_1}/T| = 2$, so $P_{G_1} = T \cup -T$, and thus $T + T \subseteq
T \cup -T$.  To show equality, we need to show that some element of $-T$ is in
$T + T$.  In this case, that amounts to showing that $T$ is not additively
closed.  Suppose that $T$ were additively closed.  Then $T$ would be a
preordering, so contained in some ordering $P_\sigma$ for some $\sigma \in
\wg$.
Further, $\sigma$ is an involution not contained in $\Phi (\wg )$, and $\sigma
\in G = G_1 \rtimes G_2$.  In particular, $\sigma$ is not a square in $G$, and
$\sigma \neq \tau$. Thus $\sigma = \sigma_1\tau$ for some $\sigma_1 \in
G_1$ and
$$\sigma^2 = \sigma_1\tau\sigma_1\tau = \sigma_1\tau\sigma_1\tau^{-1}\tau^2 =
\sigma_1\sigma_1^{-1}\tau^2 = \tau^2 \neq 1.$$
Thus $\sigma$ is not an involution, which is a contradiction, and so $-1
\in T +
T$.  Finally, since $F$ is $P_{G_1}$-rigid and $T \cup -T =
P_{G_1}$, we see that (2) holds.

Now we must show that conditions (1) - (4) are sufficient for $T$ to be an
$S(I)$-ordering. Letting $S = T \cup -T$, we see that $S$ satisfies the
condition for being a $G_1$-ordering, with $G_1 \cong \prod_{i\in
I}(C_4)_i$, as
given in Proposition 7.2.  Let $Q$ be a subgroup of index $2$ in $\dot F$ such
that $T = S \cap Q$, and let $\tau \in \wg$ such that $Q = P_\tau$.  Let $G$ be
the subgroup of $\wg$ generated by $G_1$ and $\tau$.  We need to see that $G =
G_1 \rtimes G_2$ where $G_2$ is the subgroup of $\wg$ generated by $\tau$.
Specifically, we need to show that $G_1 \cap G_2 = \{ 1\}$ and that $[\sigma
,\tau ]\sigma^2 = 1 \ \forall \sigma \in G_1$.

Since $G_1$ fixes $\sqrt{-1}$ and $\tau$ does not, we cannot have $\tau$
or $\tau^{-1}$ in $G_1$.  Suppose $\tau^2 \in G_1$.  Then it has order $2$ in
$G_1$ and hence must be a square. Let $\sigma \in G_1$ such that $\sigma^2 =
\tau^2$.  Since $P_\sigma \neq P_\tau$, there exists $a \in P_\tau
\backslash P_\sigma$, and neither $a$ nor $-a$ can be a square, since
neither is
in $P_\sigma$. Since also $-1 \notin \dot F^2$, we have a $D^{a,-a}$-extension
$L$ of $F$, and $\sigma |_L$ has order $4$ in $\Gal(L/F)$.  However, since
$\tau$ fixes $\sqrt{a}$, $\tau |_L \in \Gal(L/F(\sqrt{a})) \cong C_2 \times
C_2$, and so $\sigma^2 \neq \tau^2$, contradicting the assumption. Thus $G_1
\cap G_2 = \{ 1\}$.

To prove  $[\sigma ,\tau ]\sigma^2 = 1 \ \forall \sigma \in G_1$, it
is sufficient
to show that this condition holds for the restriction of $\sigma ,\tau$ to each
$C_4$- and $D$-extension of $F$.  Suppose $L$ is a $C_4^a$-extension of $F$.
Then $a$ is a sum of two squares, so $a \in T + T = T \cup -T = P_{G_1}$ and
$[\sigma ,\tau ]\sigma^2|_L = \sigma^2|_L$.  Since $\sigma \in G_1$, $\sigma
\in \Gal(L/F(\sqrt{a}))$ and $\sigma^2|_L = 1$.

Now suppose $L$ is a $D^{a,b}$-extension of $F$.  We may assume $\sigma \notin
Z(\Gal(L/F))$ (the centralizer), since otherwise clearly $[\sigma ,\tau
]\sigma^2 |_L = 1$.
Without loss of generality, we may assume $\sqrt{a}^\sigma = -\sqrt{a}$.  Then
$a \notin T \cup -T$, and since $1 = ax^2 + by^2$, we have $b \in T - aT$, and
by rigidity, $b \in T \cup -aT \cup -T \cup aT$.  However, if $b$ were in
$-T$ or $aT$, then we would obtain $a \in T + T = T \cup -T$, a
contradiction.  Thus $b \in T \cup -aT$.

If $b \in T$, then $\sigma$ and $\tau$ both fix
$\sqrt{b}$ and both have order $2$.  If $\tau$ does not fix $\sqrt{a}$, then
$\sigma ,
\tau$ act the same on $\sqrt{a}$ and $\sqrt{b}$ and hence commute.  If $\tau$
fixes $\sqrt{a}$ then $\tau \in Z(\Gal(L/F))$ so in either case $[\sigma
,\tau ]
\sigma^2 = \sigma^2 = 1$.

If $b \in -aT$, then $\sigma$ fixes neither $\sqrt{a}$ nor $\sqrt{-a}$, so has
order $4$.  Since $\tau$ acts differently on $\sqrt{a}$ and $\sqrt{b}$, it must
fix one of them and be of order $2$, and the same holds for $\sigma\tau$.  Then
$[\sigma ,\tau ]\sigma^2 = \sigma\tau\sigma^{-1}\tau^{-1}\sigma^2 =
\tau^{-1}\sigma^{-2}\tau^{-1}\sigma^2 = 1$ since $\sigma^2 \in Z(\Gal(L/F))$.
\qed
\enddemo

We have another convenient formulation of Proposition~7.6 as follows:

\proclaim{Corollary 7.7} A subgroup $T$ of $\dot F$ containing $\dot F^2$ is
an $S(I)$-ordering if and only if the following three conditions hold.

(a) $T$ has level $2$,

(b) $F$ is $T$-rigid, and

(c) $\dot F /T \cong \coprod_{i\in I\cup\{x\}}(C_2)_i$.

\noindent In other words the $S(I)$-orderings are exactly the rigid
orderings of
level $2$.
\endproclaim

\demo{Proof} If $I=\emptyset$, the result follows from
Definition~7.1, so we shall
assume that $I \neq \emptyset$. Assume
that $T$ satisfies (1), (2) and (3) of Proposition~7.6. We show
it is rigid. Let $a\in\dot
F\setminus(T\cup -T)$. Then $T+aT\subset (T\cup -T)+a(T\cup -T)=T\cup
-T\cup aT\cup -aT$. Take $s+at\in T+aT$ and  suppose it is not in
$T\cup aT$.  Then it is in $-T\cup -aT$. If
   $s+at=-u\in -T$ then $-a=t(u+s)\in T+T=T\cup -T$, a contradiction. If
   $s+at=-au\in -aT$ then $-a=s/(u+t)\in T+T=T\cup -T$, a contradiction.
   Thus $T$ is rigid.

    By Proposition~3.3,  a rigid ordering of finite level
greater than $1$ is exactly a rigid ordering of level $2$. This
proves (a) and (b).

Conversely, if $T$ satisfies (a) and (b), then it satisfies (1) and
(3) by Proposition~3.3. Let us show we also have (2).
Let $a\in\dot
F\setminus\pm(T\cup -T)=T\cup -T$. Then $(T\cup -T)+a(T\cup
-T)=\pm(T+aT)\cup\pm(T-aT)\subseteq \pm(T\cup aT)\cup\pm(T\cup -aT)=(T\cup
-T)\cup a(T\cup -T)$.
Since we always have $S\cup aS\subseteq S+aS$ for any subgroup $S$,
we see that $F$ is $T\cup -T$-rigid.\qed
\enddemo

\definition{Example~7.8}It is well-known that if $K\lra L$ is a field
extension and if
$T$ is a usual ordering of $L$, then $S=K\cap T$ is a usual ordering of $K$.
This need not hold for $C(\emptyset)$-orderings nor for
$S(\emptyset)$-orderings.
\comment
Consider the field extension $\Bbb F_{25}/\Bbb F_5$,
where $\Bbb F_5$ and $\Bbb F_{25}$ are fields
containing $5$ and $25$ elements respectively. Then
$\dot{\Bbb F}_{25}^2$ is a $C(\emptyset)$-ordering of $\Bbb F_{25}$. We have
$\dot{\Bbb F}_{25}^2 \cap \Bbb F_5 =
\dot{\Bbb F}_5$. This shows that an intersection of a
$C(\emptyset)$-ordering of a field extension $L/K$ with the base
field $K$ need not be a $C(\emptyset)$-ordering. (Of course in this
example we can replace $\Bbb F_5$ by any $\Bbb F_q$
such that $q$ is a power of an odd prime number $p$.)

In Chapter~9, in Notation~9.9 we consider a quadratic extension $i\colon F_1
\lra F_2$, and $i^{\star} \colon \wgb \lra
\wga$ be the associated restriction map. (See e.g. \cite{MiSm3} for
the existence of
this map. More precisely in
\cite{MiSm3} we observe that $F_1^{(3)} \subset F_2^{(3)}$ and
therefore $i^{\star}$
exists.) Let $H_2$ be a subgroup of
$\wgb$ and let $H_1 = i^{\star}(H_2)$. Then we observe $H_1$ need not
be an essential subgroup of $\wga$.

Indeed in this example we have $F_1 = \Bbb F_5$, $F_2 = \Bbb
F_{25}$ and $H_2 = \wgb \cong C_4$ but $H_1 =
i^{\star} (H_2) = \Phi(\wga)$, and therefore $H_1$ is not
an essential subgroup of $\wga$.

Next we provide an example of a quadratic field extension $K
\lra L$ with an $S(\emptyset)$-ordering $T$ on $L$ such
that $T \cap K$ is a usual ordering. (This means a $C_2$-ordering of $K$.)

Take $K=\R(X),s=1+X^2, L=K(\sqrt{s})$, and denote by $N$ the norm map
from $L$ down to $K$.
Let $O_{+}$ be the ordering $O_{+}:=\{f\in K\,\mid\,
\exists \epsilon>0 \quad f(a)>0$ for $a\in(0,\epsilon)\}$.

Consider the elements $g = X+1+\sqrt{s}$ and
$\bar{g}=X+1-\sqrt{s}$. Then $N(g)=2X$
  and $g+\bar{g}=2(X+1)$. We define an $S(\emptyset)$-ordering $T$ in the
field $L$. Let $T$ be any subgroup of
$\dot{L}$ of index $2$ in $\dot{L}$ such that (a) $T \cap \dot F =
O_{+}$ and (b) $-g,
-\bar{g} \in T$. In order to show that such a subgroup $T$ of
$\dot{L}$ exists, it is
enough to show that the smallest multiplicative subgroup $M$ of $L$ containing
$\dot{L}^2, O_{+}$ and
$\{-g, -\bar{g}\}$ does not contain $-1$. Observe that since
$g\bar{g}=2X \in O_{+} \subset M$, we can replace the set
$\{-g, -\bar{g}\}$ in our consideration above by the singleton $\{-g\}$.

 From the short exact sequence $\dot K/\dot K^2 \lra \dot L/\dot L^2 \overset
N\to\longrightarrow N(L)/\dot K^2 \lra 1$ (see \cite{L1, page 202})
and from the fact
that $N(g)=2X \notin \dot K^2$ we see that $g \notin \dot K \dot L^2$. If
$-1$ were expressible as $-1 = l^2 u(-g)$, where $l \in \dot L$ and
$u \in O_{+}$ then $g$ would belong to
$\dot K \dot L^2$, contradicting the observation above. This proves the
existence of a
  subgroup $T$ of $\dot L$ satisfying (a) and (b) above.
Since $-g -\bar{g} \in T + T$ and since $-g -\bar{g} = -2(X+1) \in
-T$ we see that the level of $T$ is $2$. 
\endcomment
Consider for example $L=K(\sqrt{\dot K})$ and assume that $L$ is
equipped with some $C_{\emptyset}$-ordering $T$.  Since $\dot L^2\cap
K=\dot K$ and $\dot L^2\subseteq T$, we also have $T\cap K=\dot K$:
the $C_{\emptyset}$-ordering $T$ ``vanishes'' under the restriction.  This
happens in particular if $K$ is the finite field $\F_{q}$ with an odd
number $q$ of elements.  With $L=\F_{q^2}$, $\dot L^2$ is a
$C_{\emptyset}$-ordering.  Observe that this cannot happen when $T$ is
an $S(\emptyset)$-ordering in an extension $L$ of $K$: since $-1$ is
not in $T$, it cannot be in $S=T\cap K$, and $S$ cannot be the trivial
index $1$ subgroup.
But $S(\emptyset)$-orderings are subject to another pathology of their
own: it may happen that the restriction of an $S(\emptyset)$-ordering
is a $C_{2}$-ordering. (Observe that this cannot happen with
$C(\emptyset)$-orderings.)
Take for example $K=\qq, L=K(\sqrt{10})$, and denote by $N$ the norm map
from $L$
down to $K$.  Let $\alpha$ be the ordering of $L$ containing
$\sqrt{10}$. Let $v$ be the discrete rank $1$ valuation on $\qq$ associated
to the
prime $3$. Define $T:=\{h\in \dot L\,\mid (-1)^{v(N(h))}h\in \alpha\}$. Then
$-1\not\in T$ and $T$ is a subgroup containing $\dot K^2$, of index $2$ in
$\dot K$ (if $x\not\in T,
-x\in T$) . It is not a usual ordering, since $-4-\sqrt{10}$ is
negative at the two orderings of $L$ but belongs to $T$, as its norm
$6$ has an odd $3$-valuation. Thus it must be an
$S(\emptyset)$-ordering.
Since $N(f)=f^2$ has an even valuation when $f\in K$, we see that $S:=T\cap
K$ is the
usual ordering of $\qq$.
\enddefinition

\proclaim{Proposition~7.9} 	A subgroup $T$ of $\dot F$ containing $\dot
F^2$
is a
$D(I)$-ordering if and only if the following three conditions hold.
\roster
\item $-1 \not\in T$,
\item $T+T=T$,
\item $F$ is $T$-rigid, and
\item $\dot F /T \cong \coprod_{i\in I\cup\{x\}}(C_2)_i$.
\endroster

In other words, $D(I)$-orderings are exactly the rigid orderings of
infinite level.
\endproclaim

\demo{Proof} It is known that conditions (1), (2), (3) are one of the
characterizations of fans~\cite{BK}, and by \cite{CrSm,  Proposition~4.1},
fans are exactly $D(I)$-orderings for some index set $I$ (possibly empty if
we think of
usual orderings as $1$-element fans). By
Proposition~3.3, we immediately see that they are the rigid orderings
of infinite level.
\qed
\enddemo

To conclude the section we may summarize the results with the following

\proclaim{Theorem 7.10} Rigid orderings are exactly $C(I)-,
S(I)-$ or $D(I)$-orderings for some (possibly empty) index set $I$.
\endproclaim

\demo{Proof} This is a straightforward application of Proposition~3.3,
Proposition~7.2, Corollary~7.7 and Proposition~7.9.
\qed
\enddemo

\head \S  8. Construction of closures for  rigid orderings
\endhead

In this section we employ valuation-theoretic techniques to construct
closures for $C(I)$-, $S(I)$- and $D(I)$-orderings.  From the preceding
section, we know that both $C(I)$- and
$S(I)$-orderings are $T$-rigid.  Then for such an ordering
we will be able to use results of Arason, Elman, Jacob \cite {AEJ},
Efrat \cite{Ef} and Ware \cite{Wa} to associate a valuation to $T$.
For $D(I)$-orderings, it is the \lq\lq Fan Trivialization Theorem" of
Br\"{o}cker~\cite{Br, Theorem~2.7} that will be used.  Since it is
well known (see \cite{Ri}) that for each algebraic extension $K/F$ we
can extend any valuation $v$ on $F$ to a valuation $w$ on $K$, we can
then use this to extend $S(I)$- or $D(I)$-orderings, and in most cases also
$C(I)$-orderings, from $F$ to $F(\sqrt{t})$, $t \in T$.  This will allow us to
prove the existence of $S(I)$- and $D(I)$-closures, and in most cases also
$C(I)$-closures.

For the reader's convenience we define here some of the valuation-theoretic
notation we will be using below.  For more detailed information, we refer the
reader to \cite{End} and \cite{Ri} as well as \cite{AEJ}, \cite{Wa} and
\cite{Br}.

Let $v : F \to \Gamma \cup \{\infty\}$ be a valuation on the field $F$, where
$\Gamma$ is some linearly ordered abelian group.  Then we set $A_v$ to be the
valuation subring of $F$, $M_v$ to be the unique maximal ideal of $A_v$
(consisting of those elements $f \in F$ such that $v(f)>0$), and $U_v$ to
be the group of invertible elements of $A_v$.  We say $T$ is {\it compatible}
with $v$ (or $A_v$) if $1 + M_v \subseteq T$.  We denote the residue field
$A_v/M_v$ by $F_v$, and we set $\pi_v : A_v \to F_v$ to denote the
canonical epimorphism from $A_v$ onto $F_v$.

The strategy of the proof is as follows:  It is easy to reduce
the problem of constructing $H$-closures to the problem of extending a given
$H$-ordering $T$ of a field $F$ to an $H$-ordering $T'$ of any quadratic
extension $L = F(\root\of{t}), t \in T$, such that $T' \cap F = T$.  (Here $H
\cong C(I), S(I)$, or $D(I)$.)  In order to extend $T$ in this manner, we
first find a suitable $T$-compatible valuation $v$ on $F$ and then extend $v$
to a valuation $w$ on $L$.  We then extend the induced ordering $\bar T$ of the
residue field $F_v$ to $\hat T$ on the residue field $L_w$ of $L$ with
respect to the valuation $w$.  Finally we lift the ordering $\hat T$ from
the residue field $L_w$ to an ordering $\tilde T$ on $L$, and then show that
$\tilde T$ is the desired extension of $T$ from $F$ to $L$.

Suppose first that we are given some $S(I)$-ordering $T$ of $F$.  In this case,
$T$ is \lq\lq not exceptional\rq\rq\ in the sense of
\cite{AEJ, Definition~2.15}. Thus we can apply \cite{AEJ, Theorem~2.16} to
obtain the following.

\proclaim {Proposition 8.1} Let $T$ be any $S(I)$-ordering of $F$.
Then there exists a $T$-compatible nondyadic valuation $v$ of $F$ such
that $U_vT = T \cup -T$.  The set
$\bar{T}:=\pi_{v}(T\cap U_{v})$ is an $S(\emptyset)$-ordering of $F_v$.
\endproclaim

\demo{Proof}  By \cite{AEJ, Theorem~2.16}, we have a $T$-compatible
valuation
$v$ such that $U_vT = T \cup -T$. The last statement of the proposition
follows from this.  Indeed we have
$$ \frac{U_v}{U_v\cap T} \cong \frac{U_vT}{T} \cong \frac{T \cup -T}{T}.$$
Since $-1 \notin T$ we see that $F_v = \bar T \cup -\bar T$ and $-1 \notin
\bar T$.  Therefore $\bar T$ has index $2$ in $\dot F_v$.

Since $T$ is an $S(I)$-ordering on $F$, we see that there exist elements $t_1,
t_2, t_3 \in T$ such that $t_1 + t_2 +t_3=0$.  Dividing through by that
element $t_i$ whose value $v(t_i)$ is  minimal among the three elements
considered (say $t_{1}$), we may assume we have
$$-1 = t_2 +t_3, v(t_2), v(t_3) \geq 0$$
Passing to the residue field we obtain
   $\bar{t_1} + \bar{t_2} = -\bar 1$ in
$F_v$.  Since $-1 \notin \bar T$ we see that $\bar{t_i} \neq 0, i =  2,3$.
Thus $-1 \in \bar T + \bar T\setminus\bar T$, and $\bar T$ is  a
$S(\emptyset)$-ordering of
$F_v$, as claimed.

Observe also that $-1\notin\bar T$ implies $-1\neq 1$ and  $\ch F_v \neq 2$.
Thus $v$ is nondyadic.
\qed \enddemo

Next suppose we have a $C(I)$-ordering $T$ of $F$.  Then we may apply
\cite{Ef, Propositions 2.1 and 2.3 and Theorem 4.1}, to yield the following
result.

\proclaim {Proposition 8.2} Let $T$ be any $C(I)$-ordering of $F$.
Then there exists a $T$-compatible valuation ring $A_v$ of $F$ such
that $[U_vT:T] \le 2$ and $\dim _{\F_2} \Gamma/2\Gamma \ge |I|$,
where $\Gamma$ is the associated value group.  The set $\bar{T}:=\pi_{v}(T\cap
U_{v})$ is either
$\dot{F_{v}}$ itself or a $C(\emptyset)$-ordering of $F_v$.
\endproclaim
\demo{Proof}  Observe again that the last statement claiming that
$\bar{T}:=\pi_{v}(T\cap
U_{v})$ is either
$\dot{F_{v}}$ itself or a $C(\emptyset)$-ordering, and also the statement
$\dim _{\F_2} \Gamma/2\Gamma \ge |I|$, are consequences of the first part
of the proposition.  We have
$\frac{U_vT}{T} \cong \frac{U_v}{U_v\cap T}$, so $[U_v : U_v\cap T] \le 2$;
hence $U_v = U_vT$ or $[U_v:U_v \cap T] = 2$.  In the latter case, we see that
$\bar T = \pi_v(T\cap U_v)$ is a
$C(\emptyset)$-ordering as $-\bar 1 \in \bar T$.  Also observe that we have
$|I| + 1 = \dim_{\F_2} \frac{\dot F}{T} = \dim_{\F_2} \frac{\dot
F}{U_vT} + \dim_{\F_2}\frac{U_vT}{T}$.  From the hypothesis $[U_vT:T] \le
2$ we see that $\dim_{\F_2}\frac{U_vT}{T} \le 1$.  Hence
$\dim_{\F_2}\frac{\dot F}{U_vT} \ge |I|$.  Therefore
$\dim_{\F_2}\Gamma_v \ge \dim_{\F_2}\frac{\dot F}{U_vT} \ge |I|$ as
claimed.
\qed
\enddemo

\proclaim {Proposition 8.3}{\rm (Fan Trivialization
Theorem~\cite{Br, Theorem 2.7})} Let $T$ be any $D(I)$-ordering of $F$.
Then there exists a $T$-compatible valuation ring $A_v$ of $F$ such
that the set $\bar{T}:=\pi_{v}(T\cap U_{v})$ is either
an ordering of $F_v$ or a $D$-ordering of $F_v$. (When $\bar{T}$ is an
ordering, $T$ is called a {\it valuation fan}.)  Moreover, the valuation
$v$ may be chosen such that $v(T)$ contains no convex subgroups of $v(F)$.
\endproclaim

Now suppose that we have an $S(I)$-ordering (respectively $C(I)$-,
$D(I)$-ordering) $T$ together with a $T$-compatible valuation $v$ on
$F$.  Assume $t \in T$, and let $K = F(\sqrt{t})$.  Our goal is to
find an $S(I)$-ordering (respectively $C(I)$-, $D(I)$-ordering) $T'$
of $K$ such that $T' \cap F = T$ and $\dot F/T \cong \dot K / T'$ is
the isomorphism of multiplicative groups induced by the inclusion $F
\hookrightarrow K$.  Note that if $T' \cap F = T$, then the map $\dot
F/T \to \dot K/T'$ is injective, so we need only worry about
surjectivity.  Then recall the well-known Krull's Theorem (\cite {Ri,
Theorem 5}):

\proclaim {Theorem 8.4}{\rm (Krull)}  Let $F$ be a field and $\tilde F$ any
overfield of $F$.  Any valuation $v$ in $F$ can be extended to a valuation
$\tilde v$ in $\tilde F$.
\endproclaim

Thus we see that there exists a valuation $w$ on $K$ which extends $v$.  In
order to proceed we make the following convenient reduction.

\proclaim {Lemma 8.5}  Assume that $T_1 \subseteq T_2$ are respectively
$S(I_1)$- and $S(I_2)$-orderings of $F$, and let $t \in
T_1\setminus\dot{F}^2$.  Let $K =
F(\sqrt{t})$.  Suppose $T_1'$ is an extension of $T_1$ to an $S(I_1)$-ordering
of $K$.  Then $T_2' := T_1'T_2$ is an $S(I_2)$-ordering of $K$ extending $T_2$.
\endproclaim
\demo{Proof} We first show that $T_2' \cap F = T_2$.  By definition,
$T_2 \subseteq T_2' \cap F$ , and if $f \in T'_{2}\cap F$ then there
exists $t'_1 \in T_1', t_{2}\in T_{2}$ such that $f = t'_1t_2$.  This
implies $t'_1 \in F \cap T_1' = T_1 \subseteq T_2$, and $f\in T_2$.
Thus $T_2' \cap F = T_2$.

Consider  the natural homomorphism $\varphi_2 : \dot F/T_2 \to \dot K /
T_2'$ induced by the inclusion map $F \hookrightarrow K$.  Because $T_2' \cap F
= T_2$ we see that $\varphi_2$ is injective. Consider the following
diagram:
$$
\CD
\dot F/T_1   @>\varphi_1>>  \dot K/T_1'\\
@VVV       @VVV\\
\dot F/T_2   @>\varphi_2>>  \dot K/T_2'
\endCD
$$
Since we know that
$\varphi _1 : \dot F/T_1 \to \dot K/T_1'$ is bijective, and since $T_1'
\subseteq T_2'$, we see that $\varphi _2$ is also surjective.

Finally we shall show that $T_2'$ is an $S(I_2)$-ordering by checking that
conditions (a),(b),(c)  of Corollary~7.7 hold.  Since $T_2' \cap F = T_2$, we
see that $-1
\notin T_2'$. As $-1\in T'_{1}+T'_{1}\subseteq T'_{2}+T'_{2}$, we see
that $T'_{2}$ satisfies condition (a).

Suppose $s=u+av\in K$ with $u,v\in T'_{2}$ and $a\not\in
   (T'_{2}\cup -T'_{2})$.  By definition of $T'_{2}$, $u,v$ can be
written $u=u'_{1}u_{2}, v=v'_{1}v_{2}$ with
$u'_{1},v'_{1}\in(T'_{1}\cup -T'_{1})$, $u_{2},v_{2}\in T_{2}$.
Then $su_{2}^{-1}=u'_{1}+(av_{2}u_{2}^{-1})v'_{1}$.
Because $av_{2}u_{2}^{-1}\not\in (T'_{1}\cup -T'_{1})$, the
$T'_{1}$-rigidity of $K$ implies
$su_{2}^{-1}\in T'_{1}\cup (av_{2}u_{2}^{-1})T'_{1}
$, and thus $s\in T'_{2}\cup aT'_{2}$, giving condition (b).

Finally, to check condition (c), observe that $\dot K/T_2' \cong \dot
F/T_2 \cong \coprod_{i\in I_2\cup\{x\}}(C_2)_i$.  Thus $T_2'$ is an
$S(I_2)$-ordering which extends $T_2$.
\qed
\enddemo

\proclaim{Lemma 8.6} Assume that $T_1 \subseteq T_2$ are respectively
$C(I_1)$- and $C(I_2)$-orderings of $F$, and let $t \in
T_1\setminus\dot{F}^2$.  Let $K =
F(\sqrt{t})$.  Suppose $T_1'$ is an extension of $T_1$ to a $C(I_1)$-ordering
of $K$.  Then $T_2' := T_1'T_2$ is a  $C(I_2)$-ordering of $K$ extending $T_2$.
\endproclaim
\demo{Proof}  The proof is identical to that of Lemma 8.5, except that one
must now
check that $-1 \in T_2'$.  Since $T_2' \cap F =
T_2$, we see
$-1 \in T_2'$.
\qed
\enddemo

\proclaim{Lemma 8.7} Assume that $T_1 \subseteq T_2$ are respectively
$D(I_1)$- and $D(I_2)$-orderings of $F$, and let $t \in
T_1\setminus\dot{F}^2$.  Let $K =
F(\sqrt{t})$.  Suppose $T_1'$ is an extension of $T_1$ to a $D(I_1)$-ordering
of $K$.  Then $T_2' := T_1'T_2$ is a  $D(I_2)$-ordering of $K$ extending
$T_2$.
\endproclaim
\demo{Proof}  Again the proof takes the same arguments as in the proof
of Lemma~8.5 to show that $T'_{2}$ extends $T_{2}$, that $-1\not\in
T'_{2}$ and  that $K$ is $T'_{2}$-rigid. Let us prove  $T'_{2}+T'_{2}=T'_{2}$.
Consider $u,v\in T'_{2}$ and write them as above,
   $u=u'_{1}u_{2},v=v'_{1}v_{2}$, with $u'_{1},v'_{1}\in T'_{1}$ and
   $u_{2},v_{2}\in T_{2}$.  Then
   $u+v=u_{2}(u'_{1}+(v_{2}u_{2}^{-1})v'_{1})$.  We know that $-1\not\in
   T'_{2}$, and this implies that $v_{2}u_{2}^{-1}\not\in -T'_{1}$.  If
   $v_{2}u_{2}^{-1}\in T'_{1}$, then $(u+v)u_{2}^{-1}\in
   T'_{1}+T'_{1}=T'_{1}$ and $u+v\in T'_{2}$.
   The remaining possibility is $v_{2}u_{2}^{-1}\not\in T'_{1}\cup
   -T'_{1}$, and by $T'_{1}$-rigidity of $K$, we have
   $(u+v)u_{2}^{-1}\in T'_{1}\cup
   (v_{2}u_{2}^{-1})T'_{1}$ and $u+v\in T'_{2}$. Hence condition (2) holds.
\qed
\enddemo

   We consider the following situation.  Assume that $v : F \to \Gamma_v
   \cup \{\infty\}$ is a valuation on the field $F$, with valuation ring $A_v$
   and maximal ideal $M_v$.  Let
   $F_v = A_v/M_v$ be the residue field, and denote by $\pi_v$ the canonical
   homomorphism of
   $A_v$ onto its quotient ring
   $F_v$.  

\proclaim{Lemma 8.8} Assume that $v$ is a valuation on the field $F$ and
that
$ T_0$ is an $S(I_{0})$-ordering of $\dot F_v$  for some (possibly
empty) set $I_{0}$.
Set $T_1 = \pi_v^{-1}(T_0)$.  Then the group $T = T_1\dot F^2$ is an
$S(I)$-ordering of $F$ with $|I| = dim_{\F_2}(\frac{\dot F}{T\cup -T})$.
\endproclaim
\demo{Proof} What is needed is to check that the conditions in
Corollary~7.7  hold for $T$.  First, suppose that $-1 \in T$.
Then $-1 = t_0f^2$ for some $t_0 \in T_1, f \in \dot F$.  Hence $f^2 =
(-t_0)^{-1} \in -T_1 \subseteq U_v$, and so $f \in U_v$ as well.
Passing to the residue field $F_v$ and knowing $\dot{F_{v}}^2\subseteq
T_{0}$
we see $-1 = \bar t_0 \bar f^2 \in T_0$, which is a contradiction.
    Thus we must have $-1 \notin T$.
Since $-1\in T_{0}+T_{0}$, we have $-1+m\in T_{1}+T_{1}$ for some $m$
in the maximal ideal of the valuation, and $-1+m\in -T_{1}\subset T$.
This shows that the level of $T$ is $2$.

To see that $F$ is $T$-rigid, let $a \in \dot F \setminus(T \cup -T),
t_1, t_2 \in T$, and consider $b := t_1 + t_2a$.  We consider various
possibilities for $v(t_1)$ relative to $v(t_2a)$.
First suppose that $v(t_1) = v(t_2a)$.  Then $b = t_1(1 +
t_1^{-1}t_2a)$, with $ u := t_1^{-1}t_2a \in U_v$.  Since $a \notin T
\cup -T$, we see that $\pi_v(u) = \bar u \notin T_0 \cup - T_0$.
(Otherwise $u \in \pi_v^{-1}( T_0) = T_1 \subseteq T$ or $u \in
-\pi_v^{-1}( T_0) = -T_1 \subseteq -T$ and hence $a \in T \cup -T$, a
contradiction.) Since we are assuming $F_v$ is $ T_0$-rigid, we see
that $1 + \bar u \in T_0 \cup \bar u T_0 $.  Hence $1 + u \in
\pi_v^{-1}( T_0 \cup \bar u T_0) = T_1 \cup uT_1$.  Thus, rewriting $u
= t_1^{-1}t_2a$ and multiplying through by $t_1$, we see $$ b = t_1 +
t_2a \in T_1 \cup aT_1 \subseteq T \cup aT $$ as required.
Now assume that $v(t_1) \neq v(t_2a)$.  If $v(t_1) < v(t_2a)$, then again
let $b
= t_1(1 + u)$, where $u = t_1^{-1}t_2a$.  Now, however, $v(u) > 0$, so $1 + u
\in 1 + M_v \subseteq T_1 = \pi_v^{-1}( T_0)$, and thus $b \in T $.
If $v(t_1) > v(t_2a)$, set $b = at_2(1 + t_1t_2^{-1}a^{-1})$.  We see
$v(t_1t_2^{-1}a^{-1}) > 0$, and therefore $b \in aT$.
In each case $b = t_1 + at_2 \in T  \cup aT $ as
desired.

It remains to see that $\dot F /T \cong \coprod_{i\in
I\cup\{x\}}(C_2)_i$.  This condition follows from the fact that $\dot F
/T$ is an $\F_2$-vector space and that $dim_{\F_2} \dot F/T$ is $1 +
|I|$.
\qed
\enddemo

We have the analogue to Lemma 8.8 for the case of $C(I)$-orderings.

\proclaim{Lemma 8.9}  Assume that $v$ is a valuation on the field $F$ such that
$[\Gamma_v : 2\Gamma_v] \geq 2$.  Let $ T_0$ be   $\dot
F_v$  or a $C(I_0)$-ordering of $F_{v}$ for some (possibly empty) set $I_0$.
Set $T_1 = \pi_v^{-1}( T_0)$.  Then the group $T = T_1\dot F^2$ is a
$C(I)$-ordering of $F$ with $|I| = dim_{\F_2}(\frac{\dot F}{T}) - 1$.
\endproclaim
\demo{Proof} We must check that the conditions of Proposition 7.2 hold
for $T$.  Clearly if $-1 \in T_0$, then $-1 \in T_1 \subseteq T$.  To
see that $F$ is $T$-rigid, one applies the same argument as in Lemma
8.8.  As in the case for $S(I)$-orderings, $\dot F /T$ is clearly an
$\F_2$-vector space.  Since
$[\Gamma_{v}\colon 2\Gamma_{v}]\ge 2$, its dimension is strictly
positive and thus may be written
$dim_{\F_2} (\dot F/T) = 1 + |I|$.
\qed
\enddemo

Again, we also have the analogue to Lemma 8.8 for the case of $D(I)$-orderings.
\proclaim{Lemma 8.10}{\rm (\cite{Br})} Assume that $v$ is a
valuation on the field $F$.  Let $ T_0$ be a fan of $\dot F_v$.  Set
$T_1 = \pi_v^{-1}( T_0)$.  Then the group $T = T_1\dot F^2$ is a fan
(i.e.  a $D(I)$-ordering) of $F$.
\endproclaim

We now formulate the key results in this section.

\proclaim{Theorem 8.11} Let $T$ be any $S(I)$-ordering of $F$ and let $L =
F(\sqrt{t}), t \in T$.  Then there exists an $S(I)$-ordering $T'$ on $L$ such
that $(L,T')$ is an $S(I)$-extension of $(F,T)$.

\endproclaim

\demo{Proof} From Proposition 8.1, we see that there exists a
nondyadic $T$-compatible valuation ring $A_v$ in $F$ such that $U_vT
= T \cup -T$ and that $\bar T := \pi_v(U_v \cap T)$ is an
$S(\emptyset)$-ordering of $F_v$.  As $\pi_{v}^{-1}(\bar{T})=(U_{v}\cap
T)(1+M_{v})$ and because $(1+M_{v})\subseteq T$, one has
$T_{1}:=\pi_{v}^{-1}(\bar{T})\dot{F}^2\subseteq T$.  By Lemma 8.8, we
see that $T_{1}$ is an $S(J)$-ordering in $F$ for a suitable set $J$.

Let $w$ be any valuation of $L$ which extends $v$.  Let $L_{w}$ denote
its residue field, and $\Gamma_v, \Gamma_w$ denote the valuation
groups of $v$ and $w$.  We may assume $\Gamma_v \subseteq \Gamma_w$,
and we set $e = [\Gamma_w: \Gamma_v]$, the ramification degree of $w$
with respect to $v$, and $f = [L_w:F_v]$, the residue class degree of $w$
with respect to $v$.  It is well known that $ef \le [L:F] = 2$ and in
particular we have $f=[L_{w}:F_{v}]\le 2$.  More precisely, one has
$L_{w}=F_{v}(\sqrt{\pi_{v}(u_{0})})$ with $u_{0}=1$ if $f=1$, and
$u_{0}/t\in\dot F^2$ if $f=2$.  By Proposition~5.6 and Remark~5.5,
$C_{4}$-orderings are known to admit $C_{4}$-closures of the same level,
and as $\pi_{v}(u_{0})\in\bar{T}$, the $S(\emptyset)$-ordering
$\bar{T}$ admits an $S(\emptyset)$-extension $\tilde{T}$ to
$F_{v}(\sqrt{\pi_{v}(u_{0})})=L_{w}$.  Calling
$T_{2}=\pi_{w}^{-1}(\tilde{T})L^2$, Lemma~8.8 implies that $T_{2}$ is
an $S(K)$-ordering of $L$ for a suitable set $K$.

Let us first show that $T_{1}=T_{2}\cap F$.  By definition of $T_{1}$,
an element $s\in T_{1}$ has the same square class as an element $u\in
U_{v}$ such that $\pi_{v}(u)\in\bar{T}\subseteq\tilde{T}$.  This
implies that $\pi_{w}(u)\in\tilde{T}$, and thus $u$ and $s$ are in
$T_{2}$. This shows $T_{1}\subseteq T_{2}\cap F$.

For the reverse inclusion, we state the following claim:

\proclaim{Claim} With notation as above, one has
$\dot{L}=U_{w}\dot{F}\cap\,\root\of{t}U_{w}\dot{F}$.
\endproclaim

\demo{Proof} We know that $e\le 2$.  If $e = 1$, then $\dot L = \dot
FU_w$ and we are done.  If $e = 2$, then $f=1$ and we may show that
$w(\root\of{t})\not\in \Gamma_v$.  Otherwise $\root\of{t}=xu$ with
$x\in F$ and $u\in U_{w}$, and denoting by $\sigma$ the nontrivial
element of the Galois group $\Gal(L/F)$, we know that
$\frac{\sigma(\root\of{t})}{\root\of{t}}=-1$ and thus
$\pi_{w}(\frac{\sigma(\root\of{t})}{\root\of{t}}) =
\pi_{w}(\frac{\sigma(u)}{u})=-1$.  Since $f=1$, $L_{w}=F_{v}$, and so
$\pi_{w}(\frac{\sigma(u)}{u})$ must also be $1$.  Since the valuation $v$ is
not dyadic, this would be a contradiction.
      Thus we see that since
$\Gamma_w \cong \dot L/U_w, \Gamma_v \cong \dot F/U_v$, and
$[\Gamma_w : \Gamma_v] = 2$, the factor group $\dot{L}/U_{w}\dot{F}$ is
$\{1,\root\of{t}\}$, and we can write
   $\dot L = U_w\dot F \cup \root\of{t}U_w\dot F$. \qed
   \enddemo

We now finish the proof of the theorem.  If $\alpha\in T_{2}\cap
F$, we may write $\alpha= u\lambda^2$ with
$u\in\pi_{w}^{-1}(\tilde{T}), \lambda\in\dot L$, and writing
$\lambda=\sqrt{t}^\eta u_{1}g$ with $u_{1}\in U_{w}, g\in\dot F,
\eta=0$ or $1$, this yields $\alpha =uu_{1}^2t^\eta g^2$.  Since
$t^\eta g^2\in T_{1}$, we may assume $\alpha =uu_{1}^2$.  Then
$\pi_{v}(\alpha)=\pi_{w}(\alpha)\in\tilde{T}\cap F_{v}=\bar{T}$ and
$\alpha\in T_{1}$.  This proves $T_{1}=T_{2}\cap F$.

We define a new subgroup $T'_{2}$ of $\dot L$ as follows.
\roster
\item If $\sqrt{t}\in(T_{2}\cup -T_{2})$, set $T'_{2}=T_{2}$.
\item If $\sqrt{t}\not\in(T_{2}\cup -T_{2})$ and
$[\Gamma_{w}:\Gamma_{v}]=1$, again set $T'_{2}=T_{2}$.
\item If $\sqrt{t}\not\in(T_{2}\cup -T_{2})$ and
$[\Gamma_{w}:\Gamma_{v}]=2$, set  $T'_{2}=T_{2}\cup\sqrt{t}T_{2}$.
\endroster
Then again $T_{1}=T'_{2}\cap F$, the only thing to prove being that in
the third case, $\sqrt{t}T_{2}\cap F\subseteq T_{1}$.  But if
$\alpha\in\sqrt{t}T_{2}\cap F$ we have $\alpha=\sqrt{t}ug^2$ with
$u\in U_{w}, g\in\dot F$ and this implies $w(\sqrt{t})\in\Gamma_{v}$,
contradicting $[\Gamma_{w}:\Gamma_{v}]=2$.  This shows that
$\sqrt{t}T_{2}\cap F=\emptyset$ in the third case.

Since $T_{2}$ is an $S(K)$-ordering, it is easy to check that
conditions (1)-(3) of Proposition~7.6 hold for $T'_{2}$ and to see that
$T'_{2}$ is also an $S(K')$-ordering for a suitable set $K'$.

We want to show that the injection
$\dot{F}/T_{1}\longrightarrow\dot{L}/T'_{2}$ is also surjective, which
reduces to showing that
   $\dot L=T'_{2}\dot F$. We already know $\dot L=U_{w}\dot
F\cup\sqrt{t}U_{w}\dot F$, and by Lemma 8.1, $U_{w}\subseteq
T_{2}\cup -T_{2}$. This gives us $U_{w}\dot F\subseteq T_{2}\dot
F\subseteq T'_{2}\dot F$.
In cases (1) and (3), one has $\sqrt{t}\in T'_{2}\cup -T'_{2}$, and so
$\dot L\subseteq T'_{2}\dot F$.  In case (2), there exists
$x_{0}\in\dot F$ such that $\sqrt{t}x_{0}\in U_{w}\subseteq T_{2}\dot
F$.  So $\sqrt{t}\in T_{2}\dot F$, finishing the proof that
$\dot{F}/T_{1}\longrightarrow\dot{L}/T'_{2}$ is an isomorphism.

We have proved so far that $(L,T'_{2})$ is an $S(J)$-extension of
$(F,T_{1})$, and that $T_{1}$ is contained in the $S(I)$-ordering $T$.
We may then apply Lemma~8.5 to show that $(L,T_{1}T'_{2})$ is an
$S(I)$-extension of $(F,T)$, and the theorem is proved.
\qed
\enddemo

\proclaim{Corollary 8.12} An $S(I)$-ordered field $(F,T)$ admits an
$S(I)$-closure.
\endproclaim
\demo{Proof}  Let $\Cal S$ be the set of extensions $(L,S)$ of $(F,T)$ inside
$F(2)$ such that $S$ is an $S(I)$-ordering on $L$.  Then by a Zorn's Lemma
argument $\Cal S$ has a maximal element $(K,T_0)$ with $\dot K/T_0 \cong \dot
F/T, T = T_0 \cap F$, and $T_0$ is an $S(I)$-ordering on $K$.  We are done
by Corollary~4.3 if we can show  $T_0 = \dot K^2$.  If not, choose $t \in
T_0 \backslash \dot K^2$.  Then by Theorem 8.11 we can extend $T_0$ to an
$S(I)$-ordering on $K(\sqrt{t})$, contradicting the maximality of $(K,T_0)$.
\qed
\enddemo

Corollary 8.12 can be reformulated in the
language of Galois theory as in the following corollary, which tells us that
a certain family of subgroups of $G_F := \Gal(F(2)/F)$ occurs whenever $G_F$
contains certain subquotients of $G_F$. Observe that in Corollary~8.13
we do not specify the action of the outer factor $\zz_{2}$ on the normal
subgroup $(\zz_{2})^{I}$ as this action depends upon a subtler analysis
of the roots of unity belonging to the fields under consideration.

\proclaim{Corollary 8.13} Let $F$ be a field of characteristic $\neq 2$.
Suppose that we have a tower of
field extensions $F \subset N_{1} \subset N_{2} \subset N_{1}^{(3)} \subset
F(2)$, where $N_1^{(3)}/N_2$ is a Galois extension and $\Gal(N_1^{(3)}/N_2)
\cong
(C_4)^I \rtimes C_4$ for $I$ some nonempty set.
Then $G_F$ = $\Gal(F(2)/F)$ contains the closed subgroup $(\zz_{2})^I
\rtimes \zz_{2}$.
\endproclaim
\demo{Proof}Let $F \subset N_{1} \subset N_{2}
\subset N_{1}^{(3)} \subset F(2)$ be a tower of field
   extensions, where $N_{1}^{(3)}/N_{2}$ is a Galois extension and
$\Gal(\tfrac{N_{1}^{(3)}}{N_{2}})
\cong (C_{4})^{I} \rtimes C_{4}$ for $I$ some nonempty set.
Set $T = \{t \in \dot{N}_{1} \mid (\sqrt{t})^{\sigma} = \sqrt{t}$ for each
$\sigma \in \Gal(N_{1}^{(3)}/N_{2})\}$. From
Definition 7.1 we see that $T$ is an $S(I)$-ordering of $N_{1}$. From Corollary
8.12 it follows that there exists a field extension $N$ of $N_{1}$ such that
$\dot{N}^{2}$ is an $S(I)$-ordering of $N$ and $\dot{N}^{2} \cap N_{1}
= T$. Then Proposition 8.1 implies the existence of an $\dot{N}^{2}$-compatible
valuation ring $A_v$ of $N$ such that $U_{v}\dot{N}^{2} = \dot{N}^{2} \cup
-\dot{N}^{2}$.

It is well known that an $\dot{N}^{2}$-compatible valuation $v$ on $N$ is
2-henselian. Moreover $N$ is a rigid field (and is $S(I)$-closed).  In
Proposition 8.1 we observed that
$v$ is a nondyadic valuation (i.e., char $F_{v} \neq 2$) and in this case it
follows from basic valuation theory (see e.g. [End, \S 20]) that we have a
split short exact sequence
$$1 \longrightarrow I_{v} \longrightarrow G_{N}(2) \longrightarrow
G_{{N}_{v}} (2) \longrightarrow 1,$$
where $I_{v}$ is the inertia subgroup of $G_{N}(2) :=
\Gal(N(2)/N) = \Gal(F(2)/N)$
   and $N_{v}$ is the residue field of $v$.
Moreover it is well known that $I_{v}$ is an abelian group. (See e.g.
\cite{EnKo}.)

Because $\dot{N}^{2}$ is an $S(I)$-ordering of $N$ we see that $s(N)=2$. In
particular $N$ is not a formally real field, and so $G_{N}(2)$
is a torsion-free group. (See [Be].) Therefore using Pontrjagin's duality and
the well-known structure of abelian divisible groups, we see that
$I_{v}
\cong (\zz_{2})^{J}$ for some set
$J$. (See e.g. \cite{RZ,  \S 4.3,
Theorem 4.3.3}.)

Because $\dot{N}^{2}$ is compatible with $v$ and
$$\frac{U_{v}}{U_{v} \cap \dot{N}^{2}} \quad \cong \quad \frac{\dot{N}^{2}
\cup
- \dot{N}^{2}}{\dot{N}^{2}},$$
we see that $\mid \dot{N}_{v}/\dot{N}_{v}^{2} \mid = 2$. Hence $G_{{N}_{v}} (2)
\cong \zz_{2}$. Since
$\dot{N}^{2}$ is an $S(I)$-ordering of $N$, it follows that the cardinality
of $I$
is the same as the cardinality of $J$. Hence
$I_{v} \cong (\zz_{2})^{I}$. Since the Galois group $G_{N}(2) = I_{v} \rtimes
\zz_{2}$ is a closed subgroup of $G_{F}$, the proof is completed.
\qed
\enddemo

In the case of $C(I)$-orderings, we cannot always find a closure.  The
problem arises from the fact that the valuation whose existence is guaranteed
by Proposition 8.2 may be dyadic, and thus the appropriate modification of
Theorem 8.11 will not go through.  For $S(I)$- and $D(I)$-orderings we do not
have this problem, as the valuation in question will be nondyadic. 
Example~8.14 below constructs a $C(1)$-ordered field which we show in
Proposition~8.15 does not admit a
$C(1)$-closure.

\example{Example 8.14} Recall that a field $K$ of characteristic $2$ is called
{\it perfect} if $K^2 = K$.  S. MacLane has shown that for any field $K$ of
characteristic $2$, there exists a field $F$ of characteristic $0$ with a
valuation $v : F \to \zz \cup\{\infty\}$ such that $F_v \cong K$ (\cite{Mac,
Theorem~2}. For some more general theorems on valued fields with prescribed
residue fields, see \cite{Ri, Chapter I}).
Then let $F$ be such a field where $F_v = K$ is a field of
characteristic $2$ which is not perfect.  Let $T_0$ be a multiplicative
subgroup of $\dot K$ of index $2$ in $\dot K$ such that $\dot K^2 \subsetneq
T_0 \subsetneq \dot K$.  Let $T = \dot F^2\pi_v^{-1}(T_0)$, a subgroup of
$\dot F$. Here $\pi_v$ is the residue map $U_v \lra \dot{K}$.
Then $\mid\dot{F}/T\mid = 4$, and one can choose as representatives
of the factor group $\dot{F}/T$ the elements $1, u, \rho, \rho u$ where
$v(\rho) = 1, u \in U_v$, and $\pi_v(u) \notin T_0$.

We claim that $F$ is $T$-rigid. Since any element in $\rho T$ or in $\rho uT$
lies outside of $U_vT$, we see that all elements of $\rho T \cup \rho uT$ are
$T$-rigid.  (See \cite{AEJ, Proposition~1.5.})  Consider an element
$\alpha =
t_1 + t_2u \in T + uT$, with $t_1, t_2 \in \dot{F}$.  Then $\alpha =
t_2(t_1t_2^{-1} +
u)$, so it is enough to show $t_1t_2^{-1} + u \in T \cup uT$.  Thus we may
restrict our attention to elements which can be written as $tf^2 + u$, where
$t \in \pi_v^{-1}(T_0), f \in \dot F$.  If $v(f) = 0$, then $tf^2 + u \in U_v
\subseteq T \cup uT$.  If $v(f) > 0$, then $tf^2 + u = u(1 + tf^2u^{-1}) \in
uT$.  Finally, if $v(f) < 0$, then $tf^2 + u = tf^2(1 + uf^{-2}t^{-1}) \in
T$.  Thus $F$ is $T$-rigid.

Since $-1 \in T_0$, we have $-1 \in T$, and
$T$ is a $C(1)$-ordering of $F$.
Observe that $T \neq \dot F^2$ and $(F,T)$ is not $C(1)$-closed.
\comment
Indeed set $\tilde{u}=\pi_v(u)\in \dot{K}$. Then 
$\tilde{u},1+\tilde{u}$ are linearly independent elements
in $\dot{K}/\dot{K}^2$. Therefore either $1+\tilde{u}$ or
$(1+\tilde{u})\tilde{u}\in T_0 \setminus \dot{K}^2$.
\endcomment

\proclaim{Proposition 8.15} The $C(1)$-ordered field $(F,T)$ does not admit a
$C(1)$-closure.
\endproclaim
\demo{Proof}  Recall that a valuation $\nu$ on a field $L$ is said to be
$T${\it -coarse} if $\nu(T)$ contains no nontrivial convex subgroups
of the valuation group $\Gamma_\nu$ of $\nu$.  Suppose that $F \subsetneq N
\subsetneq F(2), \dot N^2
\cap F = T$, and $\dot N^2$ is a $C(1)$-ordering of $N$.  Then applying
\cite{AEJ, Corollary~2.1.7} or \cite{Wa, Theorem~2.16}, we see that there
exists a
$\dot N^2$-compatible valuation $w$ on $N$ such that $[U_w\dot N^2 : \dot N^2]
\leq 2$.  This means that $\mid U_w/U_w\cap\dot{N}^2 \mid \leq 2$. We may
further choose $w$ to be the unique finest $N^2$-coarse $N^2$-compatible
valuation on
$N$ (see \cite{AEJ, Theorem~3.8}). Consider
$z := $ the restriction of the valuation $w$ to $F$.  First observe that $z$
is a $T$-compatible valuation on $F$.  Indeed, from $M_w \cap F = M_z$ we get
$(1 + M_w) \cap F = 1 + M_z$.  Thus we have
$$1 + M_z = (1 + M_w) \cap F \subseteq \dot N^2 \cap F = T.$$

Let $\Delta$ be the maximal convex subgroup of $\Gamma_z$ contained in
$z(T)$.  Then set $y$ to be the composite valuation
$$y : \dot F \overset z\to\longrightarrow \Gamma_z
\overset\rho\to\longrightarrow \Gamma_z/\Delta,$$
where the last map $\rho : \Gamma_z \to \Gamma_z/\Delta$ is the natural
projection.  Then, following the notation of \cite{AEJ, Definition~2.2}, the
valuation ring $A_y = O_F(U_zT,T)$, and $y(T)$ contains no nontrivial
convex subgroups of the value group $\Gamma_y = \Gamma_z/\Delta$
(\cite{AEJ,
Lemma~3.1 and Proposition~3.2}), so $y$ is $T$-coarse. Observe that $y$ is
also $T$-compatible.
However, since $\Gamma_v=
\zz$ and $v(T) = 2\zz \ne \zz$, the valuation $v$ is also $T$-coarse.
Hence, by \cite {AEJ, Corollary~3.7}, we see that the valuations $v$
and $y$ are
comparable.  Since $A_v$ is a maximal proper subring of $F$ (because $ \Gamma_v
= \zz$), we see that $A_v \supseteq A_y\supseteq A_z$.  However, since $M_z
\supseteq M_y \supseteq M_v$ and $2 \in M_v$, we see that both valuations $y$
and $z$ are dyadic.  Since $F_z \subseteq F_w$, it follows that $w$ is also
dyadic.  But from  \cite{AEJ, Theorem~3.8 and Lemma~4.4}, it
follows that
$w$ cannot be a dyadic valuation.  Indeed, $[D_N\langle 1, -n^2\rangle \dot
N^2 :
\dot N^2] = 4 > 2$ for all $n \in \dot N$.  Thus we have a contradiction, and
there can be no $C(1)$-closure of $(F,T)$. 
\qed
\enddemo
\endexample

\definition{Remark 8.16} Example~8.14 is analogous to Proposition~4.10.  What
makes this example striking when compared to Proposition~4.10 is that here we
have
$\mid\dot{F}/T\mid=4<\infty$, but in Proposition~4.10
$\mid\dot{F}/T\mid=\infty$.  Although this example is a relatively simple
consequence of the work in \cite{AEJ}, it seems to be the first example
where the Witt ring of a field with finitely many square classes is
realizable as a ``Witt ring of
$T$-forms over some field $F$'', but it is not realizable as an actual
Witt ring of any field extension $K$ of $F$.  We make this
last comment more precise.

First observe that, analogous to the definition of reduced Witt rings of
fields, one may define $W_T(F)$ for any subgroup $T$ of $\dot{F}$ which
contains all nonzero squares in $F$.  One possible definition is as follows:
(See also \cite{La2, Corollary~1.27} and \cite{Sc, Chapter~2,
\S~9}.)

Let $\zz[\dot{F}/T]$ be the group ring of $\dot{F}/T$ with coefficients
in $\zz$.  Let  $J$ be the ideal of $\zz[\dot{F}/T]$ generated by
\roster
\item $[T]+[-T]$,
\item $[aT]+[bT]-[(a+b)T]-[ab(a+b)T], (a,b,a+b\in\dot{F})$,
\item $[aT][bT]-[abT],(a,b\in\dot{F})$.
\endroster
Then we set $W_T(F)=\zz[\dot{F}/T]\diagup J$.

A systematic study of $W_T(F)$ for $H$-orderings $T$ of $F$ is very
desirable, but it is not pursued in this particular paper.  Here we
just point out that if $T$ is any $C(1)$-ordering of $F$ then
$W_T(F)\cong W(\qq_p)$, where $p$ is any prime such that $p\equiv
1\pmod 4$, and $\qq_p$ is the field of $p$-adic numbers.

Since $T$ is a $C(1)$-ordering in $\dot F$
and $\dot\qq_p^2$ is a $C(1)$-ordering in $\qq_p$ (see Proposition~7.2
and \cite{L1, Chapter~6}), we see that there exists a group
homomorphism $\varphi\colon\dot F/T \lra \dot\qq_p/\dot\qq_p^2$ such
that $\varphi$ takes any relation in the form (1), (2) or (3) above
again to a relation of the same type.
Using the same argument for $\varphi^{-1}$ rather than $\varphi$, we
see that $\varphi$ indeed induces an isomorphism
$\tilde{\varphi}\colon W_T(F)\cong W(\qq_p)$.  
\enddefinition

\bigskip

Similar to  Proposition~4.12, we have the following proposition.

\proclaim{Proposition 8.17} Let $(F,T)$ be the field $F$ with 
$C(1)$-ordering
$T$ constructed in Example~8.14 above.  Then there is no field extension
$K/F$ with $C(1)$-ordering $\dot K^2$ which is a $T$-extension
of $(F,T)$.  (Equivalently, $W_T(F)$ cannot be realized as $W(K)$ for any
field extension $K$ of $F$.)
\endproclaim

\demo{Proof} Suppose to the contrary that there exists a field
extension $K/F$ such that $\dot K^2$ is a $C(1)$-ordering of $K$
and $(\dot K,\dot K^2)$ is a $T$-extension of $(F,T)$.  Assume that
both $K$ and a quadratic closure $F(2)$ of $F$ are contained in
some common overfield so that we can consider the field $L=K\cap
F(2)$.  The natural isomorphism $\psi\colon\dot F/T
\lra\dot K/\dot K^2$ factors through $\theta\colon\dot F/T \lra \dot
L/(\dot K^2 \cap L)$.  Because $\psi$ is injective, so is $\theta$.
Observe that $\theta$ is also surjective.  Indeed since $\psi$ is
surjective, we see that for each $l\in\dot L$ there exists an element
$f\in\dot F$ such that $lf^{-1}\in\dot K^2\cap L$.  Thus we see that
$(L,\dot K^2\cap L)$ is a $T$-extension of $(F,T)$.

We claim that $(\dot L,\dot K^2\cap L)$ is a $C(1)$-closure
of $(F,\dot T)$.
Observe that $\dot K^2\cap L=\dot L^2$.  Indeed if $k^2\in
L,k\in\dot K$ then $k\in\dot K\cap F(2)=\dot L$.  Since $\dot L^2
\subset \dot K^2 \cap L$ is obvious, we see that $\dot K^2\cap L=\dot
L^2$.  In order to conclude the proof, it is enough to show that
$\dot L^2$ is a $C(1)$-ordering in $\dot L$.  Because
$\sqrt{-1}\in\dot K$ we see that $\sqrt{-1}\in\dot L$ as well, and
$-1\in\dot L^2$.
>From the isomorphism $\theta\colon\dot F/T \lra \dot L/\dot L^2$ we
see that $\dot L/\dot L^2=C_2 \oplus C_2$.  By Proposition~7.2 it
remains only to show that $L$ is $\dot L^2$-rigid.  Consider an element
$a\in\dot L \diagdown \dot L^2$.  For any $l\in \dot L$ we have
$l^2+a\in\dot K^2\cup a \dot K^2$ because
$\dot K$ is $\dot K^2$-rigid and $\dot L^2=\dot K^2\cap L$.  Hence
$l^2 + a \in (\dot K^2\cap L)\cup(a\dot K^2\cap L)$.  Finally since
$\dot K^2\cap L=\dot L^2$ and $a\dot K^2\cap L=a\dot L^2$ we see that
$\dot L$ is $\dot L^2$-rigid.  \qed
\enddemo

\proclaim{Theorem 8.18}  A $C(I)$-ordered field $(F,T)$ possessing a
nondyadic $T$-compatible valuation ring $A_v$ as in Proposition~8.2 admits a
$C(I)$-closure.
\endproclaim
\demo{Proof}  The proof is essentially the same as the proof of
Theorem~8.11 and Corollary~8.12, and
we will follow the same plan and the same notation.
Applying Proposition~8.2, we find a valuation $v$ on $F$ such that
$\bar{T}:=\pi_{v}(U_{v}\cap T)$ is either $\dot{F_{v}}$ or a
$C_{\empty}$-ordering.  By assumption here this valuation is nondyadic. By
Lemma~8.9, $T_{1}$ is a $C(J)$-ordering contained in $T$. Taking any valuation
$w$ on $L=F(\sqrt{t})$ extending $v$, we extend $\bar{T}$ to $\tilde{T}$ in
$L_{w}$. We obtain, by Lemma~8.9, a $C(K)$-ordering $T_{2}$ in $L$. We enlarge
it to a $C(K')$-ordering $T'_{2}$, according to the three cases (1),
(2), (3), replacing $T_{2}\cup -T_{2}$ by $T_{2}$. The only serious
change is in proving that $\dot L=T'_{2}\dot F$. For this it is enough
to show that $U_{w}\subseteq T_{2}\dot F$, which can be done as follows.
If the index $[U_{v}T:T]=[\dot{F_{v}}:\bar{T}]$ is $1$, then
$[\dot{L_{w}}:\tilde{T}]=[U_{w}T_{2}:T_{2}]=1$ and $U_{w}\subseteq
T_{2}$. If this index is $2$, there exists $a\in U_{v}$ such that
$U_{w}\subseteq T_{2}\cup aT_{2}\subseteq T_{2}\dot F$.
This shows that $(L,T'_{2})$ is a $C(J)$-extension of $(F,T_{1})$, and
we apply Lemma~8.5 to show that $(L,T_{1}T'_{2})$ is a $C(I)$-extension
of $(F,T)$. We finish by applying the same argument as in
Corollary~8.12.
\qed
\enddemo

The following observation about valuations when $F$ contains a real-closed
field was pointed out to us by J.-L. Colliot-Th\'el\`ene.

\proclaim{Corollary 8.19} Let $v$ be a $T$-compatible valuation with value
group
$\Gamma$, and denote by $U_{v}$ the units of the valuation ring.  Suppose
there exists an integer $n>1$ such that any
$n$-divisible subgroup of $\Gamma$ is trivial.  Assume that $F$
contains a real-closed field $R$.  Then $R$ is contained in $U_v$, and
in particular the valuation is nondyadic.
\endproclaim

\demo{Proof}   Assume $F$ contains a real-closed field $R$.  If $a\in R$
is positive,
for the given integer $n$ there exists $b\in R$ such that $a=b^n$, and
thus $v(a)=nv(b)$.  Thus $v(a)$, being divisible by any power
of $n$, must be $0$, and the elements of $R$ must be units.  This
implies that the residue field $F_{v}$ contains $R$, and the valuation
$v$ cannot be dyadic.  In particular, the $T$-compatible valuation
which is known to exist, cannot be dyadic, and $(F,T)$ admits a
closure by Theorem~8.18.
\qed
\enddemo

\proclaim{Theorem 8.20}  A $D(I)$-ordered field $(F,T)$ admits a
$D(I)$-closure.
\endproclaim
\demo{Proof} We have already proved that $D$-orderings admit closures,
and thus we may assume that $|I|>1$.  It has also already been shown
in~\cite{Sch} that valuation fans admit closures.  Here is a more
general situation and a different proof, that consists again in
transpositions of the proofs of Theorem~8.11 and Corollary~8.12.  As
in Theorem~8.11, if $t\in T$ and $L=F(\sqrt{t})$, applying
Proposition~8.3, we find a valuation $v$ on $F$ such that
$\bar{T}:=\pi_{v}(U_{v}\cap T)$ is either an ordering or a
$D$-ordering. By Lemma~8.10, $T_{1}$ is a $D(J)$-ordering
contained in $T$. Taking any valuation $w$ on $L=F(\sqrt{t})$
extending $v$, we extend $\bar{T}$ to $\tilde{T}$ in $L_{w}$.  By Lemma~8.10
we obtain a $D(K)$-ordering $T_{2}$ in $L$. We enlarge it to a
$D(K')$-ordering $T'_{2}$, according to the three cases (1), (2), (3),
replacing $T_{2}\cup -T_{2}$ by $T_{2}$. As in the case for
$C(I)$-ordered fields, the only serious change is in proving that $\dot
L=T'_2\dot F$, and the proof is identical to that for $C(I)$-ordered fields.
\qed
\enddemo
\head \S 9. Galois groups and additive structures (2)\endhead

Throughout this paper, we have considered a number of subgroups $H$
of $\wg$  which behave pretty well, in that we have a certain
control over the additive structure of the associated orderings, and we
are able to make closures. Actually some of these groups $H$ have an
additional property which helped us in a subtle but important way. Let us
introduce the following definition and notation.

\definition{Definition and Notation 9.1}

(1) We say that an essential subgroup $H$ of $\wg$ is {\it
liftable} if we can write $\wg=G\rtimes H$ for some normal subgroup
$G$ of $\wg$.  This means that $H$ is not only a subgroup of $\wg$,
but also a quotient $\wg\lra H$ such that $H\lra\wg\lra H$ is the
identity map. The $H$-ordering $P_{H}$ is called a {\it liftable}
ordering.  (The name {\it liftable} was chosen because such an $H$
corresponds, as a quotient of $\wg$, to a Galois extension of $F$
inside $F^{(3)}$, of group $H$, which can be lifted as a Galois subextension of
$F^{(3)}$ of same group $H$.)

(2) If we want to realize some subgroup $H$ of $\wg$ as a $\Cal G_K$ for
some field $K$, we certainly need to use an $H$ which satisfies known
properties of $W$-groups. In particular, if $H\neq\{1\},C_2$, then by
Corollary~2.18 of \cite{MiSp2}, we see that $H$ can be embedded in a
suitable product $\prod_{I}D \times\prod_{J}C_4$, where each
factor is a quotient of $H$. According to the use in universal
algebra, see e.g. \cite{Gr, p. 123}, we refer to this as
the {\it subdirect
product} condition.

(3) We say that an essential subgroup $H$ of $\wg$ is a {\it
fair} subgroup if it is liftable and if it is either $\{1\}$ or $C_2$
or a subdirect product of some $\prod_{I}D \times\prod_{J}C_4$.  The
$H$-ordering $P_H$ will be called a {\it fair} ordering
if $H$ is a fair subgroup of $\wg$.

\enddefinition

  \definition{Remark 9.2} We observed in Example 6.4
that the subgroup $H = \left< \sigma,\tau \right> \cong C_4 * C_4$ in
$\Cal G_{\qq_2}$ has associated $H$-ordering $T = \dot F^2 \cup 5 \dot
F^2, F = \qq_2$, such that $T+T$ is not multiplicatively closed.  We
now use the description of $\wg = \Bbb G_2$ as in Example 2.9,
to show that $H$ is ${\underline{not}}$ a liftable subgroup of $\wg$.
Suppose instead that $H$ is a liftable subgroup of
$\wg$.  Then there exists a subgroup $G$ of $\wg$ such that $\wg = G
\rtimes H$.  Then $G$ must contain some element of the form $\alpha =
\rho \times h \times \phi$ where $\rho$ is an element of $\wg$ such
that $\rho, \sigma, \tau$ generate $\wg$ and $\sigma^2[\rho,\tau] = 1,
h \in H$ and $\phi$ is some element in $\Phi(\wg)$.
Because $G$ is a normal subgroup of $\wg$ we see that $\alpha \in G$
implies $\alpha^{-1}(\tau^{-1} \alpha \tau) = [\alpha, \tau] \in G$ as
well.  Hence $[\alpha, \tau] = [\rho h \phi, \tau] =
[\rho, \tau][h, \tau] \in G$.  On the other hand $[\rho, \tau][h,
\tau] = \sigma^2[h, \tau] \in H$.  Because $\wg = G \rtimes H$ we see
$G \cap H = \{1\}$ and thus $\sigma^2[h, \tau] = 1$. This
equality is impossible as $H$ is a free group in  category $\Cal
C$.  Therefore $H$ is not liftable.
\enddefinition

Observe that it is sometimes fairly easy to establish the ``fairness''
of a given subgroup. For example if $H=\<\sigma\>$ is an essential
subgroup of $\wg$ of order $2$, then for $f\notin P_{H}$ the
restriction $H\lra\Gal(F(\sqrt{f})/F)$ induces an isomorphism. Since
the subdirect product condition is empty, $H$ is fair.
We can also readily check the following:

\proclaim{Proposition 9.3} Let $\varphi\colon D(I)\lra\wg$ be an essential
embedding. Then $\varphi(D(I))$ is a liftable subgroup of $\wg$. As the
subdirect product condition is also trivially satisfied, it is a fair subgroup
of $\wg$.
\endproclaim

\demo{Proof}
Consider  a $D(I)$-ordering $T$ of $F$ for some $\mid I \mid \geq
1$.  Pick a basis for $\dot F/T$ of the form $\{[-1]\} \cup
\{[a_i], i \in I \}$.  (As usual $[f]$ means the class represented by
$f$ in the factor group $\dot F/T$.) Set $K/F = F(\sqrt{-1},
\root 4 \of a_i \colon i \in I)$.  Then  $\Gal(K/F)
\cong (\prod_I C_4) \rtimes C_2$, where we can choose 
generators $\bar{\tau}_i, \, i \in I$ for factors in the inner product
and $\bar{\sigma}$ for the outer factor such that
$\bar{\sigma}(\sqrt{-1}) = -\sqrt{-1}, \,\,\bar{\sigma} \, (\root 4
\of a_i) \, = \, \root 4 \of a_i, \, \bar{\tau}_i \, (\sqrt{-1}) \, =
\, \sqrt{-1}$ and $\bar{\tau}_i \, (\root 4 \of a_i) \, = \, \sqrt{-1}
\, \root 4 \of a_i, \, \bar{\tau}_i \, (\root 4 \of a_j) \, \, = \root
4 \of a_j$ for $j \neq i$.
Moreover the action of $\bar{\sigma}$ on $\prod_I C_4$ is
described as $\bar{\sigma}^{-1} \; \, \bar{\tau}_i \; \, \bar{\sigma}
\,=\, \bar{\tau}_i^3$ for each $i \in I$.  (Or equivalently
$\bar{\sigma}^{-1} \; \, \bar{\tau} \; \, \bar{\sigma} \, = \,
\bar{\tau}^{-1}$ for each $\bar{\tau} \in \prod_I C_4.)$

Pick any elements $\sigma, \; \tau_i, \; i \in I \in H\colon =
\varphi(D(I))$
such that
their homomorphic image from $H$ to $\Gal(K/F)$ are elements
$\bar{\sigma}, \; \bar{\tau}_i, \; i \in I$.  This is possible as $H$
surjects on $\Gal(K/F)$.
Then the essential subgroup $H$ of $\wg$ is generated by the
minimal set of generators $\{\sigma, \, \tau_i, \, i \in I\}$.
Moreover the natural restriction map $r\colon H \longrightarrow
\Gal(K/F)$ is an isomorphism, as $r$ takes the generators of $H$ to
the generators of $\Gal(K/F)$ and both sets of generators satisfy the same
relations.\qed
\enddemo

Now we  consider $C_4$-orderings and determine when they
are fair orderings.  Observe that  a $C_4$-ordering is
automatically fair provided it is liftable, so it is enough to
decide when a $C_4$-ordering $T$ is liftable.

\proclaim{Proposition 9.4} Let
$T$ be a $C_{4}$-ordering of $F$.  Then $T$ is liftable if and only if
there exists an element $f \in (F^2 + F^2) \setminus (T \cup \{0\})$.
\endproclaim

\demo{Proof} Suppose
that $T$ is a $C_4$-ordering of $F, T=P_{H}$ for $H \cong C_{4}$, and $H$
is essentially embedded in $\Cal G_F$.  Suppose also that $f \in (F^2
+ F^2) \setminus (T \cup \{0\})$.  Then since $f \notin T$ and $T \supset \dot
F^2$, we see that $f \notin F^2$ and a $C_{4}^{f}$-extension $K$ of
$F$ exists.  Because $f \in \dot F \setminus T$, an element $h \in H$ exists
such that $h(\sqrt{f}) = -\sqrt{f}$.  Then the image of
$h$ in $\Gal(K/F)$ under the natural homomorphism $H \longrightarrow
\Gal(K/F)$ is a generator of $\Gal(K/F)$.  Therefore the homomorphism is
in fact an isomorphism, and $H$ is liftable as asserted.
Assume now that $H \cong C_{4}$ is a liftable subgroup of
$\Cal G_F$.  Then a surjective homomorphism $\varphi \colon \Cal G_F
\longrightarrow C_4$ exists, which induces an isomorphism
$\psi \colon H \longrightarrow C_4$.  Let $K$ be the fixed
field of the kernel of $\varphi$.  Then $K/F$ is a Galois extension on
$\Gal(K/F) \cong C_4$.  Let $F(\sqrt{f})$ be a unique
quadratic extension of $F$ contained in $K$.  Also let $T = P_H$.
Then $H$ acts nontrivially on $f$ and $f \in (F^2 + F^2) \setminus
\{0\}$.  Hence $f \in (F^2 + F^2) \setminus (T \cup \{0\})$ as
claimed.  \qed
\enddemo

\definition{Example 9.5} The following simple example shows that we
cannot drop the condition $\exists \, f \in (F^2 + F^2) \setminus (T \cup
\{0\})$ from the proposition above, and that unfair $C_4$-orderings
exist in nature.  Consider again $F = \Bbb Q_2$ and set $T =
(F^2 + F^2) \setminus \{0\}$.  Then $T$ is a subgroup of $\dot F$ of index
$2$.  Because $\Bbb Q_2$ is not a formally real field, $\Bbb Q_2$ does
not admit any usual ordering, and $T$ is a $C_4$-ordering of $F$.
However $T$ contains all sums of two squares, and therefore $T$ is not
liftable.

On the bright side, we wish to point out that for each $C_4$-ordering
there exists a quadratic extension of the base field, and an extension
of the original $C_4$-ordering on this quadratic extension where this
extended ordering become a fair ordering.  In other words an unfair
ordering may become fair in some algebraic extension.  More precisely
we have the following proposition, in which we use Definition~1.4(4) of an
$H$-extension
\enddefinition

\proclaim{Proposition 9.6} Let $T$ be a $C_4$-ordering in $F$.  If $T$
is not fair, there exists  $t \in T$ and a $C_4$-extension
$(F(\sqrt{t}),V)$ of $(F,T)$ such that $V$ is a fair ordering in $F(\sqrt{t})$.
\endproclaim

\demo{Proof}
Suppose that $T$ is a $C_4$-ordering in $F$.  Then by Proposition~5.4 there
must exist an element $t \in T$ such that $1 + t
\notin T$.   If $T$ is not a fair
ordering, we know from the characterization of fair orderings in
Proposition~9.4 that
$t \notin \dot F^2$.  Hence $K = F(\sqrt{t})$ is a quadratic extension of $F$
and $[K \colon F] = 2$.  From the proof of Proposition 4.2, we know that
there exists some subgroup $V$ in $K$ such that $\mid \dot K/V \mid = 2$ and
$V \cap \dot F = T$.  Then $V$ is a $C_4$-ordering of $K$, and $V$ is fair as
$1+(\sqrt{t})^2 \notin V$.\qed
\enddemo

In this section we merely give a few examples of fair
orderings and are not pursuing a systematic check of which
orderings considered in this paper are fair and which will
become fair after extension to a suitable $2$-extension of the base
field.  The development of a theory of
fair orderings of fields is planned for a subsequent paper.

We complete our family of examples of orderings by considering $H =
\Cal F(I)$, where $I$ is some nonempty index set and $\Cal F(I)$ is
the
free pro-$2$-group in the category $\Cal C$, on a minimal set
$\{\sigma_i \mid i \in I \}$ of generators $I$.  (We assume as usual
that each open subgroup $V$ of $\Cal F(I)$ contains all but finitely
many $\sigma_i, i \in I$.  See \cite{Koc, Chapter 4}.)

\proclaim{Proposition 9.7}
Let $K/F$ be a Galois extension such that $\Gal(K/F) \cong \Cal F(I) =
\left< \sigma_i \mid i \in I \right>$ where $\{\sigma_i, i \in I\}$ is
a family of minimal generators of the free pro-$2$-group
$\Cal F(I)$.  Then there exists a
fair $\Cal F(I)$-ordering in $F$.
\endproclaim

\demo{Proof}
We  first embed the group $\Cal F(I)$ essentially in $\Cal G_F$.
Since $F^{(3)}$ is the maximal Galois subextension of a quadratic
closure $F_q$ of $F$ such that $\Gal(F^{(3)}/F)$ belongs to the category
$\Cal C$, and since $\Cal F(I)$ also belongs to $\Cal C$, we see that $K
\subset F^{(3)}$.  Therefore there exists a surjective natural
homomorphism $\pi\colon \Cal G_F \longrightarrow \Gal(K/F)$.

It is well known that there exists a continuous map $s \colon \Gal(K/F)
\longrightarrow \Cal G_F$ such that $\pi\circ s$ is the identity map
on $\Gal(K/F)$ (See \cite{Koc, 1.3}).  (Here we use only the fact that both
groups $\Gal(K/F)$ and $\Cal G_F$ are profinite groups.)  Set
$s(\sigma_i) = \omega_i$ for each $i \in I$.  Then for each open
subgroup $V$ of $\Cal G_F$ the set $s^{-1}(V)$ is an open subset of
$\Gal(K/F)$, and because open subgroups of $\Gal(K/F)$ form a basis for
the topology of $\Gal(K/F)$ we see that all but finitely many
$\sigma_i, i \in I,$ are in $\sigma^{-1}(V)$.  Hence all but finitely
many $\omega_i$ are in $V$.

Because $\Cal F(I)$ is a free object of $\Cal C$ on the set of
generators $(\sigma_i), i \in I$ we see that there exists a continuous
homomorphism $p \colon \Gal(K/F) \longrightarrow \Cal G_F$ such that
$p(\sigma_i) = \omega_i$ for each $i \in I$.
Set $H = p(\Gal(K/F))$.  Then we have $\pi \circ p = 1$ and $\Cal G_F
\cong \ker\pi\rtimes H$.  Moreover, $\pi$ restricted to $H$
induces an isomorphism $\varphi \colon H \longrightarrow \Gal(K/F)$.
Observe that $\varphi(\omega_i) = \sigma_i$ for each $i \in I$.
Because $\sigma_i \mod{\phi}(\Gal(K/F))$ are topologically independent,
we see that $\omega_i$ must be topologically independent
$\mod{\phi}(\Cal G_F)$.  This means that $\{\omega_i, i \in I\}$
generates the essential subgroup $H$ of $\Cal G_F$.

One can check that $\Cal F(I)$ is a subdirect product of its
dihedral and $C_4$ quotients
directly from the structure of $\Cal F(I)$, but it is also possible 
simply to observe that $\Cal F(I)$ is the $W$-group of a suitable field $A$
and all $W$-groups have this property.  That each $\Cal F(I)$
is the $W$-group of a suitable field $A$ follows from the fact that for
each index set $I \neq \phi$ we can find a field $A$ such that the
Galois group of its quadratic closure is a free pro-$2$-group (see
e.g., \cite{GM, page 98}), and therefore its $W$-group is $\Cal F(I)$.\qed
\enddemo

The following corollary applies, for example, in the case of $F = \Bbb
Q_p(t)$.

\proclaim{Corollary 9.8}
Let $F$ be the quotient field of a complete local integral domain
properly contained in $F$.  Let $\Cal F(I)$ be any free object of
category $\Cal C$ on generators $I$, where $I$ is a nonempty finite
set.  Then $F$ admits a fair $\Cal F(I)$-ordering.
\endproclaim

\demo{Proof}
 From Proposition~9.7 we see that it is sufficient to show that
each group $\Cal F(I), I$ finite and nonempty, occurs as a Galois
group over $F$.  Harbater's well-known result  \cite{Har, p. 186}
says that each finite group is realizable over $F$.  (For a very nice
and rather elementary proof of this result see \cite{HaV\"ol, Theorem
4.4}.)\qed
\enddemo

Let us fix the following notation.
\definition{Notation~9.9}
Let $i\colon F_{1}\lra F_{2}$ be a quadratic extension and let
$i^{\star}\colon\wgb\lra\wga$ be the associated restriction map. (See
e.g. \cite{MiSm3} for the existence of this map.)  Let
$H_{2}$ be a subgroup of $\wgb$ and let $H_{1}=i^{\star}(H_{2})$.
   Assume $H_{1}$ is essential in $\wga$. Observe that this property is
   not automatically satisfied since the image of an essential group
   under the restriction map $i^{\star}$ need not be essential.
(See Remark 7.8 for an example exhibiting such a case.)
When this is the case, we say that the extension
   $(F_{1},H_{1})\lra(F_{2},H_{2})$ is essential.
Put $T_{1}=P_{H_{1}}, T_{2}=P_{H_{2}}$.  Then it follows that
$T_{1}=T_{2}\cap F_{1}$.
\enddefinition

If we are working with fair groups  $H$ as above, then we can show that
for an essential quadratic extension $(F_{1},H_{1})\lra(F_{2},H_{2})$,
the additive structure of the associated orderings is preserved \ifff
$i^{\star}$ induces an isomorphism between $H_{2}$ and $H_{1}$.

\proclaim{Theorem 9.10} Assume the hypotheses in Notation~9.9 hold and that
$H_{1},H_{2}$ are fair subgroups of $\wga,\wgb$ respectively.  Then the
restriction $i^{\star}$ induces an isomorphism between
$H_{2}$ and
$H_{1}$ if and only if $\dot F_{1}/T_{1}\cong\dot F_{2}/T_{2}$ and for
each $a\in F_{1}, T_{1}+aT_{1}=(T_{2}+aT_{2})\cap F_{1}$.
\endproclaim
\comment
Observe that under the assumption $\dot
F_{1}/T_{1}\cong\dot F_{2}/T_{2}$, the condition
$T_{1}+aT_{1}=(T_{2}+aT_{2})\cap F_{1}$ is just another way of writing
$i(T_{1}+aT_{1})=T_{2}+aT_{2}$. Indeed, since $i(T_1 +a T_1)\subset T_2 +a
T_2$ it is enough to prove that $T_2 +a T_2\subset i (T_1 +a T_1)$. Assume
therefore that an element $0\neq f_2=t_2+as_2\in T_2 +a T_2$, where $t_2,s_2
\in T_2\cup\{0\}$. Since $\dot F_1/T_1 cong \dot 
F_2/T_2$,  there exists an element 
$f_1\in \dot F_1$ such that $f_2 = f_1 u_2$, where $u_2 \in T_2$. 
Hence $f_1=f_2 u_{2}^{-1}=t_2 u_{2}^{-1} +a s_2 
u_{2}^{-1}\in (T_2 +a T_2)\cap F_1$. From our assumption we have 
that $(T_2 +a T_2)\cap F_1 = T_1 +a T_1$. Thus 
$f_1 \in T_1 +a T_1$. Finally $i(f_1)=f_2$. Hence 
$i(T_1 +a T_1)=T_2 +a T_2$ as desired. 
\endcomment
Since the proof is a bit long and since the
two directions are not using the same assumptions on $H_{1},H_{2}$, we split
the theorem in two parts, Proposition~9.11  and Proposition~9.12

\proclaim{Proposition 9.11} Assume that $H_{1}$ is liftable.  Following
Notation~9.9, if the restriction $i^{\star}$ induces an
isomorphism between $H_{2}$ and $H_{1}$, then $\dot
F_{1}/T_{1}\cong\dot F_{2}/T_{2}$ and for each $a\in F_{1},\quad
T_{1}+aT_{1}=(T_{2}+aT_{2})\cap F_{1}$.
\endproclaim
\demo{Proof} We know that $\dot F_{i}/T_{i}$ is the Pontrjagin dual
of $H_{i}/\Phi(H_{i})$ for $i=1,2$. Thus the natural isomorphism
$H_{2}\lra H_{1}$ yields an isomorphism $\dot F_{1}/T_{1}\cong\dot
F_{2}/T_{2}$.
In order to show that for each $a\in F_{1}$ we have
$T_{1}+aT_{1}=(T_{2}+aT_{2})\cap F_{1}$, it is enough to show that for
every $b,c\in \dot F_{1}\setminus T_1$, if there exists $s_{2},t_{2}\in
T_{2}$ such that $bs_{2}+ct_{2}=1$, then there exists  $s_{1},t_{1}\in
T_{1}$ such that $bs_{1}+ct_{1}=1$.
Indeed, assume that the latter condition involving $b,c \in \dot{F}_1
\setminus T_1$ is valid.  Consider any $a \in \dot{F}_1$ and any relation
$u_2 + av_2 = d$, where $u_2,v_2 \in T_2 \cup \{0\}$ and $d \in \dot{F}_1$.
We want to show that there exist elements $u_1,v_1\in T_1 \cup
\{0\}$ such that $u_1 + av_1=d$.  If $u_2=0$ then $v_2 \in \dot{F}_1
\cap T_2 = T_1$, and we are done.  If $v_2=0$ then $u_2=d \in \dot{F}_1
\cap T_2=T_1$, and again we are done.  Then assume
$u_2,v_2 \in T_2$.  If $-a \in T_1$, let us write $d=s^{2}-t^{2}$ for some
 elements $s,t\in\dot F_{1} $. We then have
 $d=s^{2}+a(-at^{2}/a^{2})\in T_{1}+aT_{1}$. 
 Hence we may assume that $-a \notin
T_1$.  Finally we also assume that $d \notin T_1$. From the equation
$u_2+av_2=d$ we obtain $u_2=d-av_2$, and since
$u_2,v_2\in T_2$ we can rewrite this equation as $1=ds_2-at_2$ where
$d,-a\in\dot{F}_1 \setminus T_1$.  Using our hypothesis we see that there exist
elements $s_1,t_1\in T_1$ such that $1=ds_1-at_1$.  Hence $d \in
T_1+aT_1$ as required.

Now take $b,c\in \dot{F}_{1}\setminus T_1$
and assume that $bs_{2}+ct_{2}=1$ for
some $s_{2},t_{2}\in T_{2}$. Then the quaternion algebra
$\left(\frac{bs_{2},ct_{2}}{F_{2}}\right)$ splits. We consider the
following cases.

(1) Suppose $bs_{2},ct_{2}$ are linearly independent in $\dot F_{2}/T_{2}$.
Then they are also independent modulo $\dot F_{2}^{2}$, and by
Proposition~1.5 we  have a dihedral
extension  $L_2/F_{2}$ such that
$F_{2}(\sqrt{bs_{2}},\sqrt{ct_{2}})\subset L_2$ and
$\Gal(L_2/F_{2}(\sqrt{bcs_{2}t_{2}}))\cong C_4$. In particular we have
an exact sequence
$$1\lra C_2\lra\Gal(L_2/F_{2})\cong
D\lra\Gal(F_{2}(\sqrt{bs_{2}},\sqrt{ct_{2}})/F_{2})\cong C_2\times
C_2\lra 1.$$

Let $\theta$ denote the restriction map from $H_{2}$ to
$\Gal(F_{2}(\sqrt{bs_{2}},\sqrt{ct_{2}})/F_{2})$.  We show it is
surjective.
Denote by $u_{1},u_{2}$ the two generators of
$\Gal(F_{2}(\sqrt{bs_{2}},\sqrt{ct_{2}})/F_{2})$ defined by
$u_{1}(\sqrt{bs_{2}})/\sqrt{bs_{2}}=-1,
u_{1}(\sqrt{ct_{2}})/\sqrt{ct_{2}}=1, u_{2}(\sqrt{bs_{2}})/\sqrt{bs_{2}}=1,
u_{2}(\sqrt{ct_{2}})/\sqrt{ct_{2}}=-1$. We may look at $u_{1},u_{2}$ as
linear functions on the $\F_{2}$-vector subspace of $\dot F_{2}/\dot
F_{2}^{2}$ spanned by $bs_{2},ct_{2}$, which are assumed to be
independent, and since $bT_{2}\cap cT_{2}=\emptyset$, we may  extend them
to linear functions $v_{1},v_{2}$
defined on the subspace generated by $bT_{2},cT_{2}$, by putting
$v_{i}(x)=u_{i}(b)$ if $x\in bT_{2}$ and $v_{i}(x)=u_{i}(c)$ if $x\in
cT_{2}$. Then $v_{i}$ may be viewed as a function on the
$\F_{2}$-vector subspace generated by the cosets $bT_{2},cT_{2}$ in $\dot
F_2/T_{2}$.
Again, these functions $v_{i}$'s may be extended to $w_i$ defined on
the whole vector
space $\dot F_2/T_{2}$. By duality, one has $(\dot F_2/T_{2})^{\star}\cong
H_{2}/\Phi(H_{2})$, and the $w_{i}$'s yield to elements in
$H_{2}/\Phi(H_{2})$ which may be lifted as elements $h_{1},h_{2}\in
H_{2}$. Since the duality is precisely given by the pairing
$H_{2}/\Phi(H_{2})\times\dot
F_{2}/T_{2}\lra\{\pm 1\}$ defined by $(h,f)\mapsto
h(\sqrt{f})/\sqrt{f}$, it is immediate that $h_{i}$ goes to $u_{i}$
under the restriction map $\theta\colon
H_{2}\lra\Gal(F_{2}(\sqrt{bs_{2}},\sqrt{ct_{2}})/F_{2})$.
This shows the surjectivity of $\theta$. Since  $\theta$ factors
through $\psi\colon H_{2}\lra\Gal(L_2/F_{2})\cong D$ and since the kernel of
$\Gal(L_2/F_{2})\lra\Gal(F_{2}(\sqrt{bs_{2}},\sqrt{ct_{2}})/F_{2})$ is
the Frattini subgroup of $\Gal(L_{2}/F_{2})$, we see that $\psi$ is
also surjective.
This means that $D$ may be viewed as a quotient of $H_{2}$ and that we
have inclusion maps
${F_{2}^{(3)}}^{H_{2}}\lra L'_{2}\lra F_{2}^{(3)}$ such
that $\Gal(L'_{2}/{F_{2}^{(3)}}^{H_{2}})\cong D$. Since
$i^{\star}(H_{2})=H_{1}$, applying $i^{\star}$ to this diagram gives
us another diagram ${F_{1}^{(3)}}^{H_{1}}\lra L'_{1}\lra F_{1}^{(3)}$
with $\Gal(L'_{1}/{F_{1}^{(3)}}^{H_{1}})\cong D$.

Since $H_{1}$ is
liftable, we know that there exists an $H_{1}$-extension $K/F_{1}$
inside $F_{1}^{(3)}$ containing a $D$-extension $L_{1}/F_{1}$. This
extension is a $D^{u,v}$-extension for suitable $u,v\in F_{1}$ by
Proposition~1.5. We claim that we have
$u=bs_{1},v=ct_{1}$ for suitable $s_{1},t_{1}\in T_{1}$.
Consider the surjective homomorphism
$$\theta\colon H_2 \lra \Gal(F_2(\sqrt{bs_2},\sqrt{ct_2})/F_2)$$
defined above. This homomorphism factors through the
surjective homomorphism $\psi\colon H_2 \lra
\Gal(L_2/F_2)\cong D$. Using the isomorphism $\beta\colon H_2 \lra H_1$
induced by $i^{\star}$ and our construction of
$L_1/F_1$, we see that the homomorphism $\psi\colon H_2 \lra \Gal(L_2/F_2)$
is compatible, via identification of $H_2$ with
$H_1$ using $i^{\star}$, with the restriction homomorphism
$\tilde{\psi}\colon H_1 \lra \Gal(L_1/F_1)$. Passing to the
quotients $\Gal(F_2(\sqrt{bs_2},\sqrt{ct_2})/F_2)$ and
$\Gal(F_1(\sqrt{u},\sqrt{v})/F_1)$ of $\Gal(L_2/F_2)$ and
$\Gal(L_1/F_1)$ respectively, we see that we can identify the homomorphism
$\theta\colon H_2 \lra
\Gal(F_2(\sqrt{bs_2},\sqrt{ct_2})/F_2)$ with the
restriction homomorphism $\tilde{\theta}\colon H_1 \lra
\Gal(F_1(\sqrt{u},\sqrt{v})/F_1)$ via the isomorphism
$i^{\star}\colon H_2 \lra H_1$. Finally from the natural isomorphism
$\dot{F}_1/T_1 \cong \dot{F}_2/T_2$ we may assume that $u=bs_1$ and $v=ct_1$
for suitable elements
$s_1, t_1 \in T_1$.
By Proposition~1.5, this implies that the quaternion algebra
$\left(\frac{bs_{1},ct_{1}}{F_{1}}\right)$ splits, and that there
exist $\tilde{s_{1}},\tilde{t_{1}}\in T_{1}$ such that
$b\tilde{s_{1}}+c\tilde{t_{1}}=1$.

Suppose now that $bs_{2},ct_{2}$ are linearly dependent in
$\dot F_{2}/T_{2}$. Then $b$ and $c$ are equal modulo $T_{2}$ and we may
assume $b=c$. There are
still two more cases to consider.

(2) Suppose  we have $cs_{2}+ct_{2}=1$ with
$s_{2}=t_{2}\bmod \dot F_{2}^{2}$. By Proposition~1.5,
there exists a $C_{4}^{cs_{2}}$-extension $L_{2}/F_{2}$ with
$F_{2}(\sqrt{cs_{2}})\subset L_{2}$. Using arguments similar to those
in (1), we show that the restriction $\psi\colon
H_{2}\lra\Gal(L_{2}/F_{2})$ is onto, and we find $s_{1}\in T_{1}$ such
that $\left(\frac{cs_{1},cs_{1}}{F_{1}}\right)$ splits. This implies
that there exist $\tilde{s_{1}},\tilde{t_{1}}\in T_{1}$ such that
$c\tilde{s_{1}}+c\tilde{t_{1}}=1$.

(3) Suppose  we have $cs_{2}+ct_{2}=1$ with
$s_{2}\not=t_{2}\bmod \dot F_{2}^{2}$. As in (1) we find $L_{2}$
with $\Gal(L_2/F_{2})\cong D$ and we have a tower of fields
$F_{2}\lra F_{2}(\sqrt{s_{2}t_{2}})\lra F_{2}(\sqrt{cs_{2}},\sqrt{ct_{2}})\lra
L_{2}$. Since $H_{2}$ fixes $F_{2}(\sqrt{s_{2}t_{2}})$, the
restriction map $\psi\colon H_{2}\lra\Gal(L_2/F_{2})$
induces a surjective homomorphism
$\psi'\colon H_{2}\lra\Gal(L_2/F_{2}(\sqrt{s_{2}t_{2}}))\cong C_{4}$.
We finish with arguments as in (2) and replacing $F_{2}$ by
$F_{2}(\sqrt{s_{2}t_{2}})$, we find $\tilde{s_{1}},\tilde{t_{1}}\in T_{1}$
such that
$c\tilde{s_{1}}+c\tilde{t_{1}}=1$.
\qed
\enddemo

We now prove the result in the other direction.

\proclaim{Proposition 9.12}
Let $H_{1},H_{2}$ be as in Notation~9.9 and assume they are fair
subgroups. If the inclusion $i\colon F_{1}\lra F_{2}$ induces an
isomorphism $\dot F_{1}/T_{1}\lra \dot F_{2}/T_{2}$ and if
$(T_{2}+aT_{2})\cap\dot F=T_{1}+aT_{1}$ for any $a\in F_{1}$, then
$i^{\star}$ induces an isomorphism between $H_{2}$ and $H_{1}$.
\endproclaim

\demo{Proof} If $H_2=\{1\}$ then $H_1=\{1\}$ as well. If $H_2=C_2$ then
$i^{\star}(H_2)\neq\{1\}$ because $T_2$ is a
usual ordering in $\dot{F}_2$, and it cannot contain $\dot{F}_1$. However if
$H_1$ were $\{1\}$ then $T_1=\dot{F}_1$.
Therefore $i^{\star}$ induces an isomorphism between $H_2$ and $H_1$.

For the rest of our proof we assume that
$H_2\neq \{1\},C_2$.
Call $\beta\colon H_{2}\lra H_{1}$ the restriction of
$i^{\star}$ to $H_{2}$.
Because $i^{\star}$ is a group homomorphism from $\wgb$ into $\wga$, we have
$i^{\star}(\Phi(\wgb))\subset\Phi(\wga)$.
Also we have $\beta(\Phi(H_{2}))\subset\Phi(H_{1})$.
Then the map $\beta$ induces  $\hat{\beta}\colon H_{2}/\Phi(H_{2})\lra
H_{1}/\Phi(H_{1})$, which is an isomorphism because its dual map $\dot
F_{1}/T_{1}\lra\dot F_{2}/T_{2}$ is an isomorphism. By definition $\beta$
is onto.
We want to show that $\beta$ is injective.  From
the fact that $\hat{\beta}$ is an isomorphism,
we see that
$\ker\beta\subseteq\Phi(H_{2})$. Take a fixed set of minimal
(topological) generators $(\sigma_{i})_{i\in I}$ for $H_{2}$. Then
$\gamma\in \Phi(H_{2})$ has a unique description, up to a
permutation, as $\gamma=\prod_{i\in
I}\sigma_{i}^{2}\times\prod_{(u,v)\in K}[\sigma_{u},\sigma_{v}]$ for
some possibly infinite sets $I,K$.

To complete the proof we use the following lemma.

\proclaim{Lemma 9.13} Assume that $H_{1},H_{2},T_{1},T_{2}$ are as in
Proposition~9.12, and let $\delta$ be $\sigma_{i}^{2}$ or
$[\sigma_{u},\sigma_{v}]$. Suppose that we have a surjective map
$\varphi\colon H_{2}\lra G$ where $G=D$ or $C_{4}$.
Then there exists a group $\tilde{G}$ which is again either $D$ or $C_4$
and a homomorphism
$\psi\colon H_{1}\lra
\tilde{G}$  such that $\psi(\beta(\delta))\neq 1\in \tilde{G}$ \ifff
$\varphi(\delta)\neq 1 \in G$.
Moreover $\tilde{G}$ and the homomorphism $\psi$
depend only on $G$ and on the fields $F_1$ and $F_2$, but not on $\delta$.
\endproclaim
  \demo{Proof}
(1) Assume first that $G=C_{4}$.
  Since  $H_{2}$ is liftable, there exist an $H_{2}$-extension $K_{2}/F_2$
and a $C_{4}^{u}$-extension
   $L_{2}$ of $F_{2}$ with $F_{2}\lra F_{2}(\sqrt{u})\lra L_{2}\lra
   K_{2}$. Since $\dot F_{1}/T_{1}\cong\dot F_{2}/T_{2}$, there exist
   $a\in \dot F_{1}, s_{2}\in T_{2}$ such that $u=as_{2}$. Let
   $\delta=\sigma^{2}$, which is the only case to be considered when
   $G=C_{4}$. Then $\varphi(\sigma^{2})\neq 1\in \Gal(L_{2}/F_{2})$ if
   and only if $\varphi(\sigma)$ has order $4$.
Thus $\varphi(\sigma^2)\neq 1$
 \ifff $\varphi(\sigma)$
   generates $\Gal(L_{2}/F_{2})$. This happens precisely when
   $\varphi(\sigma)(\sqrt{as_{2}})=-\sqrt{as_{2}}$. Since $H_{2}$, and
   thus $\varphi(\sigma)$, fixes $\sqrt{s_{2}}$, this is equivalent to
$\varphi(\sigma)(\sqrt{a})=-\sqrt{a}$.
On the other hand, we know by Proposition~1.5 that the quaternion algebra
$\left(\frac{as_{2},as_{2}}{F_{2}}\right)$ splits, and this implies
the existence of $s'_{2},t'_{2}\in T_{2}$ such that $as'_{2}+at'_{2}=1$.
 From the assumption on the additive structure, this implies the
existence of $s_{1},t_{1}\in T_{1}$ such that $as_{1}+at_{1}=1$.
Two cases are  to be considered.

\noindent(1.1) If $s_{1}=t_{1}\bmod \dot F_{1}^{2}$, then there is a
$C_{4}^{as_{1}}$-extension $L_{1}$ of $F_{1}$ with $F_{1}\lra
F_{1}(\sqrt{as_{1}})\lra L_{1}$. Denoting by $\psi\colon
H_{1}\lra\Gal(L_{1}/F_{1})$ the restriction, because $H_{1}$ fixes
$\sqrt{T_{1}}$ we have
$\psi(\beta(\sigma))(\sqrt{as_{1}})/\sqrt{as_{1}}
=\psi(\beta(\sigma))(\sqrt{a})/\sqrt{a}=\varphi(\sigma)(\sqrt{a})
/\sqrt{a}=-1$,
  showing $\psi(\beta(\delta))\neq 1\in C_{4}=\tilde{G}$.

  \noindent(1.2) If $s_{1}\neq t_{1}\bmod \dot F_{1}^{2}$, then there is a
  $D^{as_{1},at_{1}}$-extension $L_{1}$ of $F_{1}$ with $F_{1}\lra
  F_{1}(\sqrt{s_{1}t_{1}})\lra L_{1}$.
Here $L_1/F_1(\sqrt{s_1 t_1})$ is a $C_4$-extension.
Since $\beta(\sigma)\in H_{1}$
  fixes $F_{1}(\sqrt{s_{1}t_{1}})$, $\psi(\beta(\sigma))$ is in the
  Galois group of the latter extension, which is again a
  $C_{4}$-extension.  We then use the same argument as in (1.1) to
  conclude that $\psi(\beta(\delta))\neq 1\in \tilde{G}=C_{4}$.

(2) Assume $G=D$. Again there is an $H_{2}$-extension $K_{2}\, of
F_{2}$ and
a $D^{as_{2},bs_{2}}$-extension
   $L_{2}$ of $F_{2}$ with $F_{2}\lra F_{2}(\sqrt{abs_{2}t_{2}})\lra L_{2}\lra
   K_{2}$. Since $\varphi$ is surjective, there is an element $\tau\in
   H_{2}$ such that
   $\tau(\sqrt{abs_{2}t_{2}})/\sqrt{abs_{2}t_{2}}=-1$, or else
   $\varphi(H_{2})$ would fix $F_{2}(\sqrt{abs_{2}t_{2}})$ and would be
   contained in a proper subgroup of $\Gal(L_{2}/F_{2})\cong D$. This
   implies  $ab\notin T_{2}$. Since there exist $s'_{2},t'_{2}\in
T_{2}$ such that
$as'_{2}+bt'_{2}=1$, we also have, by the assumption on the additive
structures, $as_{1}+bt_{1}=1$ for some $s_{1},t_{1}\in T_{1}$. Since
$ab\notin T_{1}$, we see that
$as_{1}, bt_{1}$ are independent modulo $\dot F_{1}^{2}$,
and there is a $D^{as_{1},bt_{1}}$-extension $L_{1}$ of $F_{1}$
with $F_{1}\lra F_{1}(\sqrt{abs_{1}t_{1}})\lra L_{1}$. Denote by $\psi\colon
H_{1}\lra\Gal(L_{1}/F_{1})\cong D$  the restriction map.

\noindent(2.1)  Suppose $\delta=\sigma^{2}$ and $\varphi(\delta)\neq 1$. Then
$\varphi(\sigma)$ has order $4$ and must fix the quadratic extension
$F_{2}(\sqrt{abs_{2}t_{2}})$. Then it belongs to
$\Gal(L_{2}/F_{2}(\sqrt{abs_{2}t_{2}}))\cong C_{4}$. With the
same arguments as in (1), we show that
$\psi(\beta(\delta))\neq 1$.

\noindent(2.2) Suppose $\delta=[\sigma_{u},\sigma_{v}]$ and
$\varphi(\delta)\neq 1$.
Then none of $\varphi(\sigma_{u}),\varphi(\sigma_{v})$ is
  in $\Phi(D)$ (i.e. they do not fix the biquadratic extension
  $F_{2}(\sqrt{as_{2}},\sqrt{bt_{2}})$), and  they act differently on
  this biquadratic extension.
Since $\varphi(\sigma_{u})$ (respectively $\varphi(\sigma_{v})$) acts the
same way on elements in $\sqrt{\dot F}$ as $\psi(\beta(\sigma_{u}))$
(respectively
$\psi(\beta(\sigma_{v}))$, we see that $\psi(\beta(\delta))\neq 1\in
G$. 

To conclude the proof, we point out
that in all cases above, we first associated $\tilde{G}$ with
the given homomorphism
$\varphi\colon H_2\lra G$ and only then checked that
$\varphi(\delta)\neq 1 \in G$ is equivalent to
$\psi(\beta(\delta))\neq 1 \in \tilde{G}$. \qed
\enddemo

We can now finish the proof of Proposition~9.12. Suppose
$\gamma\neq 1\in\Phi(H_{2})$.  Since $H_{2}$ satisfies the subdirect
product condition, there exists a surjective map $\varphi\colon
H_{2}\lra G$ with $G\cong D$ or $C_{4}$ and with $\varphi(\gamma)\neq 1\in
G$.  Recall that the minimal set of generators $(\sigma_{i})_{i\in I}$ may be
chosen in such a way that for any open set $U$ of $H_{2}$ there are at
most finitely many $\sigma_{i}$'s outside $U$.  (See for example
\cite{Koc, Chapter~4}.) Since $\ker\varphi$ is open, we may thus
assume, when working with a given $\varphi$, that
$\gamma=\gamma_{0}\times \gamma_{1}$, with $\gamma_{0}=\prod_{i\in
I_{0}}\sigma_{i}^{2}\times\prod_{(u,v)\in
K_{0}}[\sigma_{u},\sigma_{v}]$, $\gamma_{1}=\prod_{i\in
I_{1}}\sigma_{i}^{2}\times\prod_{(u,v)\in
K_{1}}[\sigma_{u},\sigma_{v}]$, with the following properties.  The
sets $I_{0}, K_{0}$ are finite.  Any individual factor
$\sigma_{i}^{2}, [\sigma_{u},\sigma_{v}]$ of $\gamma_{0}$ is not in
$\ker\varphi$, while any individual factor of $\gamma_{1}$ is in
$\ker\varphi$. We may assume that $\gamma = \gamma_{0}$, and in particular
we have only a finite number $n$ of terms $\delta_{i}$'s with
$\delta_{i}=\sigma_{i}^{2}$ or $[\sigma_{u},\sigma_{v}]$.  The
Frattini group $\Phi(G)\cong C_{2}$ may be written $\{1,\epsilon\}$,
and each $\varphi(\delta_{i})$ must be $\epsilon$, since it is not $1$
by assumption.  Since $\varphi(\gamma)=\epsilon^{n}\neq 1$, $n$ must
be odd.  By Lemma~9.13, we know that there exists a group $\tilde{G}$ which
is again $D$ or $C_4$  and a homomorphism
$\psi\colon H_1 \lra \tilde{G}$, such that
$\varphi(\delta_{i})=\epsilon\neq 1$ is equivalent to
$\psi(\beta(\delta_{i}))=\epsilon\neq 1$.
Because $n$ is odd, this
shows that $\psi(\beta(\gamma))\neq 1$, and therefore
$\beta(\gamma)\neq 1$.  This shows the injectivity of $\beta$ and
finishes the proof of Proposition~9.12.\qed
\enddemo

\head \S 10. Concluding Remarks
\endhead

In this article we have considered all $C(I)$- and $S(I)$-orderings.  These
groups
correspond to W-groups for $p$-adic fields, for odd primes $p$.  In
particular, the
W-group $\Cal G_p$ of $\qq _p$ is $C_4 \times C_4$ for $p \equiv 1 (4)$ and
is $C_4
\rtimes C_4$ for $p \equiv 3 (4)$.  It is then natural to look for a
characterization of $\Cal G_2$-orderings, i.e. those orderings corresponding to
subgroups isomorphic to the W-group of $\qq _2$.  This is currently under
investigation \cite{MiSm4}.

For the field $\qq$, there is a unique $C_2$-ordering, which is the
unique ordering on $\qq$.  In addition there is a one-to-one
correspondence between $C_4 \times C_4$-orderings on $\qq$ and primes
$p \equiv 1 (4)$, and a one-to-one correspondence between $C_4 \rtimes
C_4$-orderings on $\qq$ and primes $p \equiv 3(4)$.  In each case the
correspondence is given by $T_p = \dot\qq_p^2 \cap \qq$.  It is not hard
to see that each such intersection gives rise to an $H$-ordering of
the appropriate type.  To see that every such orderings may be
obtained in this way, one shows that each such ordering corresponds
to a certain valuation on $\Bbb{Q}$, and the valuations on $\Bbb{Q}$
are well-known to be classified by the primes. (See e.g. \cite{End, Theorem 1.16}.)

This observation then lends itself to an alternative perspective on
the Hasse-Minkowski Theorem, which states that a quadratic form
defined over $\qq$ is isotropic over $\qq$ if and only if it is
isotropic over each $\qq_p$, including $\qq_\infty$, the real numbers.
Using Hilbert's reciprocity law, one can prove that a ternary quadratic
form is isotropic over $\qq$ \ifff it is
isotropic over all but one of these fields.
Thus we see that
a ternary quadratic form over $\qq$ is isotropic if and only if it is
isotropic with respect to all $C_2$-, ($C_4 \times C_4$)-, and ($C_4
\rtimes C_4$)-orderings on $\qq$.

We point out that the case of a ternary quadratic form over
$\qq$, together with the clever use of Dirichlet's
theorem on the existence of an infinite number of primes in an
arithmetic progression, where first term and
increment are relatively prime, are the main ingredients of a proof of
the full Hasse-Minkowski theorem over $\qq$. For a very nice
exposition of the Hasse-Minkowski theorem over $\qq$, see \cite{BS}. See
also \cite{L1, Chapter 6, Exercise 22}.

It is not difficult however, to find a quaternary quadratic form $\varphi$
over $\qq$ such that $\varphi$ is isotropic
over all $\qq_p$, $p$ is an odd prime, and $\qq_{\infty} = \Bbb{R}$ but
$\varphi$ is anisotropic over $\qq_2$. Because
we were unable to locate an explicit example of such a form in the
literature, we write down one explicit
example here:
$$\varphi=X_1^2 + X_2^2 -7 X_3^2-31 X_4^2$$
Because $-7\equiv 1 (\mod 8)$ and $-31\equiv 1 (\mod 8)$, we see that
$\varphi$ is equivalent to $\psi = X_1^2+X_2^2+X_3^2
+X_4^2$ over $\qq_2$. Because the level of $\qq_2$ is $4$, we see that
$\psi$ is anisotropic over $\qq_2$. On the other
hand using the well-known Springer theorem for local fields (\cite{L1,
Chapter 6, Proposition 1.9}), and the fact that each
ternary quadratic form over any finite field is isotropic, we see that
$\varphi$ is isotropic over each $\qq_p$, $p$
an odd prime. Since $\varphi$ is an indefinite quadratic form, $\varphi$ is
isotropic over $\qq_{\infty}$ as well.

In a subsequent paper we will present several applications of this theory
to  different kinds of local-global principles for quadratic forms. In
order to get a sense of what can be done in this direction, we
show below an example of a simple situation in which our theory
applies.

Consider a field $F$.  Recall that a $C(\emptyset)$-ordering $T$ on $F$ is an
index $2$ multiplicative subgroup of $\dot F/\dot{F}^2$ containing  $-1$.
Additively speaking, it is a hyperplane containing $-1$ in the
$\F_{2}$-vector space  $\dot{F}/\dot{F}^2$. If $f\in\dot{F}\setminus
(\dot{F}^2\cup-\dot{F}^2)$ and if $V$ is any subspace of
$\dot{F}/\dot{F}^2$ such that
$\dot{F}/\dot{F}^2=\Span\{f,-1\}\oplus V$, then
$T:=\Span\{-1\}+V$ is a $C(\emptyset)$-ordering not containing $f$. Then
the next
lemma follows immediately.

\proclaim{Lemma 10.1} Let $C_{0}(F)$ denote the set of
$C(\emptyset)$-orderings of $F$.  Then $C_{0}(F)=\emptyset$ if and
only if $\dot F=\dot F^2\cup -\dot F^2$,  and in general,
$$\bigcap_{T\in C_{0}(F)}T= \dot{F}^2\cup -\dot{F}^2.$$
\endproclaim

To every $C(\emptyset)$-ordering $T$ we associate a fixed closure
$F_{T}$ of $F$ in the quadratic closure of $F$.  Denote by
$\<\<a_{1},\ldots,a_{n}\>\>$ the Pfister form $\<1,-a_{1}\>\otimes
\ldots\otimes\<1,-a_{n}\>$. (For the basic theory of Pfister forms see e.g.
\cite{L1, Chapter 10} or \cite{Sc, Chapter 4}.
Observe that both Lam and Scharlau denote by $\<\<a_{1},\ldots,a_{n}\>\>$
the Pfister form
$\<1,a_{1}\>\otimes\ldots\otimes\<1,a_{n}\>$.)
Then we have the following.

\proclaim{Proposition 10.2} Assume $C_{0}(F)\not=\emptyset$ and denote by
$\varphi$ the  map $W(F)\lra\prod_{T\in C_{0}(F)}W(F_{T})$ induced by
the inclusions $F\lra F_{T}$.
Then $\Ker\varphi=I^2F+2W(F)$ where $IF$ denotes the fundamental ideal
of $W(F)$.
\endproclaim
\demo{Proof} For $T\in C_{0}(F)$ we have $\dot{F}/T=\{\bar{1},\bar{f}\}$
for a certain $f\in\dot{F}$, and it is easy to see that
$W(F_{T})\cong C_2 [\epsilon]/\epsilon^2$ and that the isomorphism,
say $\pi$, is defined by $\pi(\< \bar{1}\>)=1,\pi(\<
\bar{f}\>)=1+\epsilon$. If $a,b\in\dot{F}$
then the possibilities for $\bar{a},\bar{b}$ are
(1) $\bar{a}=1$ or $\bar{b}=1$, or (2) $\bar{a}=\bar{b}=\bar{f}$.
In any case the image in $W(F_{T})$ of the $2$-fold Pfister form $\<\<
a,b\>\>$ is in $2W(F_{T})=0$, and we have shown the inclusion
$I^2F+2W(F)\subseteq\Ker\varphi$.

Take $q\in\Ker\varphi$.  Then $q\in IF$, because any odd-dimensional
form is nonzero in $W(F_{T})$.  But it is known (\cite{Pf, p. 122, Kor.
to Satz 13}) that any element $q$ of $IF$ may be written
$q=\<\< u\>\>+ q_{1}$, with $q_{1}\in I^2F$.  Since $q\in\Ker\varphi$,
and $I^2F\subset\Ker\varphi$, we deduce $\<\< u\>\>\in\Ker\varphi$.
The latter is equivalent to $u\in T$ for every $T$, meaning
$u\in\dot{F^2}\cup -\dot{F^2}$, or in other words $\<\<u\>\>=0$ or $2$
in $W(F)$.\qed
\enddemo

Recall that a field $F$ is said to have virtual cohomological
dimension $n$, denoted $\vcd(F)=n$, if
$H^{d}(\Gal(F(2))/F(\sqrt{-1})), \mu_{2})=0$ for $d>n$, and
$H^{n}(\Gal(F(2))/F$ $(\sqrt{-1})), \mu_{2})\neq 0$.
(If we also considered the case of $\Bbb{F}_p, p$ an odd prime,
as coefficients of the
cohomology groups of absolute Galois groups, it would be more appropriate
to say that $F$ as above has virtual
$2$-cohomological dimension equal to $n$.)
If $\vcd(F) \le 1$, then $I^{2}F$ is torsion-free. To see this, observe
first that $\vcd(F) \le 1$ implies each binary quadratic
form over $F(\sqrt{-1})$ is universal. Then use \cite{L1, Chapter 11, Theorem
1.8 and Exercise 20}
to conclude that $I^{2}F$ is torsion free.
An example of a formally real field $F$ with $vcd(F)=1$ is $F=\Bbb{R}(X)$.
We have the following local-global
principle:

\proclaim{Theorem 10.3} Let $F$ be a field with $\vcd(F)\le 1$.  Let
$D_{0}(F)$ (resp.  $C_{0}(F)$, $S_{0}(F))$ denote the set of usual
orderings $X(F)$ (resp.  $C(\emptyset)$-orderings, $S(\emptyset)$-orderings)
of $F$.  Then $$\Lambda\colon W(F)\lra\prod_{T\in D_{0}(F)\cup
C_{0}(F)\cup S_{0}(F)}W(F_{T})$$ is injective.  If $F$ is formally
real, we may drop $S_{0}(F)$.  (If not, we may of course drop
$D_{0}(F)$.)
\endproclaim

\demo{Proof}
It is clear that a form $q\in\Ker\Lambda$  is in $IF$, and thus can be
written $q=\<\< a\>\> +q_{2}$ with $q_{2}\in I^2F$. By Pfister's
Local-Global Principle \cite{L1, Chapter 8, \S 4},
$q$ is torsion and
it is therefore the case for $\<\< a\>\>$ and $q_{2}$. (It is trivial when
$D_{0}(F)=\emptyset$, and if not, we use the fact that the signature
$\hat{q}$ of
$q$ is $0$
  and that $\hat{q_{2}}\equiv 0(\bmod 4)$.)

Since $I^2F$ is torsion-free, one has $q_{2}=0$, and $q=\<\< a\>\>$.
Since $q$ vanishes on $C_{0}(F)$, by Proposition~10.2 we have
$a\in\dot{F}^2\cup -\dot{F}^2$.  (If $C_{0}(F)=\emptyset$, this
condition is trivially satisfied.)  If the level $s(F)$ is $1$, we are
finished, and otherwise $D_{0}(F)\cup S_{0}(F)\neq\emptyset$, showing that
$q\neq\<\< -1\>\>$.  Thus $q=\<\< 1\>\>=0$.
\qed
\enddemo

\definition{Remark 10.4}
In this case we even have a strong Hasse Principle, that is a
local-global principle for detecting whether a quadratic
form is anisotropic rather than just hyperbolic. Indeed, the fact that
each ternary form over $F(\sqrt{-1})$ is isotropic and
\cite{ELP, Theorem F} give us the strong Hasse Principle for forms of rank
greater than or equal to $3$.  Then
the use of $C_{0}(F),S_{0}(F)$ and $D_{0}(F)$ provides the result for rank $2$
forms.
\enddefinition

Finally let us point out that our results are closely related to some
ideas in birational anabelian Grothendieck geometry.  In particular
there is a close connection between ideas explored in this paper and
the work of  Bogomolov,  Tschinkel and Pop (\cite{Bo},
\cite{BoT}, \cite{Po1}, and \cite{Po2}; see also Koenigsmann's
thesis~\cite{K1} and paper~\cite{K2}).
Roughly speaking, they establish that for certain fields $K$ the
isomorphy type of $K$, modulo purely inseparable extensions of $K$, is
functorially encoded in the pro-$p$-quotient of the absolute Galois
group $\tilde G : = \Gal(\bar{K}/K), \ch K \neq p$.  In fact Bogomolov
in~\cite{Bo} and also Pop in  lectures at MSRI in the fall of 1999,
considered smaller Galois groups than the Galois group defined above, namely
the maximal pro-$p$-quotient of the group
$\tilde G/[[\tilde G, \tilde G], \tilde G]$.
In this paper we consider $p=2$, because of the connections with
quadratic forms.  It is  expected however that a substantial part of
our results can be extended to any prime $p$ provided that the base
field $F$ contains a primitive $p$th root of unity.  We
allow $F$ to be any field with $\ch F \neq 2$, and we are concerned
with even smaller Galois groups than were considered by Bogomolov and
Pop.  Of course in this more general setting we cannot obtain as
precise results as Bogomolov and Pop, but we do get some interesting
information about the additive properties of multiplicative subgroups
of fields.  It would be very interesting to investigate further
relationships between our work and the quoted work of Bogomolov, Pop
and Tschinkel.


\Refs

\refstyle{A}
\widestnumber\key{MiSm3}

\ref \key AKMi \by A. Adem, D. Karagueuzian and J. Min\'a\v c \paper On the
cohomology of Galois groups determined by Witt rings \jour Adv. Math. \vol
148 \yr1999 \pages 105--160
\endref

\ref \key ArTa \by E. Artin and J. Tate \book Class field theory
\publ Benjamin \yr1967
\endref

\ref \key AEJ \by J. Arason, R. Elman and B. Jacob  \paper Rigid elements,
valuations, and realization of Witt rings  \jour J. Algebra \vol 110 \yr1987
\pages 449--467 \endref

\ref \key ArSch1 \by E. Artin and O. Schreier \paper Algebraische
Konstruktion reeller K\"orper \jour Abh. Math. Sem.
Univ. Hamburg \vol 5 \yr1927 \pages 85--99 \endref

\ref \key ArSch2 \by E. Artin and O. Schreier \paper Eine Kennzeichnung der
reell abgeschlossenen K\"orper \jour Abh. Math.
Sem. Univ. Hamburg \vol 5 \yr1927 \pages 225--231 \endref

\ref \key Be1 \by E. Becker \paper Euklidische K\"orper und euklidische H\"ullen
von K\"orpern  \jour J. reine angew. Math. \vol 268/269 \yr 1974 \pages 41--52
\endref

\ref \key Be2 \bysame \paper On the real spectrum of a ring and its
applications to semialgebraic geometry \jour Bull.
AMS \vol 15 \yr 1986 \pages 19--60 \endref

\ref \key BK \by E. Becker, E. K\"{o}pping \paper Reduzierte
quadratische Formen  und Semiordnungen reeller K\"{o}rper\jour Abh. Math.
Sem. Univ. Hamburg \vol 46 \yr 1977\pages 143--177
\endref

\ref \key BCR \by J. Bochnak, M. Coste, M. -F. Roy \book Real
Algebraic Geometry \bookinfo Ergeb. Math. \vol 36 \publ Springer-Verlag
\publaddr Berlin, Heidelberg \yr 1998 \endref


\ref \key Bo \by F. A. Bogomolov \paper On two conjectures in birational
algebraic geometry \inbook Algebraic Geometry and Analytic Geometry \bookinfo
ICM-90 Satellige Conference Proceedings \eds A. Fujiki et. al. \publ
Springer-Verlag \publaddr Tokyo \yr 1991
\endref

\ref \key BoT \by F. A. Bogomolov and Y. Tschinkel \paper Commuting elements in
Galois groups of function fields \paperinfo preprint in Ar. Xiv: Math.
AG/0012245, http://xxx.lanl.gov/  \yr 2000
\endref

\ref \key BS \by Z. I. Borevich and I. R. Shafarevich \book Number Theory
\publ Academic Press, Inc. \yr 1966
\endref

\ref  \key BEK \by Bredikhin, Ershov, and Kal'nei \paper Fields with two linear
orderings \jour Math Notes, Institute of Mathematics, Siberian Branch, Academy
of Sciences of hte USSR \yr 1970 \pages 319--325 \paperinfo Translated from
Mat. Zametki Vol. 7 (1970), 525--536
\endref

\ref \key Br \by L. Br\"ocker \paper Zur Theorie der quadratischen Formen
\"uber formal reellen K\"orpern  \jour Math. Ann.
\vol210 \yr1974 \pages233--256  \endref


\ref \key Cr1 \by T. Craven \paper Characterizing reduced Witt rings of
fields \jour J. Algebra \vol 53 \yr1978 \pages 68--77
\endref

\ref \key Cr2 \bysame \paper Fields maximal with respect to a set of orderings \jour J. Algebra \vol 115 \yr 1988
\pages 200--218 \endref

\ref \key CrSm \by T. Craven and T. Smith \paper Formally real fields from a
Galois theoretic perspective  \jour J. Pure Appl. Alg. \vol 145 \yr 2000
\pages 19--36
\endref

\ref \key Ef \by I. Efrat \paper Construction of valuations from K - theory
\jour Mathematical
Research Letters \vol 6 \yr 1999\pages  335--343
\endref

\ref \key ELP \by R. Elman, T.-Y. Lam and A. Prestel \paper On some
Hasse Principles over formally real fields \jour   Math. Z.
\yr 1973 \vol 134 \pages 291--301
\endref

\ref \key{End} \by O. Endler \book Valuation theory \publ Springer-Verlag
\publaddr New York \yr 1972
\endref

\ref\key EnKo \by  A. J. Engler and J. Koenigsmann
\paper Abelian subgroups of pro-p Galois groups
\jour Trans. A.M.S. \vol 350 \yr 1998  \pages 2473--2485
\endref

\ref \key Fr \by A. Fr\"{o}hlich \paper Orthogonal representations of
Galois groups, Stiefel-Whitney classes and Hasse-Witt
invariants \jour J. Reine Angew. Math. \vol 360 \yr 1985 \pages 84--123 \endref


\ref \key GM \by W. Gao and J. Min\'a\v c \paper Milnor conjecture and
Galois theory I \jour Fields Institute Communications \vol 16 \yr1997
\pages 95--110 \endref

\ref \key Gr \by G. Gr\"{a}tzer \book Universal Algebra \bookinfo
Second Edition \publ Springer-Verlag \yr 1979 \endref


\ref \key Hal \by M. Hall,Jr \book The theory of groups \publ The Macmillan
Company
\publaddr New York \yr 1959
\endref

\ref \key Har \by D. Harbater \paper Galois coverings of the
arithmetic line \jour Number Theory, New York 1984-1985, Lecture Notes
in Mathematics \vol 1240 \publ Springer-Verlag \publaddr New York
\yr 1987 \pages 165--195 \endref

\ref \key HaV\"ol \by D. Haran and H. V\"olklein \paper Galois groups
over complete valued fields \jour Israel Journal of Mathematics \vol
93 \yr 1996 \pages 9--27 \endref

\ref \key Jo \by D. L. Johnson \paper Non-nilpotent elements of
cohomology rings \jour Math.  Zeit.  \vol 112 \yr 1969 \pages 364--374
\endref

\ref \key Ki \by I. Kiming \paper Explicit classifications of some
$2$-extensions of a field of characteristic different
from $2$ \jour Can. J. of Math. \vol 42 \yr 1990 \pages 825--855 \endref

\ref \key Koc \by H. Koch \book Galoische Theorie der
$p$-Erweiterungen \publ Springer-Verlag
\publaddr New York \yr 1970
\endref

\ref \key K1 \by J. Koenigsmann \book Half-ordered fields
\bookinfo Dissertation Thesis \publ Universit\"at Konstanz
\publaddr Konstanz \yr 1993\endref

\ref \key K2 \bysame \paper From $p$-rigid elements to valuations (with a
Galois
characterization of $p$-adic fields)
\paperinfo with appendix by F. Pop
\jour J. Reine Angew. Math. \vol 465 \yr 1995 \pages 165--182
\endref


\ref  \key L1  \by T. Y. Lam
\book The algebraic theory of quadratic forms \publ Benjamin/Cummings
Publishing Co. \publaddr Reading, Mass. \yr 1980  \endref

\ref  \key L2  \bysame
\book Orderings, valuations and quadratic forms  \bookinfo Conference
Board of the Mathematical Sciences No. 52  \publ Amer. Math. Soc.
  \publaddr Providence, RI \yr1983  \endref

\ref \key L3 \bysame \paper An introduction to real algebra \jour Rocky
Mountain J. Math.
\vol 14 \yr 1984 \pages 767--814 \endref

\ref \key LaSm \by T. Y. Lam and T. Smith \paper On the
Clifford-Littlewood-Eckmann groups:  A new look at periodicity mod 8
\jour Rocky Mtn. J. Math. \vol 19 \yr 1989 \pages 749--786 \endref

\ref \key Mac \by S. MacLane
\paper Subfields and automorphism groups of $p$-adic fields \jour Ann. of
Math. \vol 40 \yr 1939 \pages 423--442
\endref

\ref \key Ma  \by M. Marshall  \book Abstract Witt rings  \bookinfo
Queen's Papers in Pure and Appl. Math., Vol. 57\publ Queen's University
\publaddr Kingston, Ontario \yr 1980 \endref


\ref\key Me \by A. Merkurjev \paper On the norm residue symbol of
degree two \jour Dokladi Akad. Nauk. SSSR \yr 1981 \vol 261 \pages
542-547 \paperinfo English translation:  Soviet Math. Dokladi
   {\bf 24} (1981), 546-551\endref


\ref\key Mi1 \by J. Min\'a\v c \paper Elementary 2-abelian subgroups of
W-groups,
\paperinfo unpublished manuscript \yr 1987 \endref

\ref \key MiSm1 \by J. Min\'a\v c and T. Smith \paper A characterization of
C-fields via Galois groups \jour J. Algebra  \yr 1991 \vol 137 \pages 1--11
\endref


\ref \key MiSm2 \bysame \paper Decomposition of
Witt rings and Galois groups \jour Canad. J. Math\vol47 \yr1995
\pages1274--1289 \endref

\ref \key MiSm3 \bysame \paper W-Groups under quadratic extensions of
fields \jour Canad. J. Math. \vol 52 \yr 2000
\pages 833--848 \endref

\ref \key MiSm4 \bysame \paper Formally dyadic fields \paperinfo in preparation
\endref

\ref \key MiSp1 \by J. Min\'a\v c and M. Spira \paper Formally real fields,
pythagorean fields, C-fields and W-groups \jour Math. Z. \yr1990 \vol205
\pages519--530 \endref

\ref \key MiSp2 \bysame
\paper Witt rings and Galois groups \jour Ann. of Math. \yr 1996 \vol 144
\pages 35--60 \endref

\ref \key Mo \by S. Morris  \book Pontryagin duality and the structure
of locally compact abelian groups  \bookinfo London Mathematical
Society Lecture Note Series
\vol 29 \publ Cambridge University Press
\publaddr Cambridge \yr1980 \endref


\ref\key Pf \by A. Pfister \paper Quadratische Formen in beliebigen K\"orpern
\jour Invent. Math. \vol 1 \pages 116--132 \yr 1966
\endref

\ref \key Po1 \by F. Pop \paper On Grothendieck's conjecture of birational
anabelian geometry \jour Ann. Math. \vol 139 \yr 1994 \pages 145--182
\endref

\ref \key Po2 \bysame \paper Glimpses of Grothendieck's anabelian geometry
\inbook Geometric Galois Action 1: Around Grothendieck's Esquisse d'un
Programme
\bookinfo Lecture Notes in Mathematics \vol 242 \eds L. Schneps and P. Lochak
\publ Springer-Verlag \publaddr Berlin \yr 1997 \pages 113--126
\endref


\ref \key Ri \by P. Ribenboim \book Th\acuteaccent eorie des Valuations
\bookinfo Universit\acuteaccent e de Montreal \publ Les Presses de
l'Universit\acuteaccent e Montreal \yr1965 \endref

\ref \key RZ \by L. Ribes and P. Zalesskii\book Profinite Groups
\publ Springer-Verlag \publaddr Berlin \yr 2000  \endref

\ref \key Sc \by W. Scharlau \book Quadratic and Hermitian Forms \publ
Springer-Verlag \yr 1985
\endref

\ref \key Sch \by N. Schwartz \paper Signatures and Real Closures of
Fields \jour Publications de l'Universit\'e Paris VII \vol 32 \yr 1990
\pages 65--78
\endref

\ref \key S \by E. Sperner \paper Die Ordnungsfunktionen einer Geometrie
\jour Math. Ann. \vol 121 \yr 1949 \pages 107--130
\endref

\ref \key Sp \by M. Spira
\book Witt rings and Galois groups
\bookinfo Ph.D. Thesis \publ University of California \publaddr Berkeley,
California \yr 1987
\endref

\ref \key Vo \by V. Voevodsky \paperinfo On $2$-torsion in motivic
cohomology, http:/www.math.uiuc.edu/K-theory/0502/index.
html \yr 2001 \endref

\ref \key Wa \by R. Ware \paper Valuation rings and rigid elements in fields
\jour Canad. J. Math. \vol 33 \yr1981 \pages 1338--1355 \endref

\endRefs
\enddocument